\documentclass[12pt,a4paper,reqno]{amsart}
\usepackage{tikz}
\usetikzlibrary{arrows.meta,positioning,fit}
\usepackage{amssymb,amsmath,amsfonts,amsthm}
\usepackage{yhmath}
\usepackage{url}
\usepackage{preamble}
\usepackage{tikz-cd}
\usepackage{color}
\usepackage{todonotes}
\usepackage{mathrsfs}
  
\allowdisplaybreaks[2]
\let\le\leqslant
\let\leq\leqslant
\let\ge\geqslant
\let\geq\geqslant

\declaretheorem[sibling=thm]{corollary}
\numberwithin{equation}{section}

\newtheorem{mainthm}{Theorem}

\newcommand{\Diff}{\operatorname{Diff}}
\newcommand{\wh}{\widehat}
\newcommand{\eps}{\varepsilon}
\newcommand{\dd}{\,d}
\newcommand{\dist}{\operatorname{dist}}
\newcommand{\vol}{\operatorname{vol}}
\newcommand{\conorm}{\operatorname{m}}
\newcommand{\Ent}{\operatorname H}
\newcommand{\cP}{\mathcal P}
\newcommand{\cA}{\mathcal A}
\newcommand{\cQ}{\mathcal Q}
\newcommand{\TV}{\operatorname{TV}}
\newcommand{\id}{\operatorname{id}}
\newcommand{\Se}{\mathcal{S}}

\crefname{mainthm}{Theorem}{Theorems}
\Crefname{mainthm}{Theorem}{Theorems}
\crefname{thm}{Theorem}{Theorems}
\Crefname{thm}{Theorem}{Theorems}
\crefname{theorem}{Theorem}{Theorems}
\Crefname{theorem}{Theorem}{Theorems}
\crefname{lemma}{Lemma}{Lemmas}
\Crefname{lemma}{Lemma}{Lemmas}
\crefname{proposition}{Proposition}{Propositions}
\Crefname{proposition}{Proposition}{Propositions}
\crefname{claim}{Claim}{Claims}
\Crefname{claim}{Claim}{Claims}
\crefname{corollary}{Corollary}{Corollaries}
\Crefname{corollary}{Corollary}{Corollaries}
\crefname{definition}{Definition}{Definitions}
\Crefname{definition}{Definition}{Definitions}
\crefname{remark}{Remark}{Remarks}
\Crefname{remark}{Remark}{Remarks}
\crefname{section}{Section}{Sections}
\Crefname{section}{Section}{Sections}
\crefname{subsection}{Section}{Sections}
\Crefname{subsection}{Section}{Sections}
\crefname{subsubsection}{Section}{Sections}
\Crefname{subsubsection}{Section}{Sections}
\crefname{equation}{Equation}{Equations}
\Crefname{equation}{Equation}{Equations}

\title[Entropic continuity of Lyapunov exponents]
{Entropic continuity of Lyapunov exponents\\ in arbitrary dimension}
\begin{document}
\begin{abstract}

We prove that the sum of the positive Lyapunov exponents is entropically continuous for every $C^\infty$ diffeomorphism, provided the entropy exceeds its second-largest exterior growth rate. We also establish stable strong positive recurrence and show that the topological entropy is continuous, and the number of ergodic measures of maximal entropy is upper semicontinuous. This extends the result of Buzzi, Crovisier, and Sarig~\cite{BCSLyapunov} to arbitrary dimensions. 
We give examples in both the measure-theoretic and topological senses, showing that this entropy gap condition is sharp. We also present applications related to the Newhouse phenomenon and standard maps.
\end{abstract}
\maketitle
\begingroup
\small
\tableofcontents
\endgroup
\medskip

\section{Introduction}

\subsection{Entropy and Lyapunov exponents}

Metric entropy measures orbit complexity, whereas Lyapunov exponents
measure the growth of tangent vectors.  Their relation is expressed by the
Margulis--Ruelle inequality and, for regular invariant measures, by
Pesin's entropy formula.  The question studied here is a continuity
version of this relation: when both an invariant measure and its entropy
converge, must its integrated positive Lyapunov sum converge as well?

For a smooth diffeomorphism $f$ and an invariant probability $\mu$, write
$h_\mu(f)$ for metric entropy and define
\begin{equation}\label{eq:def-Splus}
 S^+(\mu,f):=\int_M\sum_{\chi_i(x,f)>0}\chi_i(x,f)\,\dd\mu(x),
\end{equation}
where the Lyapunov exponents are counted with multiplicity.  They exist
at $\mu$-almost every point by Oseledets' theorem~\cite{Oseledets}. This integrated sum encodes the total infinitesimal expansion experienced by $\mu$-typical trajectories. 
Margulis--Ruelle's inequality gives $h_\mu(f)\le S^+(\mu,f)$
\cite{RuelleInequality}.  For absolutely continuous invariant measures,
Pesin's entropy formula gives equality
\cite{Pesin-entropyformula_1977,Ledrappier_Strelcyn_1982}; generally, Ledrappier and Young characterize equality by absolute continuity
of the conditional measures on unstable manifolds
\cite{LEDRAPPIER_YOUNG_A}.  Without this regularity, equality need not hold,
so entropy convergence alone does not give an immediate conclusion about
Lyapunov exponents.

The two invariants have different continuity properties.  For entropy,
Newhouse and Yomdin established upper semicontinuity in the smooth setting
\cite{Newhouse,yomdin1987}; finite-smoothness obstructions appear in
\cite{mis1973,New05,Buzzi2014}.  For Lyapunov exponents, continuity results
in random and hyperbolic settings include
\cite{Furstenberg1983,RuelleRandomMatrices1979,avilaviana2010,bockeviana2017,AvilaEskinViana2023}.
In the quasi-periodic setting, analytic continuity results
\cite{GoldsteinSchlag2001,BourgainJitomirskaya2002} contrast with smooth
discontinuity examples~\cite{WangYou2013}; see also
\cite{viana2020,Bochi-Viana} for other continuity and discontinuity phenomena.

Buzzi, Crovisier, and Sarig~\cite{BCSLyapunov} proved entropic continuity
for $C^\infty$ surface diffeomorphisms: on the space of ergodic invariant
measures with positive entropy, convergence of the measure and its entropy
forces convergence of the Lyapunov exponents.  Their later work on strong
positive recurrence (SPR)~\cite{SPR2025} links this continuity to uniform
control of high-entropy measures on Pesin blocks.

We establish a higher-dimensional criterion for these conclusions.  The
hypothesis compares entropy with the second largest exterior growth rate
of the derivative.  It singles out the unstable dimension relevant to
high entropy.  At a fixed unstable dimension $k$, the proof uses only the
$(k-1)$-st exterior growth rate; the other lower exterior powers enter a
local counting error whose average vanishes along large iterates.  The
same argument yields stable SPR and stability of measures of maximal
entropy (MMEs), which attain the supremum defining $h_{\mathrm{top}}(f)$
by the variational principle.

\subsection{Main results}
Let $M$ be a closed Riemannian manifold of dimension $d\ge 2$. Let $\Diff^r(M)$ ($1\le r\le \infty$) denote the space of $C^r$ diffeomorphisms of $M$, endowed with the $C^r$ topology, and let $\mathcal P(M)$ be the space of Borel probability measures on $M$ equipped with the weak$^*$ topology. For $f\in\Diff^1(M)$, let $\mathcal M(f)$ and $\mathcal M_e(f)$ denote, respectively, the sets of $f$-invariant and ergodic $f$-invariant Borel probability measures.

\bigskip

\paragraph{\it Exterior growth rates.} For
$f\in\Diff^1(M)$ and $1\le q\le d-1$, define the
$q$-th exterior growth rate by
\begin{equation}\label{eq:def-lambda-fi}
 \lambda_q(f):=
 \lim_{n\to\infty}\frac1n
 \log\|\wedge^qdf^n\|_\infty,
\end{equation}
where the limit exists by subadditivity.  We also set $\lambda_0(f):=0$.
The change of variables formula and the ordering of singular values give
$\|\wedge^qdf^n\|_\infty\ge1$, and hence $\lambda_q(f)\ge0$, for
$1\le q\le d-1$.  We call the second largest value among
$\lambda_0(f),\ldots,\lambda_{d-1}(f)$ the
\emph{second-largest exterior growth rate} of $f$: 
\begin{equation}\label{eq:def-lambda-sec}
 \lambda_{\mathrm{sec}}(f):=
 \max_{0\le q<q'\le d-1}
   \min\{\lambda_q(f),\lambda_{q'}(f)\}.
\end{equation}
With this convention, $\lambda_{\mathrm{sec}}(f)=0$ on surfaces.

 If $\mu\in\mathcal M(f)$ and $\tau_\mu$ is its
ergodic-decomposition probability on $\mathcal M_e(f)$, we say that
$\mu$ has a \emph{componentwise entropy gap} if
\begin{equation}\label{eq:def-componentwise-gap}
 h_\rho(f)>\lambda_{\mathrm{sec}}(f)
 \quad\text{for }\tau_\mu\text{-almost every }\rho\in\mathcal M_e(f).
\end{equation}
For an ergodic measure this is exactly 
$h_\mu(f)>\lambda_{\mathrm{sec}}(f)$.  If $\mu$ is a measure of maximal
entropy and $h_{\mathrm{top}}(f)>\lambda_{\mathrm{sec}}(f)$, then
\eqref{eq:def-componentwise-gap} naturally  holds because almost every ergodic
component of $\mu$ is maximal.

\begin{mainthm}[Measure-theoretic entropic continuity]
\label{thm:general-entropic-continuity}
Let $f_n,f\in\Diff^\infty(M)$ with $f_n\to f$ in $C^\infty$, and let
$\nu_n\in\mathcal M(f_n)$ and $\mu\in\mathcal M(f)$.  Assume
\begin{equation}\label{eq:general-core-assumptions}
 \nu_n\longrightarrow\mu,
 \qquad
 h_{\nu_n}(f_n)\longrightarrow h_\mu(f),
\end{equation}
 and
\begin{equation}\label{eq:general-componentwise-gap}
 h_\rho(f)>\lambda_{\mathrm{sec}}(f)
 \quad\text{for }\tau_\mu\text{-almost every }\rho\in\mathcal M_e(f).
\end{equation}
Then
\begin{equation}\label{eq:general-Splus}
 S^+(\nu_n,f_n)\longrightarrow S^+(\mu,f).
\end{equation}
\end{mainthm}

We use strong positive recurrence in the sense of
\cite[Definition~1.3]{SPR2025}.
\begin{definition}
A diffeomorphism is strongly positively recurrent if there is $\chi>0$ such
that, for every $\eps>0$, one can find a $(\chi,\eps)$-Pesin block
$\Lambda$, a threshold $h_0<h_{\mathrm{top}}(f)$, and $\tau>0$ for which
\[
 h_\nu(f)>h_0
 \quad\Longrightarrow\quad
 \nu(\Lambda)>\tau
\]
for every ergodic invariant probability $\nu$.

\end{definition} 
The entropy gap singles out one exterior degree.  
Our main result shows that
this spectral separation is strong enough to imply stable strong positive
recurrence and, at the same time, to make the maximal entropy part of the
Lyapunov spectrum continuous.  In the non-invertible setting, Delplanque and Li~\cite{DelplanqueLi2026}
introduced the expanding-SPR property, which ensures that
all high-entropy ergodic measures have full unstable index, and established entropic continuity for such measures.

\begin{definition} A  diffeomorphism $f$ is called  entropy-dominant  if 
\[
 h_{\mathrm{top}}(f)>\lambda_{\mathrm{sec}}(f).
\]
\end{definition} 

This property is $C^1$-open in the space  $\Diff^{1+}(M)$ of  $C^1$ diffeomorphisms with H\"older continuous derivative (see Corollary \ref{cor:openess} for a proof). For $1\le r\le \infty$, we denote 
$$\Se^{r}(M)=\{f\in  \Diff^r(M): \,\,f \, \,\text{is   entropy-dominant} \}.$$ 
Note that for diffeomorphisms on $M$, the notation $f_n\to f$ in $C^\infty$ means
that $f_n^{\pm1}\to f^{\pm1}$ in $C^s$ for every finite integer $s\ge1$.

\begin{mainthm}[Stable SPR]\label{thm:main-intro}  Every
continuity or upper semicontinuity assertion below is understood with respect
to the $C^\infty$ topology.
\begin{enumerate}
\item every $g\in \Se^{\infty}(M)$ is strongly positive recurrent;
\item the map $g\mapsto h_{\mathrm{top}}(g)$ is continuous on $\Se^{\infty}(M)$;
\item every $g\in \Se^{\infty}(M)$ has finitely many ergodic measures of maximal entropy,
      and
      \[
       n_{\mathrm{MME}}(g)
       :=\#\{\mu\in\mathcal M_e(g):h_\mu(g)=h_{\mathrm{top}}(g)\}
      \]
      is upper semicontinuous on $\Se^{\infty}(M)$;
\item if $g_n,g\in \Se^{\infty}(M)$, $g_n\to g$ in $C^\infty$,
      $\mu_n\in\mathcal M(g_n)$, $\mu_n\to\mu$, and
      \[
       h_{\mu_n}(g_n)\longrightarrow h_{\mathrm{top}}(g),
      \]
      then $\mu$ is a measure of maximal entropy for $g$ and
      \[
       S^+(\mu_n,g_n)\longrightarrow S^+(\mu,g).
      \]
      Consequently, with
      \[
       X:=\{(\mu,g):g\in\Se^{\infty}(M),\ \mu\in\mathcal M(g),\
                    h_\mu(g)=h_{\mathrm{top}}(g)\},
      \]
      the $X$ is closed in $\mathcal P(M)\times \Se^{\infty}(M)$, and
      $(\mu,g)\mapsto S^+(\mu,g)$ is continuous on $X$;
\item if $f\in\Se^{\infty}(M)$ has a unique measure of maximal entropy $\mu_f$, then  every  $g\in\Se^{\infty}(M)$ in a neighborhood of $f$  has a unique measure of maximal
      entropy $\mu_g$, and the maps
      \[
       g\longmapsto\mu_g,\qquad
       g\longmapsto h_{\mu_g}(g),\qquad
       g\longmapsto S^+(\mu_g,g)
      \]
      are continuous on this neighborhood.
\end{enumerate}
\end{mainthm}

\subsection{Sharpness of the entropy gap}
\label{subsec:intro-sharpness}
For an ergodic invariant probability $\nu$, write
$\operatorname{Lyap}(\nu,f)$ for its Lyapunov spectrum, counted with
multiplicity.  Set
\[
 A=\begin{pmatrix}2&1\\1&1\end{pmatrix}:\mathbb T^2\to\mathbb T^2,
 \qquad
 h_A=\log\frac{3+\sqrt5}{2}.
\]
The following two sharpness examples show that the strict entropy inequalities in
the main theorems cannot be weakened to
nonstrict ones.

\begin{proposition}[Sharpness at the measure-theoretic threshold]
\label{prop:sharp-measure}
There exist positive numbers $H>\gamma>0$, a fixed volume-preserving
$F\in\Diff^\infty(\mathbb T^2\times S^2)$, ergodic invariant probabilities
$\mu_n,\mu$, and numbers $\gamma_n\ge\gamma/4$ with
$\gamma_n\to\gamma/2$ such that
\[
 \mu_n\longrightarrow\mu,
 \qquad
 h_{\mu_n}(F)=h_\mu(F)=\lambda_{\mathrm{sec}}(F)=H,
\]
and
\[
 \operatorname{Lyap}(\mu_n,F)=\{H,\gamma_n,-\gamma_n,-H\},
 \qquad
 \operatorname{Lyap}(\mu,F)=\{H,\gamma,-\gamma,-H\}.
\]
Thus all these measures are hyperbolic of unstable dimension two, every
Lyapunov exponent has absolute value at least
$\min\{H,\gamma/4\}>0$, and
\[
 S^+(\mu_n,F)\longrightarrow H+\frac\gamma2
 <H+\gamma=S^+(\mu,F).
\]
\end{proposition}

Since the limiting measure is ergodic, the componentwise hypothesis of
\cref{thm:general-entropic-continuity} reduces here to
$h_\mu(F)>\lambda_{\mathrm{sec}}(F)$; the proposition shows that equality
cannot be admitted, even with a fixed map, a common index, and a uniform
Lyapunov gap.  The construction uses a conservative quadratic return and a
sufficiently high power of the cat map.

\begin{proposition}[Sharpness at the topological-entropy threshold]
\label{prop:sharp-main}
For every $c\in(0,h_A)$, there exist a fixed volume-preserving
$F\in\Diff^\infty(\mathbb T^4)$, ergodic measures of maximal entropy
$\mu_n,\mu$, and numbers $\kappa_n>0$ with $\kappa_n\to0$ such that
\[
 h_{\mathrm{top}}(F)=\lambda_{\mathrm{sec}}(F)=h_A,
 \qquad
 \mu_n\longrightarrow\mu,
 \qquad
 h_{\mu_n}(F)=h_\mu(F)=h_A,
\]
and
\begin{align*}
 \operatorname{Lyap}(\mu_n,F)&=\{h_A,\kappa_n,-\kappa_n,-h_A\},\\
 \operatorname{Lyap}(\mu,F)&=\{h_A,c,-c,-h_A\}.
\end{align*}
In particular, all these measures are hyperbolic of unstable dimension two,
but
\[
 S^+(\mu_n,F)=h_A+\kappa_n\longrightarrow h_A
 \ne h_A+c=S^+(\mu,F).
\]
\end{proposition}

Thus the strict inequality in \cref{thm:main-intro} is necessary already for
Item~\textnormal{(4)}, even for a fixed volume-preserving map and hyperbolic
ergodic maximal measures of a common index.  The construction couples the cat
map with a zero-entropy Hamiltonian surface map whose resonant cycles are
hyperbolized inside annuli accumulating on a saddle.

\subsection{Examples and applications}
The two sharpness propositions are proved in
\cref{subsec:sharp-main-proof,subsec:sharp-measure-proof}.  We then apply the
entropy-gap criterion to Newhouse domains and products of standard maps.
Both applications concern high-entropy dynamics without requiring uniform
hyperbolicity of the ambient map.

\subsubsection{Newhouse domains}

A \emph{Newhouse domain} is a $C^\infty$-open set in which tangencies
between stable and unstable manifolds of hyperbolic basic sets persist.  The
next application places such a domain entirely inside the entropy-gap region
of \cref{thm:main-intro}.

Let $L\in\mathrm{GL}(3,\mathbb Z)$ be hyperbolic and let
$0<r_1\le r_2\le r_3$ be the moduli of its eigenvalues.  Put
\begin{equation}\label{eq:intro-newhouse-chi}
 \chi_L:=\min_{\lambda\in\operatorname{Spec}(L)}
              \bigl|\log|\lambda|\bigr|
          =|\log r_2|.
\end{equation}
For $\Delta>0$, we say $L$ is \emph{$\Delta$-hyperbolic} if
$\chi_L>\Delta$.

\begin{proposition}[Newhouse domains in $\Delta$-hyperbolic classes]
\label{prop:newhouse-domains}
There exists a universal constant $\Delta_{\mathrm{NH}}>0$ with the
following property.  Let $L\in\mathrm{GL}(3,\mathbb Z)$ be
$\Delta_{\mathrm{NH}}$-hyperbolic and put
$k=\dim E_L^u\in\{1,2\}$.  Then the homotopy class of $L$ contains a
nonempty $C^\infty$-open set
\[
 \mathcal N_L\subset\Diff^\infty(\mathbb T^3)
\]
such that: 
\begin{enumerate}[label=\textnormal{(\roman*)},leftmargin=10mm]
\item every $f\in\mathcal N_L$ is homotopic to $L$ and satisfies
      \[
       h_{\mathrm{top}}(f)>\lambda_{\mathrm{sec}}(f);
      \]
\item $\mathcal N_L$ is a Newhouse domain.  If $k=1$, maps with infinitely
      many nonhyperbolic H\'enon-like strange attractors are dense in
      $\mathcal N_L$; if $k=2$, the analogous statement holds for strange
      repellers;
\item every measure of maximal entropy of every
      $f\in\mathcal N_L$ is hyperbolic with unstable dimension $k$;
\item all conclusions of \cref{thm:main-intro} hold in
      $\mathcal N_L$.
\end{enumerate}
\end{proposition}

The condition measures the distance of the middle spectral modulus from the
unit circle.  Its role is to absorb the universal $C^1$ cost of compressing
one multiplier and bending one transverse homoclinic intersection, while
Manning's homological lower bound supplies the entropy.  Since
$\chi_{L^m}=m\chi_L$, the proposition applies to every sufficiently large
power of every hyperbolic $L\in\mathrm{GL}(3,\mathbb Z)$.

\medskip

\subsubsection{Products of standard maps}
\label{subsec:intro-standard-products}

On $\mathbb T^2=(\mathbb R/\mathbb Z)^2$, consider the standard family:
\begin{equation}\label{eq:standard-map-definition}
 f_k(x,y)=(2x-y+k\sin(2\pi x),x)\pmod1,\qquad k>0.
\end{equation}
These maps preserve area.  Duarte~\cite{Duarte1994} constructed basic sets
with entropy growing like $\log k$ for all large parameters, and showed
that for a residual set of large parameters these sets are accumulated by
elliptic islands.  Obata~\cite[Theorem~A]{obatastandardmap} proved that
$f_k$ nevertheless has a unique Bernoulli measure of maximal entropy for
every sufficiently large $k$.

The next result places arbitrary finite products of large standard maps
inside the entropy-gap region, with no comparability requirement on the
parameters.

\begin{theorem}[Products of standard maps]
\label{thm:standard-map-products}
For every $\alpha>0$ and integer $m\ge1$, there exists
$k_0=k_0(\alpha,m)>1$ such that every product
$F=f_{k_1}\times\cdots\times f_{k_m}$ with $k_i>k_0$ for $1\le i\le m$
satisfies
\[
 h_{\mathrm{top}}(F)-\lambda_{\mathrm{sec}}(F)>\alpha.
\]
In particular, $F\in\Diff^\infty(\mathbb T^{2m})$ satisfies the hypothesis
of \cref{thm:main-intro}.
\end{theorem}

Consequently, stable SPR, continuity of topological entropy, upper
semicontinuity of the number of ergodic MMEs, and entropic continuity of
$S^+$ near maximal entropy hold in $C^\infty$ neighborhoods in the full space
$\Diff^\infty(\mathbb T^{2m})$, not only within the product family.
Combining this gap with Obata's Bernoulli property gives uniqueness and
stability of the maximal measure, even outside the product family.
Both results are proved in \cref{subsec:standard-products-proof}.

\begin{corollary}[Stable uniqueness for products of standard maps]
\label{cor:standard-product-uniqueness}
For every integer $m\ge1$, there exists $K_m>1$ such that whenever
$k_i>K_m$ for all $i$, the product $F=f_{k_1}\times\cdots\times f_{k_m}$
has a unique measure of maximal entropy
\[
 \mu_F=\mu_{k_1}\otimes\cdots\otimes\mu_{k_m},
\]
which is Bernoulli.  Here $\mu_{k_i}$ is the unique maximal measure of
$f_{k_i}$.  Moreover, $F$ has a $C^\infty$ neighborhood
$U_F\subset\Diff^\infty(\mathbb T^{2m})$ on which every $g$ has a unique
maximal measure $\mu_g$, and the maps
\[
 g\longmapsto\mu_g,\qquad
 g\longmapsto h_{\mu_g}(g),\qquad
 g\longmapsto S^+(\mu_g,g)
\]
are continuous.
\end{corollary}

\bigskip
\noindent{\bf Strategy and organization.}\,\,
We extend to higher dimensions the entropic continuity of Lyapunov exponents established for smooth surface diffeomorphisms in \cite{BCSLyapunov, Burguet}.
To this end we use a slightly modified version of the Yomdin
reparametrizations built in \cite{BOB}. In their work, the dynamics
of a $k$-disk $D$ is described by a tree whose nodes at level $n$
represent elementary pieces of $f^nD$. The red and blue nodes
correspond, respectively, to $k$-dimensional growth and
$l$-dimensional growth with $l<k$. 
 
 The  argument is carried out at a fixed unstable index $k$.
Sections~\ref{sec:preliminaries}--\ref{sec:plaque-selection} combine the
closed red--blue reparametrization tree with record times and entropy on
unstable plaques.  The resulting decomposition separates an invariant
hyperbolic component from a neutral remainder.  Its counting estimate keeps
the lower-dimensional exterior growth local: a continuous nonnegative
error records the excess over a prescribed threshold.  Along large
iterates, this error has vanishing average for a measure of unstable index
$k$ as soon as the threshold exceeds $\lambda_{k-1}(f)$.

In \cref{sec:ergodic-continuity}, a suitable choice of cyclic component and
backward contraction identify the limiting unstable lift.  An
entropy--Lyapunov graph-measure argument then treats arbitrary invariant
probabilities.  This proves the fixed-index criterion
\cref{thm:fixed-index-continuity}; the spectral index-selection lemma gives
\cref{thm:general-entropic-continuity} as a consequence.
Section~\ref{sec:main-theorem-proof} proves a fixed-index criterion for
stable SPR and maximal entropy stability, and then derives
\cref{thm:main-intro}.  Section~\ref{sec:examples} proves the sharpness
examples, the Newhouse domain application, and the standard map products.

\subsection*{Acknowledgment}
GL was partially supported by the National Key R\&D Program of China (2022YFA1005802), NSFC(12671228).
YT was partially supported by the National Key R\&D Program of China (2022YFA1005801).

\section{Reparametrization and spectral preliminaries}
\label{sec:preliminaries}

We first describe the reparametrization tree and its pointwise ownership
rule.  We then formulate the fixed-index continuity theorem and explain how
the exterior growth gap selects that index.

\subsection{Bounded couples}

We use the conventions
\[
 \log^+t=\max\{\log t,0\},
 \qquad
 \wedge^0A=1.
\]
For an invertible linear map $A$ and a $k$-plane $E$, let
\[
 \conorm(A):=\|A^{-1}\|^{-1},
 \qquad
 J_k(g,x,E):=\|\wedge^k(dg_x|_E)\|.
\]
Let $g$ be a diffeomorphism on $M$. If $D\subset M$ is an embedded $k$-disk, we abbreviate
\[
 J_k(g^n|_{T_xD}):=J_k(g^n,x,T_xD).
\]
We write $\operatorname{Leb}_k$ for Lebesgue measure on a parameter cube and
$\vol_k$ for the Riemannian $k$-volume on an embedded disk.

We shall repeatedly use terminal source maps
$\psi:[0,1]^k\to M$ satisfying the orbit bound
\[
 \|d(g^j\circ\psi)\|_\infty\le K_0,
 \qquad 0\le j\le N,
\]
where the constant $K_0$ will be precisely supplied by
\cref{prop:yomdin-tree-restart}.



We now recall the terminology of Ben Ovadia and Burguet in \cite{BOB}.  A
$C^{r-1,1}$ map is a $C^{r-1}$ map whose derivatives of order $r-1$ are
Lipschitz.

\begin{definition}[Bounded couples and strongly $\eps$-bounded couples]
\label{def:strongly-bounded-front}
  For a $C^{r-1,1}$ map $u$, set
\[
 \|u\|_{\mathrm{BO},r}:=\max_{1\le s\le r}[u]_s   \,\,\, \text{with}\,\,\,  [u]_s=\max_{|\alpha|=s-1}\operatorname{Lip}(\partial^\alpha u).
\]
A \emph{couple} is a pair $(\sigma,\theta)$ consisting of a $C^r$ embedded
disk $\sigma:[0,1]^k\to M$ and a $C^{r-1,1}$ map
$\theta:[0,1]^k\to[0,1]^k$.  It is a \emph{bounded couple} if
\[
 \left[t\mapsto\wedge^k d_{\theta(t)}\sigma\right]_s
 \le \frac1{10k}\sup_t\left|\wedge^k d_{\theta(t)}\sigma\right|,
 \qquad 1\le s\le r-1.
\]
A bounded couple is \emph{strongly $\eps$-bounded} if
\[
 \|\sigma\circ\theta\|_{\mathrm{BO},r}
 \le\frac{\eps}{\sqrt{k}}.
\]
There exists a constant $C_M>0$ such that, for every strongly
$\eps$-bounded couple,
\begin{equation}\label{eq:strong-chart-volume-front}
 \vol_k\bigl((\sigma\circ\theta)([0,1]^k)\bigr)
 \le C_M k^{-k/2}\eps^k.
\end{equation}
\end{definition}

\subsection{The one-step construction}
\label{The one-step language}

The construction in \cite{BOB} organizes successive images of a disk
as a red--blue tree of reparametrized subdisks.  We recall only the  estimates used below.

Put $Q=[0,1]^k$.  For two maps
$\theta_+,\theta_-:Q\to Q$, write $\theta_+>\theta_-$ if
\[
 \theta_-=\theta_+\circ h
 \quad\text{for some }h:Q\to Q
 \quad\text{with}\quad \|dh\|_\infty<\tfrac12.
\] 
That is,  $\theta_-$ parametrizes a strictly contracted part of the image
parametrized by $\theta_+$.  For a couple $(\sigma,\theta)$, we distinguish
between the \emph{carrier disk} $\sigma(Q)$ and the \emph{constituent image}
$(\sigma\circ\theta)(Q)$.  Tangent spaces are always taken in the carrier;
this matters at points represented by the boundary of a parameter cube.

We also recall the geometric content of an admissible family \cite[Definition~3.6]{BOB}.  Let
$(\sigma,\theta')$ be a bounded couple.  A finite family $\Theta$ of maps
$\theta:Q\to Q$ is called a
$(\sigma,\theta',\eps)$-\emph{admissible family} if:
\begin{enumerate}[label=\textnormal{(A\arabic*)},leftmargin=11mm]
\item every $(\sigma,\theta)$, $\theta\in\Theta$, is strongly
      $\eps$-bounded;
\item the parameter images $\theta(Q)$ are disjoint up to their
      $k$-dimensional boundaries and their union is contained in
      $\theta'(Q)$;
\item writing
      \[
       \sigma(\Theta):=\bigcup_{\theta\in\Theta}
       (\sigma\circ\theta)(Q),
      \]
      there are $x\in(\sigma\circ\theta')(Q)$, a $k$-plane
      $H\subset T_xM$ tangent to the carrier, and a $1$-Lipschitz map
      $\psi:H\to H^\perp$ with $\psi(0)=0$ such that
      \[
       \Gamma_\psi(\eps)\subset\sigma(\Theta)
       \subset\Gamma_\psi(2\eps),
      \]
      where
      \[
       \Gamma_\psi(\delta)
       :=\exp_x\{v+\psi(v):v\in H,\ |v|\le\delta\}
      \]
      is the corresponding centered $\delta$-graph.  The scale $\eps$ is
      chosen so that these graphs lie in normal coordinate balls on which
      the exponential maps and their inverses have uniform Lipschitz bounds.
\end{enumerate}
\smallskip
\smallskip
\smallskip
\smallskip

\paragraph{\bf The Yomdin tree.}
A \emph{blue node} consists of one strongly $\eps$-bounded couple, and a
\emph{red node} consists of an admissible family.  A red
node carries a definite amount of $k$-volume, of order $\eps^k$, whereas a
blue node is a single residual piece whose offspring can be counted by
lower-dimensional exterior growth.  At the next step each constituent is
reparametrized separately; the relation $>$ records the nesting of the
parameter maps along a branch.

The raw reparametrizations in \cite{BOB} are disjoint on the interiors of
the parameter cubes and cover the parent image modulo a $k$-null set.  Our
entropy argument must also cover points on the artificial algebraic
boundaries, since an unstable conditional measure need not ignore them.  We
therefore assign every point to one child by a fixed Borel rule.  The next
lemma performs this passage from a modulo null refinement to an exact
pointwise refinement.  The resulting pieces will be called \emph{owned
sets}.

\begin{lemma}[Closed refinement with Borel ownership]
\label{lem:closed-ownership}
Let $Q=[0,1]^k$.  Let $\theta:Q\to\mathbb R^k$ be Lipschitz and injective
on $(0,1)^k$, and let
$\theta_1,\ldots,\theta_t:Q\to\theta(Q)$ be Lipschitz maps with the same
property.  Set
\[
 U:=\theta((0,1)^k),
 \qquad
 U_i:=\theta_i((0,1)^k),
 \qquad
 K_i:=\theta_i(Q).
\]
Assume that the open sets $U_i$ are pairwise disjoint and that
\[
 \operatorname{Leb}_k\!\left(U\setminus\bigcup_{i=1}^tU_i\right)=0.
\]
Then
\begin{equation}\label{eq:closed-refinement-cover}
 \theta(Q)=\bigcup_{i=1}^tK_i,
 \qquad
 \operatorname{Leb}_k(K_i\cap K_j)=0
 \quad(i\ne j).
\end{equation}

More generally, let $B\subset\theta(Q)$ be Borel and satisfy
$\operatorname{Leb}_k(\theta(Q)\setminus B)=0$.  After fixing an order on the
children, define
\begin{equation}\label{eq:closed-ownership-rule}
 B_1:=B\cap K_1,
 \qquad
 B_i:=B\cap\left(K_i\setminus\bigcup_{h<i}K_h\right)
 \quad(2\le i\le t).
\end{equation}
Then the sets $B_i$ are Borel and pairwise disjoint, and
\begin{equation}\label{eq:closed-ownership-conclusions}
 B=\bigsqcup_{i=1}^tB_i,
 \qquad
 B_i\subset K_i,
 \qquad
 \operatorname{Leb}_k(K_i\setminus B_i)=0.
\end{equation}
The same conclusions hold after composition with a $C^1$ embedded disk:
parameter-null sets and $k$-dimensional null sets in the embedded image
correspond under the parametrization.
\end{lemma}

\begin{proof}
By continuity and invariance of domain,
$\theta(Q)=\overline U$ and $K_i=\overline{U_i}$, with all interiors open.
The modulo-null union $\bigcup_iU_i$ is dense in $U$, so finiteness gives
$\theta(Q)=\bigcup_iK_i$.  Distinct closures can meet only in
$\theta_i(\partial Q)\cup\theta_j(\partial Q)$, a Lipschitz image of a
$k$-null set.  The ownership rule therefore partitions $B$ exactly, and
\[
 K_i\setminus B_i\subset (\theta(Q)\setminus B)
                 \cup\bigcup_{h<i}(K_i\cap K_h)
\]
is $k$-null.  A $C^1$ embedding of the compact cube and its local inverses
are Lipschitz on a finite cover, so the same conclusions hold in the
embedded image.
\end{proof}

\subsection{The closed red--blue Yomdin tree}
\label{sec:closed-yomdin-tree}

We use a single carrier embedding $\sigma:Q\to M$ throughout each tree and
write
\[
 \mathscr D:=\sigma(Q),\qquad
 T_{\sigma(u)}\mathscr D:=d_u\sigma(\mathbb R^k).
\]
The initial node $o$ is either a blue node carrying one strongly
$\eps$-bounded couple or a red admissible family.  Denote its constituent
index set by $\Theta_o$ and its couples by
$(\sigma,\theta_{o,\alpha})$, $\alpha\in\Theta_o$.  For a blue root,
$|\Theta_o|=1$.  For a red root, only the parameter maps vary; the carrier
is still $\sigma$.  Let $\theta'_o$ be the bounded envelope of the root:
it is the unique parameter map for a blue root and the common envelope
in the definition of an admissible family for a red root.

For every level-$n$ node $v$, let $\Theta_v$ be its finite constituent
index set and let $\theta_{v,\vartheta}:Q\to Q$ be its parameter maps.
The corresponding raw source images and terminal couples are
\[
 D_{v,\vartheta}:=(\sigma\circ\theta_{v,\vartheta})(Q),
 \qquad
 \bigl(g^n\circ\sigma,\theta_{v,\vartheta}\bigr),
 \qquad \vartheta\in\Theta_v.
\]
Thus $g^n(D_{v,\vartheta})$ is the raw terminal image, whereas
$\mathscr D$ and $g^n(\mathscr D)$ are the source and terminal carriers.
A red node has a common bounded terminal envelope $\theta'_v$ containing
the parameter images of all its constituents; for a blue node, its unique
parameter map serves as $\theta'_v$.

The initial owned sets are pairwise disjoint Borel sets
$A_{o,\alpha}\subset D_{o,\alpha}$ satisfying
\begin{equation}\label{eq:initial-owned-data}
 \vol_k(D_{o,\alpha}\setminus A_{o,\alpha})=0,
 \qquad
 A_o:=\bigsqcup_{\alpha\in\Theta_o}A_{o,\alpha}.
\end{equation}
For a fresh construction, fix an order on $\Theta_o$ and take
\[
 A_{o,\alpha}
 :=D_{o,\alpha}\setminus\bigcup_{\beta\prec\alpha}D_{o,\beta}.
\]
At a restart we instead retain the inherited owned sets.  In particular,
$A_o$ need not be the whole carrier.  The initial plaque construction is
the special case of a single couple $(\sigma,\id)$ with $A_o=\mathscr D$.

Ownership is propagated separately from each parent constituent, before
regrouping its offspring into red and blue nodes.  More precisely, for a
parent constituent $P=(v,\vartheta)$, write $A_P=A_{v,\vartheta}$ and order
\emph{all} of its child constituents as $(w_i,\vartheta_i)$, $1\le i\le t$.
With $K_i:=D_{w_i,\vartheta_i}$, set
\begin{equation}\label{eq:tree-recursive-ownership}
 A_{w_i,\vartheta_i}
 :=A_P\cap\left(K_i\setminus\bigcup_{h<i}K_h\right).
\end{equation}
The order includes constituents belonging to different child nodes, not
just those in one red family.  No child node combines constituents from
different parent constituents.  Consequently, every positive depth node
has a unique initial constituent as an ancestor, without requiring an
additional carrier index.

For $|v|=n$, write
\[
 A_v:=\bigsqcup_{\vartheta\in\Theta_v}A_{v,\vartheta},
 \qquad R_{v,\vartheta}:=g^n(A_{v,\vartheta}),
 \qquad R_v:=g^n(A_v).
\]
In particular, $R_{o,\alpha}=A_{o,\alpha}$ and $R_o=A_o$.  The raw charts
carry the geometric estimates; the owned sets carry the exact pointwise
partition.  Restricting an active subset of a node never changes its raw
charts or its owned sets.

\begin{theorem}[Closed red--blue Yomdin tree]
\label{thm:closed-yomdin-tree}
Fix an integer $r\ge2$ and $1\le k\le d-1$.  There are constants
\[
 A_*,C_{\rm Y}>1,
 \qquad
 c_*,\rho_*>0,
\]
depending only on $M,r,d,k$, with the following property.  Let
$\mathscr U\subset\Diff^r(M)$ be a family for which both $\mathscr U$ and
$\mathscr U^{-1}$ are bounded in $C^r$.  There exists
$\eps_*=\eps_*(\mathscr U,r)>0$ such that, for every
$g\in\mathscr U$ and $0<\eps\le\eps_*$, the red--blue construction may be
iterated to any finite depth from a blue node, or from a red admissible
family with at most $A_*$ constituents.  Use the common carrier $\sigma$
and initial owned sets in \eqref{eq:initial-owned-data}, and assume that
every initial parameter map is injective on $(0,1)^k$.  Then all descendant
parameter maps have the same property, every terminal couple
\[
 \bigl(g^n\circ\sigma,\theta_{v,\vartheta}\bigr)
\]
is strongly $\eps$-bounded, and the recursive owned sets satisfy the
following properties.

Set
\begin{equation}\label{eq:tree-complexity}
 \mathfrak M(g):=\max\{2,\|dg\|_\infty,\|dg^{-1}\|_\infty\},
 \qquad
 D_r(g):=C_{\rm Y}\mathfrak M(g)^{\frac{3k^2}{r-1}},
\end{equation}
and
\begin{equation}\label{eq:tree-Lk}
 \ell_k(g,x):=\log\max_{0\le q<k}\|\wedge^q dg_x\|,
 \qquad L_k(g):=\sup_{x\in M}\ell_k(g,x).
\end{equation}
\begin{enumerate}[label=\textnormal{(CP\arabic*)},leftmargin=12mm]
\item For every level $n\ge1$,
\[
 \{A_{v,\vartheta}:|v|=n,\ \vartheta\in\Theta_v\}
\]
is an exact Borel partition of $A_o$.  The owned sets are nested along
constituent branches, and
$\vol_k(D_{v,\vartheta}\setminus A_{v,\vartheta})=0$.
For $y\in A_{v,\vartheta}$,
\[
 R_{v,\vartheta}=g^n(A_{v,\vartheta}),
 \qquad
 T_{g^ny}g^n(\mathscr D)=dg_y^n(T_y\mathscr D).
\]
\item If $v$ is red, then $|\Theta_v|\le A_*$,
\[
 \vol_k(R_v)\ge c_*\eps^k,
\]
and distinct raw constituent targets have zero $k$-volume intersection.
The corresponding owned targets are disjoint.
\item Every node has at most $D_r(g)e^{L_k(g)}$ blue children.
\item If $\operatorname{diam}Q_0<\rho_*\eps$, then $Q_0$ meets the target of
at most $D_r(g)$ children of any fixed node.
\end{enumerate}
All constants are uniform on $\mathscr U$.  The constant
$C_{\rm Y}=C_{M,r,d,k}$ is independent of $g$ and therefore of the iteration
power $p$ when $g=f^p$.
\end{theorem}

\begin{proof}
We use the one-step construction of \cite[Lemma~4.1]{BOB}, applied
separately to each parent constituent.  The additional points to verify
are exact ownership on closed cubes and local target incidence.

\medskip\noindent\emph{One-step estimates and uniformity.}
The scale conditions in \cite[Section~3.2]{BOB} involve finitely many
$C^r$ bounds and a lower singular-value bound.  These are uniform on
$\mathscr U$ and in a fixed finite atlas, so a common $\eps_*>0$ is
available.  Decrease it so that all constituent images and the one-step
envelopes being estimated lie in the chosen coordinate neighborhoods and
\begin{equation}\label{eq:eps-normalization}
 \eps_*\le\sqrt{k}.
\end{equation}
For a strongly $\eps$-bounded parent $(\tau,\theta)$, the 
construction in \cite[Lemma~4.1]{BOB} gives envelopes $\phi$ and full-dimensional maps $\eta$ with
\begin{align}
 \#\{\phi\}&\le B_0\mathfrak M(g)^{3k^2/(r-1)},
       \label{eq:BO-envelope-count}\\
 \theta&>\phi>\phi\circ\eta,
       \label{eq:BO-nesting-order}\\
 (g\circ\tau,\phi)&\text{ bounded},\qquad
 (g\circ\tau,\phi\circ\eta)\text{ strongly }\eps\text{-bounded}.
       \label{eq:BO-nesting}
\end{align}
At level $j$ the parent carrier is $\tau=g^j\circ\sigma$, so the next
carrier is $g^{j+1}\circ\sigma$.  The full dimensional child parameter
images are pairwise disjoint on their interiors and cover the parent
parameter image modulo a $k$-null set.  The red images are grouped into
admissible families of at most $A_0$ members, and the remaining blue images
satisfy
\begin{equation}\label{eq:BO-blue-count}
 \#\{\text{blue images}\}
 \le B_0\mathfrak M(g)^{3k^2/(r-1)}
    \max_{0\le q<k}\|\wedge^q dg\|_{\operatorname{Im}(\tau\circ\theta)}.
\end{equation}
Here $A_0,B_0$ depend only on $M,r,d,k$, and we  take $A_*=A_0$.
These are the estimates of \cite[Lemma~4.1 and its Step~1]{BOB}; all
higher derivative bounds affect the admissible scale, not these constants.

\medskip\noindent\emph{Closed refinement and ownership.}
The algebraic maps in \cite[Lemma~4.2]{BOB} are analytic and injective on
the open cube, with bounded $C^r$ seminorms.  Their derivatives through
order $r-1$ extend continuously to the closed cube, with Lipschitz
bounds at order $r-1$.  Full-dimensional open images lie in $(0,1)^k$
by invariance of domain.  The affine and algebraic factorizations
therefore preserve interior injectivity under composition with the
initial parameter maps.  Lower-dimensional images are $k$-null, as are
their images under the preceding Lipschitz parameter maps.
Thus \cref{lem:closed-ownership} applies at every step.

For a parent constituent $P$, pull $A_P$ back by the carrier embedding
$\sigma$.  This gives a Borel set of full $k$-measure in its raw parameter
image.  Apply \cref{lem:closed-ownership} to all of its child constituents
in the fixed order, and push the resulting sets forward by $\sigma$.
The resulting rule is exactly \eqref{eq:tree-recursive-ownership}, and  in
particular,
\[
 A_P=\bigsqcup_{i=1}^t A_{w_i,\vartheta_i},
 \qquad
 \vol_k(K_i\setminus A_{w_i,\vartheta_i})=0.
\]
Starting with \eqref{eq:initial-owned-data} proves by induction both
exact coverage of $A_o$ and nesting.  Grouping constituents into nodes
only groups already assigned sets and does not reassign their boundaries.
This proves \textnormal{(CP1)}, including its tangent identity by the
chain rule.

A red raw target contains an $\eps$-graph, so its volume is at least
$c_*\eps^k$.  Its owned target differs only by a $k$-null set, since
carrier embeddings and finite iterates preserve such null sets.  Hence
the same lower bound holds for $R_v$.  The raw intersections are null by
the one-step construction, while injectivity of $g^n$ preserves exact
disjointness of the owned constituents.  This proves \textnormal{(CP2)}.
Summing \eqref{eq:BO-blue-count} over at most $A_*$ parent constituents
proves \textnormal{(CP3)} once $C_{\rm Y}\ge A_*B_0$.

\medskip\noindent\emph{Local target incidence.}
For each envelope $\phi$, retain the unit lattice cube labels of the raw
target charts in \cite[proof of Lemma~4.1, Steps~2--4]{BOB}.  In the
normalized target coordinates $\mathcal R_{\phi,\eps}$, at most $B_0$
charts are assigned to one cube.  The polynomial-approximation error
places a chart assigned to $\mathbf m+[0,1]^d$ inside
$\mathbf m+[-2,3]^d$.  The subsequent red grouping does not duplicate
charts.  The normalization consists of a target coordinate chart, a
translation and an orthogonal change of coordinates, and rescaling by
$\eps^{-1}$.  Thus, after enlarging the atlas constant,
\begin{equation}\label{eq:target-rescaling-Lipschitz}
 \operatorname{Lip}(\mathcal R_{\phi,\eps})\le C_{\rm at}\eps^{-1}.
\end{equation}
Choose $C_{\rm at}\rho_*<1/2$.  A set of diameter $<\rho_*\eps$ then has
normalized diameter $<1/2$ and meets at most $7^d$ of these enlarged cubes.
Hence at most $7^dB_0$ raw charts from one envelope meet it.  Each child
whose target meets the set contains a constituent whose raw target meets
it, and different children have disjoint constituent families.  Using
\eqref{eq:BO-envelope-count} and at most $A_*$ parent constituents bounds
the number of such children by
\[
 A_*7^dB_0^2\mathfrak M(g)^{3k^2/(r-1)}.
\]
Taking $C_{\rm Y}\ge\max\{A_*B_0,A_*7^dB_0^2\}$ proves
\textnormal{(CP4)}.  Only $\eps_*$ depends on the common $C^r$ bound.
\end{proof}

\begin{remark}[Carriers and parameter maps]\label{rem:tree-carrier-convention}
The carrier at depth $n$ is $g^n\circ\sigma$, which remains a $C^r$
embedding.  The terminal source map
$\sigma\circ\theta_{v,\vartheta}$ need not be an embedding: interior
injectivity of a parameter map does not imply that its differential has
full rank.  All tangent spaces are those of the carrier.  At a restart,
the inherited carrier and parameter maps remain separate, as specified
in \cref{prop:yomdin-tree-restart} below.
\end{remark}

\begin{lemma}[Local blue child bound]\label{lem:localized-blue-count}
Use the construction of \cref{thm:closed-yomdin-tree}.  Let $a\in\mathbb R$
and $w\in C(M,[0,\infty))$ satisfy
\[
 \ell_k(g,x)\le a+w(x)\qquad(x\in M).
\]
Assume that $w$ has oscillation at most $1$ on every raw terminal
constituent image, including those at level zero.
Fix a node at level $j$, and let $F$ be a subset of its owned source such
that $w(g^jx)\le c$ for all $x\in F$.  Then the number of its blue children
whose owned source meets $F$ is at most
\begin{equation}\label{eq:localized-blue-count}
 D_r(g)e^{a+c+1}.
\end{equation}
\end{lemma}

\begin{proof}
A child whose owned source meets $F$ descends from a constituent with the
same property.  For each such parent constituent choose $x$ in that
intersection.  The raw target of that constituent contains $g^jx$, so
$\ell_k(g,z)\le a+c+1$ throughout that target by the oscillation bound.
Apply \eqref{eq:BO-blue-count} to this complete constituent and sum over at
most $A_*$ constituents.  
\end{proof}

\begin{proposition}[Distortion, orbit control, and restart]
\label{prop:yomdin-tree-restart}
In the setting of \cref{thm:closed-yomdin-tree}, there are
constants $C_{\rm bd},K_0>1$, depending only on $M,r,d,k$, with the following
properties at positive depth.
\begin{enumerate}[label=\textnormal{(CP\arabic*)},leftmargin=12mm,start=5]
\item For every node $v$ at level $n$,
\[
 C_{\rm bd}^{-1}
 \le
 \frac{J_k(g^n|_{T_y\mathscr D})}
      {J_k(g^n|_{T_z\mathscr D})}
 \le C_{\rm bd}
 \qquad(y,z\in\overline{A_v}).
\]
\item Every terminal source map
$\psi=\sigma\circ\theta_{v,\vartheta}$ at depth $N$ satisfies
\[
 \|d(g^j\circ\psi)\|_\infty\le K_0,
 \qquad 0\le j\le N.
\]
\end{enumerate}
The construction may be restarted at any node $v$ of depth $s$ by setting
\begin{equation}\label{eq:tree-restart-data}
 \widetilde\sigma:=g^s\circ\sigma,\qquad
 \widetilde\theta_{o,\vartheta}:=\theta_{v,\vartheta},\qquad
 \widetilde A_{o,\vartheta}:=R_{v,\vartheta},
 \qquad \vartheta\in\Theta_v.
\end{equation}
Here $\widetilde\Theta_o=\Theta_v$ and the initial bounded envelope is
$\widetilde\theta'_o=\theta'_v$.  The initial owned union is
$\widetilde A_o=R_v$, not necessarily the whole new carrier
$\widetilde\sigma(Q)$.  With inherited raw charts and ordering,
the restarted tree is the original descendant subtree in the new time
coordinates and satisfies \textnormal{(CP1)}--\textnormal{(CP6)} with the
same constants.  These constants are uniform on $\mathscr U$.
\end{proposition}

\begin{proof}
Let $v$ be a node at depth $n\ge1$.  The initial bounded envelope is
$\theta'_o$ and the terminal bounded envelope is $\theta'_v$.
Nesting and compactness give
\[
 \sigma^{-1}(\overline{A_v})
 \subset\theta'_o(Q)\cap\theta'_v(Q).
\]
 Choose a finite coordinate atlas of $M$.  For $j=0,n$, let $\kappa_{j}:U_{j}\to\mathbb R^d$
be the coordinate chart used for the corresponding
bounded envelope. Define
\[
a_j(u):=
\left\|
\bigwedge\nolimits^k
D(\kappa_j\circ g^j\circ\sigma)_u
\right\|,
\]
the Euclidean $k$-dimensional volume density
of the parametrized disk in these coordinates.
Let
\[
J_j(u):=
\left\|
\bigwedge\nolimits^k
D(g^j\circ\sigma)_u
\right\|_{\mathrm{Riem}}
\]
be its Riemannian volume density, and put
\[
b_j(u):=\frac{J_j(u)}{a_j(u)}.
\]
By uniform equivalence of the coordinate and
Riemannian metrics on the fixed finite atlas,
there exists $C_{\mathrm{met}}\geq1$ such that
\[
C_{\mathrm{met}}^{-1}\leq b_j(u)
\leq C_{\mathrm{met}}.
\]
For $y=\sigma(u)$ and $z=\sigma(w)$ in $\overline{A_v}$,
\begin{equation}\label{eq:two-end-distortion}
 \frac{J_k(g^n|_{T_y\mathscr D})}
      {J_k(g^n|_{T_z\mathscr D})}
 =\frac{a_n(u)}{a_n(w)}\frac{a_0(w)}{a_0(u)}
  \frac{b_n(u)}{b_n(w)}\frac{b_0(w)}{b_0(u)}.
\end{equation}
Apply \cite[Lemma~3.4]{BOB} to the bounded couples
$(\sigma,\theta'_o)$ and $(g^n\circ\sigma,\theta'_v)$.
The first two factors are at most $\sqrt2$ each.  Interchanging $u,w$
proves \textnormal{(CP5)} with $C_{\rm bd}=2C_{\mathrm{met}}^4$.

For a terminal source map $\psi_N=\sigma\circ\theta_N$, let $\theta_j$
be its ancestral constituent parameter map at time $j$, including its
initial constituent at $j=0$.  The edge contractions in
\eqref{eq:BO-nesting-order} give
\[
 \theta_N=\theta_j\circ h_{j,N},\qquad
 \|dh_{j,N}\|_\infty\le2^{j-N}.
\]
At $j=N$ take $h_{N,N}=\id$.  Strong boundedness of every ancestral couple,
metric comparison, and \eqref{eq:eps-normalization} imply
\[
 \|d(g^j\circ\psi_N)\|_\infty
 \le C_{\rm at}\frac{\eps}{\sqrt{k}}2^{j-N}\le C_{\rm at}
 \qquad(0\le j\le N).
\]
After enlarging $C_{\rm at}>1$, this proves \textnormal{(CP6)}.

At a restart, $\widetilde\sigma=g^s\circ\sigma$ is a $C^r$ embedding,
and its initial couples are exactly the inherited terminal couples.
They are strongly $\eps$-bounded and their parameter maps are internally
injective.  A red root retains its admissible family and the bounded
envelope $\theta'_v$; a blue root retains its unique constituent as its
envelope.  
The inherited owned sets are disjoint and satisfy
\[
 \vol_k\bigl((\widetilde\sigma\circ\theta_{v,\vartheta})(Q)
                \setminus R_{v,\vartheta}\bigr)=0.
\]
Thus they satisfy \eqref{eq:initial-owned-data}.  Keep the original raw
offspring and the same order at each parent constituent.  Since $g^s$ is
injective, it commutes with the intersections and set differences in
\eqref{eq:tree-recursive-ownership}.  For a descendant $w$ with
$|w|=s+n$ and $\vartheta\in\Theta_w$, induction identifies its new owned
source set with $g^s(A_{w,\vartheta})$.  Its relative terminal target is
$g^{s+n}(A_{w,\vartheta})$, exactly the original target.

The one-step volume and incidence estimates are consequently unchanged.
For restarted distortion, the initial comparison is made on the
inherited bounded envelope, not on all of $Q$ with the identity
parameter map.  The preceding two endpoint calculation applies with the
same constant.  The edge factorizations also remain unchanged.  If
$\widetilde\psi_N$ is a restarted terminal source map at relative depth
$N$, then $\widetilde\psi_N=g^s\circ\psi_{s+N}$ for the corresponding
original source map.  Applying the original factorization along the whole
branch gives \textnormal{(CP6)} over the full time interval, rather than a
product of block bounds.  Likewise, distortion over the full branch uses
only the original and final envelopes.  At depth zero the Jacobian ratio
is one.  This proves the restart assertion.
\end{proof}

\subsection{Spectral gap, unstable lifts, and the ergodic theorem}
\label{sec:core-statement}

At every Oseledets-regular point, Lyapunov exponents are listed with
multiplicity in nonincreasing order,
\[
 \chi_1(x,f)\ge\cdots\ge\chi_d(x,f).
\]
The integrand in \eqref{eq:def-Splus} is measurable and bounded above by
$d\log^+\|df\|_\infty$, so $S^+(\nu,f)$ is finite for every invariant
probability $\nu$.

For $1\le k\le d-1$, let
\[
 \wh M_k:=\{(x,E):x\in M,\ E\in\operatorname{Gr}_k(T_xM)\}
\]
be the Grassmann bundle of $k$-planes, with projection $\pi(x,E)=x$.
The derivative of $f$ induces the map and the continuous potential
\[
 F_f(x,E):=(f(x),df_xE),
 \qquad
 \Phi_f(x,E):=\log J_k(f,x,E).
\]
Suppose that an invariant probability $\nu$ has no zero Lyapunov exponent
and that its positive Oseledets space $E^u_{\nu}$ has dimension $k$ almost everywhere.
Its \emph{unstable lift} is the graph measure
\begin{equation}\label{eq:unstable-lift}
 \wh\nu^u:=\bigl(x\mapsto(x,E^u_\nu(x))\bigr)_*\nu.
\end{equation}
Then
\begin{equation}\label{eq:positive-sum}
 S^+(\nu,f)=\int_{\wh M_k}\Phi_f\,\dd\wh\nu^u.
\end{equation}

For $1\le k\le d-1$, write
\begin{equation}\label{eq:hyperbolic-index-class}
 \mathscr H_k(g):=
 \{\rho\in\mathcal M_e(g):
       \chi_k(\rho,g)>0>\chi_{k+1}(\rho,g)\}.
\end{equation}
Thus $\mathscr H_k(g)$ consists of the hyperbolic ergodic measures of
unstable index $k$.  This is a Borel set.  

We first prove continuity for ergodic approximating measures of a fixed
index.  This identifies the full unstable graph measure, not just its
Jacobian integral.  The arbitrary measure extension is
\cref{thm:fixed-index-continuity} below.

\begin{theorem}[Fixed-index continuity of unstable lifts]\label{thm:core}
Let $f_n\to f$ in $C^\infty$, fix $1\le k\le d-1$, and let
$\nu_n\in\mathscr H_k(f_n)$.  Assume
\begin{equation}\label{eq:core-assumptions}
 \nu_n\longrightarrow\mu,
 \qquad h_{\nu_n}(f_n)\longrightarrow h_\mu(f),
\end{equation}
and, for $\tau_\mu$-almost every $\rho$,
\begin{equation}\label{eq:core-componentwise-gap}
 \rho\in\mathscr H_k(f),
 \qquad h_\rho(f)>\lambda_{k-1}(f).
\end{equation}
Then
\begin{equation}\label{eq:core-lift-convergence}
 \wh\nu_n^u\longrightarrow\wh\mu^u
\end{equation}
on $\wh M_k$, and consequently
\begin{equation}\label{eq:core-Splus}
 S^+(\nu_n,f_n)\longrightarrow S^+(\mu,f).
\end{equation}
\end{theorem}

The proof occupies \cref{sec:record-times,sec:plaque-selection,sec:ergodic-continuity}. 
The index selection lemma below is used only when deriving the
entropy gap statements in the introduction from the fixed-index criteria.

\bigskip

\paragraph{\it The entropy gap and the unstable index.}
The next lemmas convert an entropy threshold above the second-largest
exterior growth rate into a robust unstable index selection, and record how a
componentwise entropy gap behaves under powers and invariant submeasures.

\begin{lemma}[Finite iterate bound for exterior growth]
\label{lem:finite-iterate}
Let $I\subset\{0,\ldots,d-1\}$ be finite, and let
\[
 A>\max_{q\in I}\lambda_q(f).
\]
There exists $p_0\ge1$ such that, for every $p\ge p_0$,
\begin{equation}\label{eq:finite-iterate-bound}
 \max_{q\in I}\log\|\wedge^qdf^p\|_\infty<pA.
\end{equation}
For each fixed $p\ge p_0$, there exist a $C^1$ neighborhood $\mathscr V_p$
of $f$ and a number $\delta_p>0$ such that
\[
 \max_{q\in I}\log\|\wedge^qdg^p\|_\infty
 \le pA-\delta_p
 \qquad(g\in\mathscr V_p).
\]
\end{lemma}

\begin{proof}
For every $q$, the sequence
$\log\|\wedge^qdf^n\|_\infty$ is subadditive, and thus $\frac{1}{n}\log\|\wedge^qdf^n\|_\infty$
converges to $\lambda_q(f)$ as $n\to \infty$.  Since $I$ is finite, the convergence gives
\eqref{eq:finite-iterate-bound} simultaneously for all $q\in I$.  For fixed
$p$, the left hand side depends continuously on $g$ in the $C^1$ topology, and 
the strict inequality therefore persists with a positive margin.
\end{proof}

\begin{lemma}[An entropy threshold determines the unstable index]
\label{lem:index}Let $f\in  \mathcal{S}^1(M)$.
Then there is a (unique)
$k\in\{1,\ldots,d-1\}$ such that
\begin{equation}\label{eq:unique-index}
 \lambda_k(f)=\max_{\substack{0\le q\le d-1}}\lambda_q(f)\geq h_{\mathrm{top}}(f)>
 \max_{\substack{0\le q\le d-1\\q\ne k}}\lambda_q(f)
 =\lambda_{\mathrm{sec}}(f).
\end{equation}
Moreover any $\eta\in \mathcal M_e(f)$ with $h_\eta(f)> \lambda_{\mathrm{sec}}(f)$ belongs to $\mathscr H_k(f)$ with  $\chi_k(\eta,f), -\chi_{k+1}(\eta, f)\geq h_\eta(f)-\lambda_{\mathrm{sec}}(f)>0$.

\end{lemma}

\begin{proof}
Take  $\eta\in \mathcal M_e(f)$ such that $h_{\eta}(f)>\lambda_{\mathrm{sec}}(f)$. Let $u$ be the number of positive Lyapunov exponents of $\eta$.  The cases
$u=0$ and $u=d$ imply  $h_\eta(f)=0$ by Margulis-Ruelle's inequality for $f$ and
$f^{-1}$, respectively.  Hence $1\le u\le d-1$, and Margulis-Ruelle's inequality
and Oseledets' theorem give
\[
 h_\eta(f)
 \le\sum_{i=1}^{u}\chi_i(\eta,f)
 \le\lambda_u(f).
\]
Thus $\lambda_u(f)>\lambda_{\mathrm{sec}}(f)$. But at  most one exterior degree can have growth larger
than $\lambda_{\mathrm{sec}}(f)$.  Together with the variational principle, this proves the existence and uniqueness of $k(=u)$  in 
\eqref{eq:unique-index}. We may also write 
\begin{align*}
 h_\eta(f)&
 \le \sum_{i=1}^{u-1}\chi_i(\eta,f) +\chi_{u}(\eta,f),\\
 &\le\lambda_{u-1}(f)+\chi_{u}(\eta,f),\\
 & \leq \lambda_{\mathrm{sec}}(f)+\chi_{u}(\eta,f)
 \end{align*}
 and similarly for $u\leq d-2$
 \begin{align*}
 h_\eta(f)&
 \le \sum_{i=1}^{u+1}\chi_i(\eta,f) -\chi_{u+1}(\eta,f),\\
 &\le\lambda_{u+1}(f)-\chi_{u+1}(\eta,f),\\
 & \leq\lambda_{\mathrm{sec}}(f)-\chi_{u+1}(\eta,f).
 \end{align*}
 For $u=d-1$, apply Ruelle's inequality to $f^{-1}$:
\begin{equation}\label{eq:weakest-stable-uniform-edge}
 h_\eta(f)\le  -\chi_d(\eta,f).
\end{equation}

Therefore in any case we have (recall $k=u$) $\chi_k(\eta,f), -\chi_{k+1}(\eta, f)\geq h_\eta(f)-\lambda_{\mathrm{sec}}(f)>0$, in particular $\eta\in \mathscr H_k(f)$.


\end{proof}

\begin{corollary} \label{cor:openess}
For any $r>1$ the map $f\mapsto \max\left(0,h_{\mathrm{top}}(f)-\lambda_{\mathrm{sec}}(f)\right)$ is lower semicontinuous on $\mathrm{Diff}^r(M)$ for the $C^1$ topology. Therefore $\mathcal S^r(M)$ is $C^1$-open in $\mathrm{Diff}^r(M)$.  
\end{corollary}
\begin{proof} It is enough to show the lower semicontinuity of $g\mapsto h_{\mathrm{top}}(g)-\lambda_{\mathrm{sec}}(g)$ at some $f\in \mathcal S^r(M)$ for the $C^1$ topology. 
By  Lemma \ref{lem:index} any $f$-invariant ergodic measure with large enough entropy is  hyperbolic. It follows then from Katok's horseshoe theorem that $$\liminf_{g\xrightarrow{C^1}f}h_{\mathrm{top}}(g)\geq h_{\mathrm{top}}(f).$$ Moreover by  Lemma \ref{lem:index} there is a unique $k\in \{1,\cdots,d-1\}$ such that $$\lambda_k(f)=\max_{\substack{0\le q\le d-1}}\lambda_q(f)\geq h_{\mathrm{top}}(f)>
 \max_{\substack{0\le q\le d-1\\q\ne k}}\lambda_q(f)
 =\lambda_{\mathrm{sec}}(f).$$

By  Lemma \ref{lem:finite-iterate} the map $g\mapsto \max_{\substack{0\le q\le d-1\\q\ne k}}\lambda_q(g)$ is upper semicontinuous for the $C^1$ topology. Therefore for $g$ in a $C^1$ neighborhood of  $f$, we have 
$h_{\mathrm{top}}(g)>\max_{\substack{0\le q\le d-1\\q\ne k}}\lambda_q(g)$. As already observed,  $h_{\mathrm{top}}(g)\leq \max_{\substack{0\le q\le d-1}}\lambda_q(g)$, so that we  necessarily have $\lambda_k(g)= \max_{\substack{0\le q\le d-1}}\lambda_q(g)$ and $\lambda_{\mathrm{sec}}(g)=\max_{\substack{0\le q\le d-1\\q\ne k}}\lambda_q(g)$.  Therefore 
\begin{align*}
\liminf_{g\xrightarrow{C^1}f}\left(h_{\mathrm{top}}(g)-\lambda_{\mathrm{sec}}(g)\right)&\geq \liminf_{g\xrightarrow{C^1}f}h_{\mathrm{top}}(g)-\limsup_{g\xrightarrow{C^1}f}\lambda_{\mathrm{sec}}(g),\\
&\geq \liminf_{g\xrightarrow{C^1}f}h_{\mathrm{top}}(g)-\limsup_{g\xrightarrow{C^1}f}\max_{\substack{0\le q\le d-1\\q\ne k}}\lambda_q(g),\\
 & \geq h_{\mathrm{top}}(f)-\max_{\substack{0\le q\le d-1\\q\ne k}}\lambda_q(f)= h_{\mathrm{top}}(f)-\lambda_{\mathrm{sec}}(f).
\end{align*}
\end{proof}

\begin{lemma}[Componentwise gaps under powers and invariant domination]
\label{lem:componentwise-gap-inheritance}
Let $\ell\in\mathbb R$ and let $\mu\in\mathcal M(f)$ satisfy
\begin{equation}\label{eq:abstract-componentwise-gap}
 h_\rho(f)>\ell
 \quad\text{for }\tau_\mu\text{-almost every }\rho\in\mathcal M_e(f).
\end{equation}
Then:
\begin{enumerate}[label=\textnormal{(\roman*)},leftmargin=10mm]
\item for every $p\ge1$, almost every $f^p$-ergodic component $\eta$ of
      $\mu$ satisfies
      \begin{equation}\label{eq:componentwise-gap-power}
       h_\eta(f^p)>p\ell;
      \end{equation}
\item if $\alpha\in\mathcal M(f)$ and $0<\beta\le1$ satisfy
      $\beta\alpha\le\mu$, then
      \begin{equation}\label{eq:componentwise-gap-domination}
       h_\rho(f)>\ell
       \quad\text{for }\tau_\alpha\text{-almost every }\rho\in\mathcal M_e(f).
      \end{equation}
\end{enumerate}
\end{lemma}

\begin{proof}
Let $\mu=\int\mu_\omega\,\dd\tau(\omega)$ be its $f$-ergodic
 decomposition.  For $\tau$-almost every $\omega$, the $f^p$-ergodic
 decomposition of $\mu_\omega$ is a finite cyclic family, and every member
 of that family has $f^p$-entropy $p h_{\mu_\omega}(f)$.  Integrating these
 finite decompositions and using \eqref{eq:abstract-componentwise-gap}
 proves \eqref{eq:componentwise-gap-power}.

For \textnormal{(ii)}, let $v=\dd(\beta\alpha)/\dd\mu$.  Invariance gives
$v\circ f=v$ $\mu$-almost everywhere.  Hence $v$ is $\mu_\omega$-almost everywhere constant for
$\tau$-almost every $\omega$; denote this constant by
$c(\omega)\in[0,1]$.  Therefore
\[
 \beta\alpha=\int c(\omega)\mu_\omega\,\dd\tau(\omega),
 \qquad
 \int c(\omega)\,\dd\tau(\omega)=\beta.
\]
The ergodic decomposition probability of $\alpha$ is obtained by reweighting
this same family by $c/\beta$, so it is supported on the full measure set in
\eqref{eq:abstract-componentwise-gap}.
\end{proof}

\subsection{Local exterior growth and finite-time margins}

The pointwise quantity $\ell_k(g,x)$ in \eqref{eq:tree-Lk} controls the blue
children.  A different, smaller quantity suffices to turn $k$-volume growth
into expansion of every tangent vector:
\begin{equation}\label{eq:L-U}
 A_{k-1}(g):=\log\|\wedge^{k-1}dg\|_\infty.
\end{equation}
For $k=1$ this is zero.  For any $a>\lambda_{k-1}(f)$,
\cref{lem:finite-iterate} gives
\begin{equation}\label{eq:L-margin}
 A_{k-1}(f^p)<pa
\end{equation}
for all sufficiently large $p$, with a positive margin under perturbation
once $p$ is fixed.  We do not replace all lower exterior powers in the blue
count by the $(k-1)$-st one.  Instead, we retain their local excess in the
covering estimate and show in \cref{sec:ergodic-continuity} that its average,
divided by $p$, tends to zero.

\section{Record times and neutral coverings}\label{sec:record-times}
\label{sec:record-neutral}

In this section, by taking a hyperbolic threshold above the second growth rate, we introduce a record-time construction that separates an orbit segment into selected times,
carrying backward $k$-Jacobian estimates, and neutral blocks, whose prefixes
grow below the threshold.  The additional entropy  of recording the selected pattern tends to zero with
the coarse time scale, whereas neutral blocks admit closed red--blue
Yomdin covers.

\subsection{Record times and backward hyperbolic inequalities}

Fix $N\ge1$, a diffeomorphism $g$, and a threshold
$\Theta\in\mathbb R$.  Let
\[
 \wh x,G\wh x,\ldots,G^N\wh x,
 \qquad G:=F_g,
\]
be an orbit segment in $\wh M_k$.  Define the centered Jacobian sums by
\begin{equation}\label{eq:record-potential}
 S_0(\wh x):=0,
 \qquad
 S_n(\wh x):=\sum_{j=0}^{n-1}
 \bigl(\Phi_g(G^j\wh x)-\Theta\bigr)
 \quad(1\le n\le N).
\end{equation}
Starting from $t_0:=0$, define recursively
\[
 t_{i+1}:=\min\{t\in\{t_i+1,\ldots,N\}:S_t\ge S_{t_i}\},
\]
whenever this set is nonempty.  We call the resulting times
$0=t_0<t_1<\cdots<t_s\le N$ the \emph{record times}.  If no further record
exists after $t_s<N$, the remaining segment after $t_s$ is called the
\emph{final tail}.

The two elementary consequences below are the reason for this definition.
Between consecutive records, every proper prefix is neutral; at a record,
the entire past satisfies a backward hyperbolic inequality.

\begin{lemma}[Record-time inequalities]
\label{lem:record-properties}
Let $u<b$ be two consecutive record times.  Then
\begin{equation}\label{eq:neutral-prefix}
 \sum_{j=u}^{u+\ell-1}\Phi_g(G^j\wh x)<\ell\Theta
 \qquad(1\le\ell<b-u).
\end{equation}
If the construction stops at $t_s<N$, then
\begin{equation}\label{eq:neutral-final-tail}
 \sum_{j=t_s}^{t_s+\ell-1}\Phi_g(G^j\wh x)<\ell\Theta
 \qquad(1\le\ell\le N-t_s).
\end{equation}
Every positive record time $b$ also satisfies
\begin{equation}\label{eq:finite-hyperbolic-past}
 \sum_{j=b-\ell}^{b-1}\Phi_g(G^j\wh x)\ge\ell\Theta
 \qquad(1\le\ell\le b).
\end{equation}
\end{lemma}

\begin{proof}
Minimality of the next record gives $S_{u+\ell}<S_u$ for
$1\le\ell<b-u$; the same strict inequality holds through the final tail
when no next record exists.  Inductively, every record is a past maximum:
$S_b=\max_{0\le q\le b}S_q$.  Hence $S_b-S_{b-\ell}\ge0$ for
$1\le\ell\le b$, which is the backward inequality.
\end{proof}

For $M\ge1$, a complete record gap $[t_i,t_{i+1}]$ is called \emph{short}
when $t_{i+1}-t_i\le M$ and \emph{long} otherwise.  We select every positive
record time and, in addition, every interior time of a short gap:
\begin{equation}\label{eq:selected-times}
 E_N^M(\wh x):=
 \{t_i:1\le i\le s\}
 \cup
 \bigcup_{t_{i+1}-t_i\le M}
 \{t_i+1,\ldots,t_{i+1}-1\}
 \subset\{1,\ldots,N\}.
\end{equation}
Thus a short gap contributes all times after its entrance, including its
terminal record, whereas a long gap contributes only its terminal record.
The omitted times form the neutral blocks used in the covering argument.

\begin{lemma}[Selected times and neutral blocks]
\label{lem:exact-time-decomposition}
Let $0=t_0<\cdots<t_s\le N$ be the record sequence and set
\[
 \mathcal L_M:=\{0\le i<s:t_{i+1}-t_i>M\}.
\]
For $i\in\mathcal L_M$, define the interior of the long gap by
\[
 I_i:=\{t_i+1,\ldots,t_{i+1}-1\},
 \qquad
 \ell_i:=|I_i|=t_{i+1}-t_i-1.
\]
The final tail is encoded by
\[
 (I_*,\ell_*):=
 \begin{cases}
 (\{t_s+1,\ldots,N\},N-t_s),&t_s<N,\\
 (\varnothing,0),&t_s=N.
 \end{cases}
\]
Then:
\begin{enumerate}[label=\textnormal{(\roman*)},leftmargin=10mm]
\item the following is an exact disjoint decomposition:
      \begin{equation}\label{eq:exact-time-decomposition}
       \{1,\ldots,N\}
       =E_N^M(\wh x)\sqcup
         \bigsqcup_{i\in\mathcal L_M}I_i
         \sqcup I_*;
      \end{equation}
\item every complete neutral block has length at least $M$, and
      \begin{equation}\label{eq:complete-neutral-length}
       \ell_i\ge M\quad(i\in\mathcal L_M),
       \qquad
       \sum_{i\in\mathcal L_M}\ell_i+\ell_*
       =N-|E_N^M(\wh x)|;
      \end{equation}
\item if $I=\{u+1,\ldots,u+\ell\}$ is one of the sets $I_i$ or $I_*$,
      then every sum beginning at its entrance time $u$ is neutral:
      \begin{equation}\label{eq:anchored-neutral-prefix}
       \sum_{j=0}^{q-1}\Phi_g(G^{u+j}\wh x)<q\Theta
       \qquad(1\le q\le\ell).
      \end{equation}
\end{enumerate}
If $E_N^M$ is constant on a set $F$, then the interval components of its
complement are common to all $\wh x\in F$, and
\eqref{eq:anchored-neutral-prefix} holds for every $\wh x\in F$ with the
common entrance time of the corresponding component.
\end{lemma}

\begin{proof}
Each complete gap contributes all of $\{t_i+1,\ldots,t_{i+1}\}$ when
short, and only $t_{i+1}$ when long.  The final tail contributes no selected
time.  This proves the disjoint decomposition and the length identities;
a long gap has $t_{i+1}-t_i\ge M+1$.  The prefix inequalities are
\cref{lem:record-properties}, applied at the entrance record.
If $E_N^M$ is fixed, its complementary interval components and their
entrances are fixed.  Each component is a long gap interior or the final
tail for every point, so the same inequalities hold with these common
entrances.
\end{proof}

\begin{remark}[Endpoint conventions]\label{rem:record-boundary-cases}
A gap of length $M$ is entirely selected, while a gap of length $M+1$ has a
neutral interior of length $M$.  Moreover $I_*=\varnothing$ when $t_s=N$,
and if there is no positive record time then
$E_N^M=\varnothing$ and $I_*=\{1,\ldots,N\}$.  Thus every transition time is
assigned exactly once.
\end{remark}

\begin{lemma}[Combinatorics of the selected time sets]
\label{lem:selected-time-sets}
For fixed $N$ and $M$, the selected time map
$\wh x\mapsto E_N^M(\wh x)$ has the following properties:
\begin{enumerate}[label=\textnormal{(\roman*)},leftmargin=10mm]
\item it is Borel, and $E_N^M(\wh x)\subset E_N^{M+1}(\wh x)$;
\item every set $E_N^M(\wh x)$ is a union of at most $N/M+2$ integer
      intervals;
\item there are $N_0(M)$ and $\kappa_M\ge0$ such that
      \begin{equation}\label{eq:selected-set-count}
       \#\{E_N^M(\wh x):\wh x\in\wh M_k\}
       \le e^{N\kappa_M}
       \qquad(N\ge N_0(M)),
       \qquad
       \kappa_M\xrightarrow[M\to\infty]{}0;
      \end{equation}
\item if $j\in E_N^M(\wh x)$, then
      \[
       j=b-r
      \]
      for some record time $b$ and some integer $0\le r<M$.
\end{enumerate}
\end{lemma}

\begin{proof}
Each proposed record sequence is specified by finitely many inequalities
$S_{t_{i+1}}\ge S_{t_i}$ and $S_q<S_{t_i}$ between continuous functions,
including the final tail inequalities when $t_s<N$.  There are finitely
many proposed sequences, so the selected set map is Borel.  Monotonicity
in $M$ follows from \eqref{eq:selected-times}.

By \cref{lem:exact-time-decomposition}, the complement of $E_N^M$ has at
most $Q_N:=\lfloor N/M\rfloor+1$ interval components.  Hence $E_N^M$ has
at most $N/M+2$ components.  The left and right endpoint sets of its
complement determine $E_N^M$, giving
\begin{equation}\label{eq:selected-endpoint-count}
 \#\{E_N^M\}
 \le\left(\sum_{j=0}^{Q_N}\binom Nj\right)^2.
\end{equation}
For $1\le Q\le N/2$, the binomial theorem with $x=Q/(N-Q)$ gives
$\sum_{j=0}^Q\binom Nj\le x^{-Q}(1+x)^N$, or
\begin{equation}\label{eq:exact-binomial-tail}
\begin{gathered}
 \sum_{j=0}^Q\binom Nj\le e^{N h_{\rm bin}(Q/N)},\\
 h_{\rm bin}(t):=-t\log t-(1-t)\log(1-t).
\end{gathered}
\end{equation}
For $M\ge4$ and all sufficiently large $N$, we have $Q_N\le2N/M\le N/2$.
Since $h_{\rm bin}$ is increasing on $[0,1/2]$, we may take
$\kappa_M=2h_{\rm bin}(2/M)\to0$.  For $M<4$, take $\kappa_M=\log2$.

Finally, take $b=j$ for a selected record time.  Every other selected time
lies in a short gap ending at a record $b$, so $0<b-j<M$.
\end{proof}

The backward inequalities at record times will be encoded by the following
compact sets.  For $R\ge1$, define
\begin{equation}\label{eq:HThetaR}
 H_{\Theta,R}(g):=
 \left\{\wh y:
 \sum_{j=1}^{\ell}\Phi_g(G^{-j}\wh y)\ge\ell\Theta
 \quad\text{for }1\le\ell\le R\right\}
\end{equation}
and
\begin{equation}\label{eq:HTheta}
 H_\Theta(g):=\bigcap_{R\ge1}H_{\Theta,R}(g).
\end{equation}
For fixed $R$, the set $H_{\Theta,R}(g)$ is compact: it is the intersection
of the finitely many closed sets defined by the inequalities in
\eqref{eq:HThetaR}.  Hence $H_\Theta(g)$ is compact as well. Recall  $A_{k-1}(g):=\log\|\wedge^{k-1}dg\|_\infty$.\\ 

\begin{lemma}[Conorm estimate at a backward $k$-volume hyperbolic point, cf.\ the proof of Lemma~5.4 in \cite{BOB}]
\label{lem:conorm-expansion}
Assume  $\gamma:=\Theta-A_{k-1}(g)>0$.  If $\wh y=(y,E)\in H_\Theta(g)$, then for every
$\ell\ge1$,
\begin{equation}\label{eq:inverse-contraction}
 \left\|
 \left(dg^\ell_{g^{-\ell}y}|_{E_{-\ell}}\right)^{-1}
 \right\|
 \le e^{-\gamma\ell},
 \qquad E_{-\ell}=dg^{-\ell}_yE.
\end{equation}
Equivalently,
\begin{equation}\label{eq:conorm-expansion}
 \conorm\bigl(dg^\ell_{g^{-\ell}y}|_{E_{-\ell}}\bigr)
 \ge e^{\gamma\ell}.
\end{equation}
\end{lemma}

\subsection{Neutral block counting}

We now estimate the branching of the red--blue tree during a time interval
on which every initial prefix has $k$-Jacobian growth below the threshold
$\Theta$.  Red nodes are controlled by their definite target volume, while
blue nodes satisfy a one-step recurrence.  Their combination is the basic
counting estimate for a neutral block.

Fix an integer $r\ge2$.  Every initial node in this subsection and the
next satisfies the hypotheses of \cref{thm:closed-yomdin-tree}: it has
one common carrier, at most $A_*$ constituents, internally injective
parameter maps, and initial owned sets as in \eqref{eq:initial-owned-data}.
This includes the inherited data at a restart.
Take $a\in\mathbb R$ and
$w\in C(M,[0,\infty))$ such that
\begin{equation}\label{eq:local-counting-setup}
 \ell_k(g,x)\le a+w(x)\quad(x\in M),
 \qquad \Theta>a+1.
\end{equation}
Choose the admissible scale so that $w$ has oscillation at most $1$ on each
raw terminal constituent image, including the initial images.  This is
possible because strongly $\eps$-bounded constituent images have diameter
at most $C_M\eps$.
For a $C^r$-bounded family of maps and an equicontinuous family of weights,
the scale can be chosen uniformly.  Put
\begin{equation}\label{eq:D-L-gamma}
 D:=D_r(g),\qquad \gamma:=\Theta-a-1>0.
\end{equation}
The records and the selected-time sets are still defined using the
unmodified potential $\Phi_g-\Theta$.

\begin{lemma}[Red nodes in a neutral block]\label{lem:red-neutral-count}
Let $o$ be an initial node as in \cref{thm:closed-yomdin-tree}, with common
carrier $\mathscr D=\sigma(Q)$ and initial owned source $A_o$.
Let $\varnothing\ne F\subset A_o$, fix an integer $s\ge1$, and assume
\begin{equation}\label{eq:neutral-on-F}
 \sum_{j=0}^{s-1}\Phi_g\bigl(G^j(x,T_x\mathscr D)\bigr)<s\Theta
 \qquad(x\in F).
\end{equation}
Set
\[
 \mathcal R_s(F):=
 \{v:|v|=s,\ v\text{ is red},\ A_v\cap F\ne\varnothing\}.
\]
Then
\begin{equation}\label{eq:red-count-neutral}
 |\mathcal R_s(F)|\le C_1e^{s\Theta},
 \qquad
 C_1:=A_*C_Mk^{-k/2}C_{\rm bd}c_*^{-1}.
\end{equation}
\end{lemma}

\begin{proof}
For each active red
node $v$ (i.e., $A_v$ has nonempty intersection with $F$), choose $x_v\in A_v\cap F$.  The tangent identity in
\textnormal{(CP1)}, the Jacobian cocycle, and \eqref{eq:neutral-on-F} give
$J_k(g^s|_{T_{x_v}\mathscr D})<e^{s\Theta}$.
By \textnormal{(CP5)}, this Jacobian is at most $C_{\rm bd}e^{s\Theta}$
throughout $A_v$.  The volume formula and \textnormal{(CP2)} imply
\[
 c_*\eps^k\le\vol_k(g^sA_v)
 \le C_{\rm bd}e^{s\Theta}\vol_k(A_v).
\]
The sets $A_v$ at level $s$ are disjoint subsets of $A_o$.  Moreover,
\eqref{eq:strong-chart-volume-front} gives
\begin{equation}\label{eq:initial-owned-volume-bound}
 \vol_k(A_o)
 \le\sum_{\alpha\in\Theta_o}\vol_k(D_{o,\alpha})
 \le A_*C_Mk^{-k/2}\eps^k.
\end{equation}
Summing over the active red nodes yields
\[
 |\mathcal R_s(F)|c_*\eps^k
 \le C_{\rm bd}e^{s\Theta}\vol_k(A_o)
 \le A_*C_Mk^{-k/2}C_{\rm bd}e^{s\Theta}\eps^k,
\]
which proves the estimate. 
\end{proof}

\begin{lemma}[Reparametrizations of one neutral block]
\label{lem:neutral-block}
Let $o$ be an initial node as in \cref{thm:closed-yomdin-tree}, with common
carrier $\mathscr D=\sigma(Q)$ and initial owned source $A_o$.
Fix an integer $\ell\ge1$ and let $\varnothing\ne F\subset A_o$.  Assume
\begin{equation}\label{eq:neutral-block-all-prefixes}
 \sum_{j=0}^{q-1}\Phi_g\bigl(G^j(x,T_x\mathscr D)\bigr)<q\Theta
 \qquad(x\in F,\ 1\le q\le\ell),
\end{equation}
and suppose that nonnegative integers $c_0,\ldots,c_{\ell-1}$ satisfy
$w(g^jx)\le c_j$ for all $x\in F$ and $0\le j<\ell$.
Let $\mathcal V_q(F)$ be the level-$q$ nodes whose owned source meets $F$,
and put $n_q=|\mathcal V_q(F)|$.  Then
\begin{equation}\label{eq:neutral-block-proper}
 n_q\le C_\gamma D^q
       \exp\!\left(q\Theta+\sum_{j=0}^{q-1}c_j\right),
 \qquad 0\le q\le\ell,
 \qquad C_\gamma:=1+\frac{C_1}{1-e^{-\gamma}}.
\end{equation}
\end{lemma}

\begin{proof}
For $q\ge1$, \cref{lem:red-neutral-count} bounds the number of all active red
nodes at level $q$ by $C_1e^{q\Theta}$.  Its proof uses a point of $F$ to
control the Jacobian on the complete owned source of each active node;
these owned sets have not been cut by $F$.

An active blue child has an active parent.  For each active parent $v$
at level $q-1$, apply \cref{lem:localized-blue-count} to $F\cap A_v$ with
$c=c_{q-1}$.  Thus
\begin{equation}\label{eq:neutral-recurrence}
 n_q\le D e^{a+c_{q-1}+1}n_{q-1}+C_1e^{q\Theta},
 \qquad n_0=1.
\end{equation}
The equality $n_0=1$ counts the single initial node, whether it is blue or
red, not its constituents.  Set
\[
 \widetilde n_q:=
 \frac{n_q}{D^q\exp(q\Theta+\sum_{j=0}^{q-1}c_j)}.
\]
Since $D\ge1$ and $c_j\ge0$, the recurrence gives
\[
 \widetilde n_q\le e^{-\gamma}\widetilde n_{q-1}+C_1D^{-q},
 \qquad \widetilde n_0=1.
\]
Iteration and the geometric-series bound imply
$\widetilde n_q\le C_\gamma$, including $q=0$ with the empty sum.
\end{proof}

\subsection{Neutral block covering}

The preceding block estimate can now be multiplied over the exact
selected neutral decomposition of the time interval.  At a selected time,
the fixed $\cQ$-name controls the number of children by local target
incidence; on a neutral block, \cref{lem:neutral-block} controls the entire
subtree.

\begin{theorem}[Neutral block Yomdin covering]
\label{thm:neutral-covering}
Fix $M,N\ge1$, and let $\cQ$ be a finite partition of $M$ whose atoms
have piecewise smooth boundaries and diameter less than $\rho_*\eps$.
Let $(\sigma,\theta)$ be a strongly $\eps$-bounded couple with $\theta$
injective on $(0,1)^k$, and set
\[
 \mathscr D_0:=\sigma([0,1]^k),
 \qquad D_0:=(\sigma\circ\theta)([0,1]^k).
\]
Let $\varnothing\ne F_N\subset D_0$.  For $x\in F_N$, apply the record
construction to the lifted point $(x,T_x\mathscr D_0)$.  Assume that
\[
 E=E_N^M\bigl((x,T_x\mathscr D_0)\bigr)
\]
is independent of $x\in F_N$, and that the $\cQ$-name
$(Q_j)_{j\in E}$ is also independent of $x$, meaning that
$g^jx\in Q_j$ for every $x\in F_N$ and $j\in E$.  Assume also that
there is a common sequence of nonnegative integers $c_0,\ldots,c_{N-1}$
satisfying
\begin{equation}\label{eq:cover-weight-word}
 w(g^jx)\le c_j
 \qquad(x\in F_N,\ 0\le j<N).
\end{equation}
The parameters $a,w,\Theta$ and the admissible scale are as fixed in
\eqref{eq:local-counting-setup}.

Then there is a finite family $\Psi(F_N)$ of terminal source maps
$\psi:[0,1]^k\to M$ whose images cover $F_N$ and satisfy
\begin{equation}\label{eq:bounded-orbit-chart}
 \|d(g^j\circ\psi)\|_\infty\le K_0,
 \qquad 0\le j\le N,
\end{equation}
and
\begin{equation}\label{eq:fiber-cover-main}
 \log\#\Psi(F_N)
 \le (N-|E|)\Theta+N\Delta_r(g)
       +\sum_{j=0}^{N-1}c_j+N\omega_M+C_{\rm end},
\end{equation}
where
\begin{equation}\label{eq:Delta-r}
 \Delta_r(g):=\log D_r(g),
\end{equation}
\begin{equation}\label{eq:omega-M}
 \omega_M:=\frac{\log C_\gamma}{M}
 \xrightarrow[M\to\infty]{}0,
\end{equation}
and
\begin{equation}\label{eq:C-end}
 C_{\rm end}:=\log C_\gamma+\log A_*.
\end{equation}
Here $C_\gamma$ is the constant defined in
\eqref{eq:neutral-block-proper}.  The cover contains every point,
including points on artificial algebraic boundaries.  For fixed $r$,
$\Theta$, and a common admissible scale $\eps$, the statement is uniform
in a $C^r$-bounded family for which the local counting setup holds at a
common scale and $\Theta-a-1\ge\gamma_0>0$.
\end{theorem}

\begin{proof}
Start the tree with the single couple $(\sigma,\theta)$ and
$A_o=D_0$.  
By
\cref{lem:exact-time-decomposition},
\[
 \{1,\ldots,N\}=E\sqcup\bigsqcup_{i=1}^q I_i,\qquad
 I_i=\{s_i+1,\ldots,s_i+\ell_i\},
\]
where
\begin{equation}\label{eq:number-neutral-blocks}
 q\le N/M+1,\qquad \sum_i\ell_i=N-|E|.
\end{equation}
We count active descendants across these consecutive pieces.

At a selected target time $j\in E$, every active child has target meeting
the prescribed atom $Q_j$.  Its diameter is $<\rho_*\eps$, so
\textnormal{(CP4)} gives at most $D$ active children of each parent.

At the entrance $s_i$ of a neutral block, restart at each active node
$v$ using \eqref{eq:tree-restart-data}.  The new carrier is
$\widetilde\sigma=g^{s_i}\circ\sigma$ and the new initial owned source
is $\widetilde A_o=R_v$.  Its active subset is
\[
 \widetilde F_v
 :=\bigsqcup_{\vartheta\in\Theta_v}
       g^{s_i}(F_N\cap A_{v,\vartheta})
 =g^{s_i}(F_N\cap A_v)\subset R_v.
\]
This set is nonempty.  By \textnormal{(CP1)}, the tangent space of the new
carrier at $g^{s_i}x$ is $dg_x^{s_i}(T_x\mathscr D_0)$.  Thus the lifted
orbit of $g^{s_i}x$ is the original lifted orbit starting at time $s_i$,
and \eqref{eq:anchored-neutral-prefix} gives every neutral prefix
inequality required by \cref{lem:neutral-block} on $\widetilde F_v$.
The weight word is $c_{s_i},\ldots,c_{s_i+\ell_i-1}$, including the
entrance source time.  The lemma, applied to this entire initial node
rather than separately counting its constituents as roots, bounds the
active descendants of $v$ by
\[
 C_\gamma D^{\ell_i}
 \exp\left(\ell_i\Theta+\sum_{j=s_i}^{s_i+\ell_i-1}c_j\right).
\]
By \cref{prop:yomdin-tree-restart}, the restarted tree is the original
descendant subtree with the inherited ordering.  An owned source set $B$
in the restart pulls back to $g^{-s_i}B\subset A_v$, preserving activity.
The parameter factorizations remain those of the full branch, so
\textnormal{(CP6)} holds on the entire time interval, not just separately
on its blocks.

Multiplying from the single root gives
\begin{equation}\label{eq:terminal-node-product}
\begin{split}
 \#\{\text{active terminal nodes}\}
 &\le D^{|E|}\prod_i
 \left[C_\gamma D^{\ell_i}
 e^{\ell_i\Theta+\sum_{j=s_i}^{s_i+\ell_i-1}c_j}\right]\\
 &\le D^N C_\gamma^q
       e^{(N-|E|)\Theta+\sum_{j=0}^{N-1}c_j}.
\end{split}
\end{equation}
The last inequality uses nonnegativity of the weights and disjointness of
the neutral source-time sets, including the final tail.  For each active
terminal node take the source maps $\psi=\sigma\circ\theta_{v,\vartheta}$
of all its constituents.  Nesting places their images in $D_0$;
\textnormal{(CP1)} makes these images cover every point of $F_N$, and
\textnormal{(CP6)} gives the orbit-derivative bound.
Each terminal node has at most $A_*$ constituents.  Taking logarithms and
using \eqref{eq:number-neutral-blocks} gives
\eqref{eq:fiber-cover-main}, with the stated $\omega_M$ and $C_{\rm end}$.
Uniformity follows from $\gamma\ge\gamma_0>0$.
\end{proof}

\section{Entropy decomposition for a fixed iterate}\label{sec:plaque-selection}

Let $g$ be a smooth diffeomorphism and let $\nu$ be an ergodic hyperbolic
invariant probability of unstable index $k$.  Write $d^u$ for the intrinsic
Riemannian distance on each full unstable manifold, and set
\begin{equation}\label{eq:intrinsic-Bowen-ball}
 B_N^u(x,\rho):=\{y\in W^u_g(x):
 d^u(g^jx,g^jy)<\rho\text{ for }0\le j<N\}.
\end{equation}
Choose a measurable partition $\xi$ subordinate to the full unstable
manifolds, with conditional probabilities $(\nu_x^\xi)$.  Ledrappier--Young's
conditional local-entropy theorem gives
\begin{equation}\label{eq:conditional-local-entropy-positive-set}
 \lim_{\rho\downarrow0}\liminf_{N\to\infty}
 -\frac1N\log\nu_x^\xi(B_N^u(x,\rho))=h_\nu(g)
 \quad\text{for }\nu\text{-almost every }x;
\end{equation}
see \cite[Proposition~7.2.1 and Corollary~7.2.2]{LEDRAPPIER_YOUNG_B}.
The equality of full unstable entropy and metric entropy is
\cite[Corollary~5.3]{LEDRAPPIER_YOUNG_A}.

\begin{theorem}[Joint upper semicontinuity]\label{thm:joint-usc}
Let $f_n\to f$ in $C^\infty$ and let $\nu_n\to\nu$ weakly, with
$\nu_n\in\mathcal M(f_n)$.  Then $\nu\in\mathcal M(f)$ and
\begin{equation}\label{eq:joint-usc}
\limsup_{n\to\infty}h_{\nu_n}(f_n)\le h_\nu(f).
\end{equation}
\end{theorem}
\begin{proof}  The invariance of measures passes to the limit by
the uniform convergence $f_n\rightrightarrows  f$. The joint upper semicontinuity of metric entropy  is obtained by Newhouse \cite[Theorem 1]{Newhouse} together with the uniform local  entropy  estimate of \cite[Theorem 2, Theorem 3]{Newhouse} and  Yomdin \cite{yomdin1987} in the convergence, see also Buzzi \cite{BuzziIsrael}. \end{proof} 

We combine the conditional entropy formula with the neutral covering.  The only
partition used in the limiting argument will be a fixed small partition
recorded at selected times.

\subsection{Entropy on one unstable plaque}

\begin{lemma}[Uniform set on an unstable plaque]\label{lem:plaque-preparation}
Let $g\in\Diff^\infty(M)$, let $G=F_g$, and let
$\nu\in\mathscr H_k(g)$.  Fix an integer $r\ge2$, a target scale
$\eps>0$, an error $\iota>0$, and a finite family
$\mathscr F\subset C(\widehat M_k)$.  There exist a $C^r$ unstable plaque
$D=\sigma([0,1]^k)$ with $(\sigma,\id)$ strongly $\eps$-bounded, a compact
set $F\subset D$ of regular points, a probability $\zeta$ obtained by
restricting one unstable conditional to $F$ and normalizing, and numbers $\rho>0$ and
$N_0\ge1$ such that, for all $x\in F$ and $N\ge N_0$,
\begin{equation}\label{eq:prepared-Bowen-mass}
 \zeta(B_N^u(x,\rho))\le e^{-N(h_\nu(g)-\iota)},
\end{equation}
and
\begin{equation}\label{eq:uniform-empirical}
 \left|\frac1N\sum_{j=0}^{N-1}\varphi(G^j(x,T_xD))
       -\int\varphi\,\dd\widehat\nu^u\right|<\iota
 \qquad(\varphi\in\mathscr F).
\end{equation}
\end{lemma}

\begin{proof}
By \eqref{eq:conditional-local-entropy-positive-set}, some $\rho>0$ and
$L\ge1$ give a positive-measure set $E$ on which
\[
 \nu_x^\xi(B_N^u(x,\rho))
 \le e^{-N(h_\nu(g)-\iota/2)}\qquad(N\ge L).
\]
Intersect $E$ with the full-measure regular and Birkhoff-generic sets for
$\mathscr F$, and choose an atom $\xi(x_*)$ on which this intersection has
positive conditional measure.  On this atom the conditionals agree almost
everywhere.  Egorov's theorem and inner regularity give a compact subset
$F_0$ of positive conditional measure on which the prescribed empirical
averages converge uniformly.

The containing unstable manifold is $C^r$ by
\cite[Theorem~6.3]{Ruelle1979}, applied to $g^{-1}$.
Work in a compact subchart of a local $C^r$ embedding $u$ whose image
meets $F_0$ in positive conditional measure.  An affine cube map $\ell$
of side $\delta$ in this subchart satisfies
\begin{gather*}
 \|u\circ\ell\|_{\mathrm{BO},r}=O(\delta),
 \qquad \inf|\wedge^k d(u\circ\ell)|\ge c\delta^k,\\
 [\wedge^k d(u\circ\ell)]_s=O(\delta^{k+s})
 \qquad(1\le s<r).
\end{gather*}
A countable cover by sufficiently small cubes compactly contained in the
chart therefore contains one strongly $\eps$-bounded plaque $D$ with
$\nu_{x_*}^\xi(F_0\cap D)>0$.  Set $F=F_0\cap D$ and
$\zeta=\nu_{x_*}^\xi|_F/\nu_{x_*}^\xi(F)$.  Restriction multiplies the
Bowen ball estimate by the fixed factor $\nu_{x_*}^\xi(F)^{-1}$, which is
absorbed by increasing $N_0$.  Uniform empirical convergence is preserved,
and $T_xD=E^u_\nu(x)$ for $x\in F$.
\end{proof}

For a finite
Borel partition $\cP$ and a probability $\zeta$, we use the standard notation
\[
 \Ent_\zeta(\cP):=-\sum_{P\in\cP}\zeta(P)\log\zeta(P),
\]
with $0\log0=0$; conditional partition entropy of $\cP$ relative to a finite
Borel partition $\cA$ is denoted by
\[\;
\Ent_\zeta(\cP\mid\cA)
:=-\sum_{A\in\cA}\sum_{P\in\cP}
\zeta(A\cap P)\log\frac{\zeta(A\cap P)}{\zeta(A)}
\;=\Ent_\zeta(\cP\vee\cA)-\Ent_\zeta(\cA),\]
where $ \cP\vee\cA=\{P\cap A:P\in\cP,\,A\in\cA\}$.

\begin{lemma}[Entropy from terminal plaque covers]
\label{lem:plaque-cover-entropy}
In the setting of \cref{lem:plaque-preparation},  let $N\ge N_0$.
Suppose a finite Borel partition $\mathcal B$ of $F$ has, for each nonempty
atom $B$, a cover by $L_B$ $C^1$ source maps $\psi:[0,1]^k\to D$ satisfying
\[
 \|d(g^j\circ\psi)\|_\infty\le K_0\qquad(0\le j<N).
\]
Then, for a constant $C=C(k,K_0,\rho)$ independent of $N$ and the covers,
\begin{equation}\label{eq:plaque-cover-entropy}
 h_\nu(g)-\iota
 \le\frac1N\Ent_\zeta(\mathcal B)
   +\frac1N\sum_{B\in\mathcal B}\zeta(B)\log L_B
   +\frac{\log C}{N}.
\end{equation}
\end{lemma}

\begin{proof}
Subdivide the parameter cube into $C$ closed cubes of diameter
$<\rho/K_0$.  For any such cube $Q$ and any $x\in F\cap\psi(Q)$,
\[
 F\cap\psi(Q)\subset B_N^u(x,\rho).
\]
  Thus \eqref{eq:prepared-Bowen-mass} gives
$\zeta(B)\le C L_B e^{-N(h_\nu(g)-\iota)}$.  Take logarithms for the
positive-mass atoms, multiply by $\zeta(B)$, and sum.
\end{proof}

\begin{lemma}[Regular small partitions, Lemma 8.5 of \cite{Walters}]\label{lem:regular-small-partition}
Let $\lambda$ be a probability on $M$ and let $\delta>0$.  There is a
finite partition $\cQ$ with piecewise smooth boundaries such that
\begin{equation}\label{eq:regular-small-partition}
 \operatorname{mesh}\cQ<\delta,
 \qquad \lambda(\partial\cQ)=0,
\end{equation}
and every connected subset of $M\setminus\partial\cQ$ lies in one atom.
If $\lambda$ is $g$-invariant, then
\[
 \lambda\!\left(\bigcup_{j\in\mathbb Z}g^{-j}\partial\cQ\right)=0.
\]
\end{lemma}

\subsection{Finite-partition entropy and hyperbolic components} 
We next isolate the part of an empirical measure carried by selected times.
The finite-set entropy estimate below converts the entropy of a selected
name into the entropy of an averaged measure, while the subsequent
proposition proves that the limiting selected part is invariant and lies on
backward hyperbolic segments.

For a finite partition $\cQ$ and a finite set
$E\subset\{0,\ldots,N\}$, write
\[
 \cQ_g^E:=\bigvee_{j\in E}g^{-j}\cQ,
 \qquad
 \cQ_g^m:=\bigvee_{i=0}^{m-1}g^{-i}\cQ.
\]

\begin{lemma}[Entropy along a finite set]\label{lem:interval-entropy} 
Let $E\subset\mathbb N$ be finite and nonempty, put
\[
 \zeta^E:=\frac1{|E|}\sum_{j\in E}g_*^j\zeta,
 \qquad E+1:=\{j+1:j\in E\},
\]
and let $m\ge1$.  Then
\begin{equation}\label{eq:interval-entropy}
 \frac1{|E|}\Ent_\zeta(\cQ_g^E)
 \le\frac1m\Ent_{\zeta^E}(\cQ_g^m)
 +6m\frac{|(E+1)\triangle E|}{|E|}\log\#\cQ.
\end{equation}
In particular, if $E$ is a union of at most $q$ integer intervals, then the
last term after multiplication by $|E|$ is at most
$12mq\log\#\cQ$.
\end{lemma}

\begin{proof}
The displayed estimate is exactly \cite[Section~2.4, Lemma~8]{Burguet}, applied to the
topological system $(M,g)$, the probability $\zeta$, and the partition
$\cQ$.  If $E$ is a union of $q$ integer intervals, each interval contributes
at most its two endpoints to $(E+1)\triangle E$, so
$|(E+1)\triangle E|\le2q$.
\end{proof}

\begin{lemma}[Finite-time partition entropy for varying maps]
\label{lem:varying-partition-entropy}
Let $g_n\to g$ uniformly, $\rho_n\to\rho$ weakly, and fix $m\ge1$. If
\begin{equation}\label{eq:Q-boundary-null}
 \rho\left(\bigcup_{j=0}^{m-1}g^{-j}\partial\cQ\right)=0,
\end{equation}
then
\begin{equation}\label{eq:varying-static-entropy}
 \Ent_{\rho_n}(\cQ_{g_n}^m)\longrightarrow
 \Ent_\rho(\cQ_g^m).
\end{equation}
\end{lemma}

\begin{proof}
 This is a standard consequence of weak convergence on continuity sets (boundary of measure zero) and of
the continuity of Shannon entropy on a finite probability simplex.
\end{proof}

\bigskip
\paragraph{\it Hyperbolic empirical components.}
For fixed $M$, the measure below keeps only the orbit mass carried by the
selected times.  Its total mass is the selected-time frequency; passing first
to $n\to\infty$ and then to $M\to\infty$ produces the invariant hyperbolic
component used in the entropy decomposition. 

\begin{proposition}[Hyperbolic empirical component]\label{prop:selected-package}
Assume $g_n^{\pm1}\to g^{\pm1}$ in $C^1$, write $G_n=F_{g_n}$ and $G=F_g$, and let
$N_n\to\infty$.  Let $\widehat\zeta_n$ be probabilities on $\widehat M_k$ and
assume
\begin{equation}\label{eq:prepared-full-empirical}
 \widehat e_n:=\frac1{N_n}\int\sum_{j=1}^{N_n}
 \delta_{G_n^j\widehat x}\,\dd\widehat\zeta_n(\widehat x)
 \longrightarrow\widehat\eta,
 \qquad \mu:=\pi_*\widehat\eta.
\end{equation}
Use the record construction with the same threshold $\Theta$ for every $n$,
and define, for fixed $M$,
\begin{equation}\label{eq:selected-submeasure}
 \widehat\xi_n^M:=\frac1{N_n}
 \int\sum_{j\in E_{N_n}^M(\widehat x)}
 \delta_{G_n^j\widehat x}\,\dd\widehat\zeta_n(\widehat x).
\end{equation}
After passage to one subsequence, there are finite measures
$\widehat\xi^M$ such that
\begin{equation}\label{eq:selected-limits}
 \widehat\xi_n^M\longrightarrow\widehat\xi^M
 \quad(M\ge1),
 \qquad
 0\le\widehat\xi^M\le\widehat\xi^{M+1}\le\widehat\eta.
\end{equation}
Put
\begin{equation}\label{eq:selected-limit-infinite}
 \widehat\xi^\infty:=\lim_{M\to\infty}\widehat\xi^M,
 \qquad
 \beta:=\widehat\xi^\infty(\widehat M_k).
\end{equation}
The measure $\widehat\xi^\infty$ is $G$-invariant.  If $\beta>0$, define
\begin{equation}\label{eq:mu1-definition}
 \widehat\mu_1:=\beta^{-1}\widehat\xi^\infty,
 \qquad
 \mu_1:=\pi_*\widehat\mu_1.
\end{equation}
Then $\widehat\mu_1$ and $\mu_1$ are $G$- and $g$-invariant, respectively,
\begin{equation}\label{eq:selected-domination}
 \beta\widehat\mu_1\le\widehat\eta,
 \qquad
 \beta\mu_1\le\mu,
\end{equation}
and
\begin{equation}\label{eq:selected-support}
 \widehat\mu_1\left(\bigcup_{s\ge0}G^{-s}H_\Theta(g)\right)=1.
\end{equation}
\end{proposition}

\begin{proof}
\medskip\noindent\emph{Step 1: simultaneous limits in $M$.}
A diagonal extraction over the countable parameter $M$ gives
$\widehat\xi_n^M\to\widehat\xi^M$ for every $M$.  Since
$E_{N_n}^M\subset E_{N_n}^{M+1}$,
\[
 0\le\widehat\xi_n^M\le\widehat\xi_n^{M+1}\le\widehat e_n;
\]
passing to the limit proves \eqref{eq:selected-limits}.  Hence
$\widehat\xi^M\nearrow\widehat\xi^\infty\le\widehat\eta$.

\medskip\noindent\emph{Step 2: invariance of the increasing limit.}
Fix $\varphi\in C(\widehat M_k)$ with $\|\varphi\|_\infty\le1$.  For
$\widehat x$, write the maximal interval decomposition as
\[
 E_{N_n}^M(\widehat x)=
 \bigsqcup_{h=1}^{q_n(\widehat x)}
 \{u_h(\widehat x),\ldots,v_h(\widehat x)\},
 \qquad q_n(\widehat x)\le N_n/M+2.
\]
When $E_{N_n}^M(\widehat x)=\varnothing$, put
$q_n(\widehat x)=0$ and interpret the sums below as zero.  For every
$\widehat x$ one has the exact  identity
\begin{align}
 \sum_{j\in E_{N_n}^M(\widehat x)}
 \bigl(\varphi(G_n^{j+1}\widehat x)-
       \varphi(G_n^j\widehat x)\bigr)=
 \sum_{h=1}^{q_n(\widehat x)}
 \bigl(\varphi(G_n^{v_h(\widehat x)+1}\widehat x)-
       \varphi(G_n^{u_h(\widehat x)}\widehat x)\bigr).
 \label{eq:selected-telescoping}
\end{align}
By \eqref{eq:selected-submeasure},
\begin{align*}
 &\int(\varphi\circ G_n-\varphi)\,\dd\widehat\xi_n^M\\
 &=\frac1{N_n}\int
   \sum_{j\in E_{N_n}^M(\widehat x)}
   \bigl(\varphi(G_n^{j+1}\widehat x)-
         \varphi(G_n^j\widehat x)\bigr)
   \,\dd\widehat\zeta_n(\widehat x).
\end{align*}
Insert \eqref{eq:selected-telescoping}.  Since each summand on its right has
absolute value at most $2$,
\begin{equation}\label{eq:approx-invariance}
 \left|\int\varphi\circ G_n\,\dd\widehat\xi_n^M
       -\int\varphi\,\dd\widehat\xi_n^M\right|
 \le\frac{2}{N_n}\left(\frac{N_n}{M}+2\right)
 =\frac2M+\frac4{N_n}.
\end{equation}
For fixed $M$, let $n\to\infty$.  Uniform convergence $G_n\to G$ and
weak convergence give
\[
 \left|\int\varphi\circ G\,\dd\widehat\xi^M
       -\int\varphi\,\dd\widehat\xi^M\right|\le\frac2M.
\]
Since $\widehat\xi^M\nearrow\widehat\xi^\infty$, for every bounded
continuous $\psi$,
\[
 \left|\int\psi\,\dd(\widehat\xi^\infty-\widehat\xi^M)\right|
 \le\|\psi\|_\infty
 \bigl(\widehat\xi^\infty(\widehat M_k)-
       \widehat\xi^M(\widehat M_k)\bigr)\longrightarrow0.
\]
Letting $M\to\infty$ in the preceding invariance defect yields
$G_*\widehat\xi^\infty=\widehat\xi^\infty$.  When $\beta>0$, normalization
and projection give invariance of $\widehat\mu_1$ and $\mu_1$;
$\widehat\xi^\infty\le\widehat\eta$ then gives
\eqref{eq:selected-domination}.

\medskip\noindent\emph{Step 3: concentration on backward
hyperbolic segments.}
Fix $R,M\ge1$.  By \cref{lem:selected-time-sets}, each selected time has
an associated record $b$ with
\begin{equation}\label{eq:selected-record-association}
 j=b-s,\qquad 0\le s<M.
\end{equation}
At most $R-1$ selected times have $b<R$.  For every remaining time,
\eqref{eq:finite-hyperbolic-past} gives
$G_n^b\widehat x\in H_{\Theta,R}(g_n)$.  Thus $\widehat\xi_n^M$ gives
mass at most $(R-1)/N_n$ to the complement of
\[
 K_{n,R,M}:=\bigcup_{s=0}^{M-1}G_n^{-s}H_{\Theta,R}(g_n).
\]
Set $K_{R,M}:=\bigcup_{s=0}^{M-1}G^{-s}H_{\Theta,R}(g)$.
If $y_n\in K_{n,R,M}$ converges to $y$, pass to a subsequence on which
the associated shift $s<M$ is fixed.  Uniform convergence of the finitely
many maps and potentials in the defining inequalities gives
$G^sy\in H_{\Theta,R}(g)$.  Hence
\begin{equation}\label{eq:moving-support-inclusion}
 \limsup_{n\to\infty}K_{n,R,M}\subset K_{R,M}.
\end{equation}
If $K_{R,M}$ is empty, compactness makes $K_{n,R,M}$ eventually empty,
and $\widehat\xi^M=0$.  Otherwise, for each $t>0$ the closed
$t$-neighborhood of $K_{R,M}$ contains every sufficiently late
$K_{n,R,M}$.  Portmanteau's theorem gives full $\widehat\xi^M$-mass to
that neighborhood.  Letting $t\downarrow0$ shows that $K_{R,M}$ carries
$\widehat\xi^M$.

The sets $H_{\Theta,R}(g)$ decrease with $R$.  Since there are only
finitely many shifts,
\begin{equation}\label{eq:finite-union-intersection}
 \bigcap_{R\ge1}K_{R,M}
 =\bigcup_{s=0}^{M-1}G^{-s}H_\Theta(g).
\end{equation}
Indeed, for a point in the left side, some shift occurs for arbitrarily
large $R$ and hence satisfies every defining inequality.  Continuity from
above now gives full $\widehat\xi^M$-mass to the right side.  Finally,
$\widehat\xi^M\nearrow\widehat\xi^\infty$, so normalization yields
\eqref{eq:selected-support}.
\end{proof}

\subsection{Identification of the unstable lift}

The selected component lies on backward-expanding $k$-planes.  Since its
projection has unstable dimension $k$, the two-sided Oseledets theorem
identifies these planes with the unstable spaces.  We first record how
invariant domination passes to ergodic decompositions and graph measures.

\begin{lemma}[Invariant domination and unstable lifts]
\label{lem:dominated-lifts}
Let $\alpha$ and $\mu$ be finite $f$-invariant measures with
$0\le\alpha\le\mu$, and let $\tau_\alpha$ and $\tau_\mu$ be the finite
measures on $\mathcal M_e(f)$ in their ergodic decompositions.  Then
\begin{equation}\label{eq:dominated-ergodic-decompositions}
 \tau_\alpha\le\tau_\mu.
\end{equation}
In particular, every property holding for $\tau_\mu$-almost every ergodic
measure holds for $\tau_\alpha$-almost every ergodic measure.  If almost every ergodic component of $\mu$ is hyperbolic with
exactly $k$ positive Lyapunov exponents, define the unstable lifts as finite
graph measures with the same masses as their base measures.  Then
\begin{equation}\label{eq:dominated-lifts}
 \widehat\alpha^u\le\widehat\mu^u.
\end{equation}
\end{lemma}

\begin{proof}
The Radon--Nikodym derivative $d\alpha/d\mu$ is invariant and takes values
in $[0,1]$.  It is therefore constant on almost every ergodic component.
Uniqueness of ergodic decomposition gives
$d\tau_\alpha=c\,d\tau_\mu$ with $0\le c\le1$.  Applying the common
unstable graph map to $\alpha\le\mu$ proves the lift inequality.
\end{proof}

\begin{proposition}[Identification of the hyperbolic component]
\label{prop:poincare-selection}
Assume $\beta>0$, almost every ergodic component of $\mu$ is hyperbolic
with exactly $k$ positive Lyapunov exponents, and $A_{k-1}(g)<\Theta$.  Then the measure $\wh\mu_1$ in
\eqref{eq:mu1-definition} is the unstable lift of its projection:
\begin{equation}\label{eq:mu1-unstable}
 \wh\mu_1=\wh\mu_1^u.
\end{equation}
Consequently
\begin{equation}\label{eq:unstable-domination}
 \beta\wh\mu_1^u\le\wh\eta.
\end{equation}
\end{proposition}

\begin{proof}
Let $\mathcal R$ be an invariant full $\mu$-measure set of two-sided
Oseledets-regular points with unstable dimension $k$.  Since
$\beta\mu_1\le\mu$, it also has full $\mu_1$-measure.
For $\widehat\mu_1$-almost every $\widehat x=(x,E)$, we have
$x\in\mathcal R$ and, by \eqref{eq:selected-support},
$\widehat y=G^s\widehat x=(g^sx,dg_x^sE)\in H_\Theta(g)$ for some
$s\ge0$.  Set $E_y=dg_x^sE$ and $\gamma=\Theta-A_{k-1}(g)>0$.
By \cref{lem:conorm-expansion},
\[
 \|dg_y^{-m}|_{E_y}\|\le e^{-m\gamma}\qquad(m\ge1).
\]
At a two-sided regular point, the vectors with negative backward
Lyapunov exponent are precisely those in the positive Oseledets space
\cite{Oseledets}.  Hence $E_y\subset E^u_\mu(y)$, and the common dimension
$k$ gives equality.  Pulling back by $dg_x^s$ yields $E=E^u_\mu(x)$.
Thus $\widehat\mu_1$ is the unstable graph measure of its projection, and
\eqref{eq:selected-domination} gives \eqref{eq:unstable-domination}. 
\end{proof}

\subsection{The fixed-iterate entropy decomposition}

The next theorem is the main technical output of the preceding construction.
It decomposes the limiting entropy into an invariant hyperbolic part and a
neutral remainder.  Besides the loss of finite smoothness by Yomdin estimates, it retains
an averaged localization error.  The latter permits a single exterior
growth threshold in the large-iterate argument. Recall  $\ell_k(g,x):=\log\max_{0\le q<k}\|\wedge^q dg_x\|$.

\begin{theorem}[Fixed-iterate entropy decomposition]
\label{thm:fixed-iterate-decomposition}
Fix an integer $r\ge2$.  Let $g_n,g\in\Diff^\infty(M)$, and write
$G_n:=F_{g_n}$ and $G:=F_g$.  Assume the following.
\begin{enumerate}[label=\textnormal{(\roman*)},leftmargin=10mm]
\item The maps and their inverses converge in $C^r$:
\[
 g_n^{\pm1}\longrightarrow g^{\pm1}.
\]
\item Each $\nu_n$ is an ergodic hyperbolic $g_n$-invariant probability
with exactly $k$ positive Lyapunov exponents, and
\begin{equation}\label{eq:fixed-thm-limits}
 \nu_n\longrightarrow\mu,
 \qquad
 \wh\nu_n^u\longrightarrow\wh\eta.
\end{equation}
Almost every ergodic component of $\mu$ is hyperbolic with exactly $k$ positive
Lyapunov exponents.
\item There is $a\in\mathbb R$ such that, for all
sufficiently large $n$,
\begin{equation}\label{eq:fixed-thm-margin}
 A_{k-1}(g_n)\le a.
\end{equation}
\end{enumerate}
Discard these finitely many indices.  With $t^+:=\max\{t,0\}$, set
\begin{equation}\label{eq:localization-weights}
 w_n(x):=(\ell_k(g_n,x)-a)^+,
 \qquad w(x):=(\ell_k(g,x)-a)^+.
\end{equation}
Choose an integer $R\ge2$ such that
\begin{equation}\label{eq:weight-alphabet-bound}
 1+\left\lceil\sup_n\|w_n\|_\infty\right\rceil\le R.
\end{equation}
Then for  $\Theta>a+1$ there exist $\beta\in[0,1]$ and, when $\beta>0$, a $g$-invariant
probability $\mu_1$ with the following properties:
\begin{equation}\label{eq:fixed-thm-domination}
 \beta\mu_1\le\mu,
 \qquad
 \beta\wh\mu_1^u\le\wh\eta.
\end{equation}
Moreover,
\begin{equation}\label{eq:fixed-thm-entropy}
\begin{split}
 \limsup_{n\to\infty}h_{\nu_n}(g_n)
 \le{}&\beta h_{\mu_1}(g)+(1-\beta)\Theta+\Delta_r(g)\\
      &+\int w\,\dd\mu+1+\log R.
\end{split}
\end{equation}
Here
\begin{equation}\label{eq:fixed-Delta}
 \Delta_r(g):=\log C_{\rm Y}
 +\frac{3k^2}{r-1}\log\mathfrak M(g).
\end{equation}
When $\beta=0$, the expressions $\beta\mu_1$,
$\beta\wh\mu_1^u$, and $\beta h_{\mu_1}(g)$ are understood to be zero.
The constant $C_{\rm Y}=C_{M,r,d,k}$ is the complexity constant from
\cref{thm:closed-yomdin-tree}; it is independent of $g$ and, in particular,
of an iterate $g=f^p$.
\end{theorem}

\begin{proof}
\medskip\noindent\emph{Step 1: plaque preparation and the limiting
empirical lift.}
Pass to a subsequence realizing the limsup.  Finite-time $C^1$ continuity
gives $w_n\to w$ uniformly and
\begin{equation}\label{eq:limit-L-margin}
 A_{k-1}(g)\le a<\Theta.
\end{equation}
In particular the family $(w_n)$ is equicontinuous.
Uniform convergence and invariance give
\begin{equation}\label{eq:fixed-invariant-limits}
 g_*\mu=\mu,
 \qquad G_*\widehat\eta=\widehat\eta,
 \qquad \pi_*\widehat\eta=\mu.
\end{equation}
Choose a common scale
$0<\eps\le\eps_*$ for the $C^r$-bounded family $\{g,g_n:n\ge n_0\}$,
small enough that each $w_n$ has oscillation at most $1$ on every raw
terminal constituent image.  Thus \eqref{eq:local-counting-setup} holds uniformly with
$a,w_n,\Theta$.
Let $(\varphi_j)$ be dense in the unit ball of $C(\widehat M_k)$, put
$\mathscr F_n=\{\varphi_1,\ldots,\varphi_n,w_n\circ\pi\}$, and choose
$\eta_n\downarrow0$.  Apply \cref{lem:plaque-preparation} to obtain
$D_n,F_n,\zeta_n,\rho_n,N_{0,n}$ with error $\eta_n$ and target scale
$\eps$.  If $C_n=C(k,K_0,\rho_n)$ is the constant from
\cref{lem:plaque-cover-entropy}, choose
\[
 N_n\ge\max\{n,N_{0,n}\},\qquad \frac{\log C_n}{N_n}<\eta_n.
\]
All prescribed lifted empirical averages are then $\eta_n$-close to
$\widehat\nu_n^u$ uniformly on $F_n$.  In particular,
\begin{equation}\label{eq:prepared-weight-average}
 \frac1{N_n}\sum_{j=0}^{N_n-1}w_n(g_n^jx)
 \le\int w_n\,\dd\nu_n+\eta_n
 \qquad(x\in F_n).
\end{equation}

Lift $\zeta_n$ by $x\mapsto(x,T_xD_n)$ and put
\[
 \widehat e_n=\frac1{N_n}\int\sum_{j=1}^{N_n}
 \delta_{G_n^j\widehat x}\,\dd\widehat\zeta_n(\widehat x).
\]
For every fixed $j$ and every $n\ge j$, the plaque preparation gives
\begin{equation}\label{eq:prepared-empirical-test}
 \left|\int\varphi_j\,\dd\widehat e_n
       -\int\varphi_j\,\dd\widehat\nu_n^u\right|
 \le\eta_n+\frac{2\|\varphi_j\|_\infty}{N_n}.
\end{equation}
Since $(\varphi_j)$ is dense in the unit ball of
$C(\widehat M_k)$ and $\widehat\nu_n^u\to\widehat\eta$,
\eqref{eq:prepared-empirical-test} implies
$\widehat e_n\to\widehat\eta$.  After passing to the subsequence supplied by
\cref{prop:selected-package}, write
\[
 \widehat\xi_n^M\to\widehat\xi^M,
 \qquad
 \widehat\xi^M\nearrow\widehat\xi^\infty,
 \qquad
 \beta:=\widehat\xi^\infty(\widehat M_k).
\]
When $\beta>0$, we write
$\widehat\xi^\infty=\beta\widehat\mu_1$.  Set
\begin{equation}\label{eq:selected-projections}
 \xi_n^M:=\pi_*\widehat\xi_n^M,
 \qquad
 \xi^M:=\pi_*\widehat\xi^M.
\end{equation}
Choose once and for all a partition $\cQ$ from
\cref{lem:regular-small-partition} with
$(\lambda,\delta)=(\mu,\rho_*\eps)$.  Since $\mu$ is $g$-invariant,
\begin{equation}\label{eq:Q-null-mu}
 \mu\left(\bigcup_{j\in\mathbb Z}g^{-j}\partial\cQ\right)=0.
\end{equation}

\medskip\noindent\emph{Step 2: the selected-name
partition.}
Fix $M$ and put $\gamma_0:=\Theta-a-1>0$.  The constants in
\cref{thm:closed-yomdin-tree,prop:yomdin-tree-restart,thm:neutral-covering}
are uniform at the chosen scale.  In particular, with
\begin{equation}\label{eq:uniform-block-errors}
 C_{\gamma_0}:=1+\frac{C_1}{1-e^{-\gamma_0}},\qquad
 \omega_M:=\frac{\log C_{\gamma_0}}M,\qquad
 C_{\rm end}^*:=\log C_{\gamma_0}+\log A_*,
\end{equation}
we have
\begin{equation}\label{eq:uniform-block-error-bounds}
 \omega_M\longrightarrow0\quad(M\to\infty),
 \qquad C_{\rm end}^*<\infty.
\end{equation}
Let $\mathcal E_n^M$ be the partition of $F_n$ by the selected set
$E_{N_n}^M$.  Since $N_n\to\infty$, after discarding finitely many
indices (depending on the fixed $M$), we have $N_n\ge N_0(M)$.  Therefore
\eqref{eq:selected-set-count} gives
\begin{equation}\label{eq:selected-set-entropy}
 \Ent_{\zeta_n}(\mathcal E_n^M)
 \le\log\#\mathcal E_n^M
 \le N_n\kappa_M.
\end{equation}
For an atom $C$ of $\mathcal E_n^M$, let $E(C)$ be its selected set and
write
\[
 \cQ_{g_n}^{E(C)}=\bigvee_{j\in E(C)}g_n^{-j}\cQ.
\]
Define the finite partition
\begin{equation}\label{eq:selected-name-partition}
 \cA_n^M:=\left\{C\cap Q:
 C\in\mathcal E_n^M,\ Q\in\cQ_{g_n}^{E(C)},\
 C\cap Q\ne\varnothing\right\}.
\end{equation}
Thus an atom of $\cA_n^M$ fixes both the selected set and its $\cQ$-name.
Let $\mathcal C_n$ be the finite Borel partition of $F_n$ determined by the
source-time word
\begin{equation}\label{eq:quantized-weight-word}
 c_j(x):=\lceil w_n(g_n^jx)\rceil,\qquad 0\le j<N_n.
\end{equation}
There are at most $R^{N_n}$ such words.  Set
$\mathcal B_n^M:=\cA_n^M\vee\mathcal C_n$.
The entropy chain rule gives
\begin{equation}\label{eq:entropy-three-terms}
 \Ent_{\zeta_n}(\mathcal B_n^M)
 \le \Ent_{\zeta_n}(\mathcal E_n^M)
    +\Ent_{\zeta_n}(\cA_n^M\mid\mathcal E_n^M)+N_n\log R.
\end{equation}

\medskip\noindent\emph{Step 3: entropy generated on neutral blocks.}
For each atom $B$ of $\mathcal B_n^M$, use the terminal source cover of
\cref{thm:neutral-covering}, and denote its cardinality by $L_B$.  It gives
\begin{equation}\label{eq:cover-size-on-B}
 \log L_B\le (N_n-|E(B)|)\Theta+N_n\Delta_r(g_n)
   +\sum_{j=0}^{N_n-1}c_j(B)+N_n\omega_M+C_{\rm end}^*.
\end{equation}

Put
\begin{equation}\label{eq:beta-nM}
 \beta_n^M:=\frac1{N_n}\int|E_{N_n}^M(\widehat x)|
                      \,\dd\widehat\zeta_n(\widehat x)
             =\widehat\xi_n^M(\widehat M_k).
\end{equation}
Average \eqref{eq:cover-size-on-B} over $\mathcal B_n^M$ and use
$c_j\le w_n(g_n^jx)+1$, \eqref{eq:prepared-weight-average}, and
\cref{lem:plaque-cover-entropy}.  Together with
\eqref{eq:selected-set-entropy} and \eqref{eq:entropy-three-terms}, this gives
\begin{equation}\label{eq:entropy-before-selected}
\begin{split}
 h_{\nu_n}(g_n)\le{}&
 \frac1{N_n}\Ent_{\zeta_n}(\cA_n^M\mid\mathcal E_n^M)
 +(1-\beta_n^M)\Theta+\Delta_r(g_n)\\
 &+\int w_n\,\dd\nu_n+1+\log R
   +\omega_M+\kappa_M+3\eta_n+\frac{C_{\rm end}^*}{N_n}.
\end{split}
\end{equation}

\medskip\noindent\emph{Step 4: entropy of the selected names.}
For an atom $C$ of $\mathcal E_n^M$ with $p_C:=\zeta_n(C)>0$, let
$\zeta_{n,C}=\zeta_n|_C/p_C$ and $E=E(C)$.   Definition
\eqref{eq:selected-name-partition} gives the exact identity
\begin{equation}\label{eq:selected-conditional-entropy-sum}
 \Ent_{\zeta_n}(\cA_n^M\mid\mathcal E_n^M)
 =\sum_{C\in\mathcal E_n^M}p_C
   \Ent_{\zeta_{n,C}}(\cQ_{g_n}^{E(C)}).
\end{equation}
 Put
\[
 \zeta_{n,C}^E:=
 \begin{cases}
 |E|^{-1}\displaystyle\sum_{j\in E}(g_n^j)_*\zeta_{n,C},
      &E\ne\varnothing,\\[1mm]
 \zeta_{n,C},&E=\varnothing.
 \end{cases}
\]
The second choice is immaterial because its coefficient $|E|$ is zero.
Since $E$ is a union of at most $N_n/M+2$ intervals,
\cref{lem:interval-entropy} gives the following estimate when
$E\ne\varnothing$; when $E=\varnothing$, the same estimate is immediate
because the left-hand side is zero:
\begin{equation}\label{eq:selected-interval-C}
 \Ent_{\zeta_{n,C}}(\cQ_{g_n}^{E})
 \le\frac{|E|}{m}\Ent_{\zeta_{n,C}^{E}}(\cQ_{g_n}^m)
 +12m(N_n/M+2)\log\#\cQ.
\end{equation}
From the definition of the conditional probabilities and of
$\widehat\xi_n^M$,
\begin{align}
 \sum_{C\in\mathcal E_n^M}
 p_C\frac{|E(C)|}{N_n}\zeta_{n,C}^{E(C)}
 &=\frac1{N_n}\int
   \sum_{j\in E_{N_n}^M(\widehat x)}
   \delta_{g_n^j x}\,\dd\zeta_n(x)
 =\xi_n^M, \label{eq:selected-convex-combination}\\
 \sum_Cp_C\frac{|E(C)|}{N_n}&=\beta_n^M. \label{eq:selected-weight-sum}
\end{align}
  We obtain
\begin{equation}\label{eq:selected-static}
 \frac1{N_n}\Ent_{\zeta_n}(\cA_n^M\mid\mathcal E_n^M)
 \le\frac{\beta_n^M}{m}
 \Ent_{\xi_n^M/\beta_n^M}(\cQ_{g_n}^m)
 +\frac{12m}{M}\log\#\cQ+\frac{24m}{N_n}\log\#\cQ,
\end{equation}
with the first term equal to zero when $\beta_n^M=0$. Here, if $\beta_n^M>0$, concavity of the entropy of the fixed partition
$\cQ_{g_n}^m$ applied to the probability weights
$p_C|E(C)|/(N_n\beta_n^M)$ gives the first term above.  If
$\beta_n^M=0$, then $|E(C)|=0$ for $p_C$-almost every $C$, and that term
is defined to be zero.

\medskip\noindent\emph{Step 5: passage to the limits.}
Let $n\to\infty$ in
\eqref{eq:entropy-before-selected} and \eqref{eq:selected-static}, and put
$\beta_M=\widehat\xi^M(\widehat M_k)$.  If $\beta_M>0$, then
$\xi_n^M/\beta_n^M\to\xi^M/\beta_M$.  Moreover,
$\xi^M\le\mu$ by \eqref{eq:selected-limits}, so
\eqref{eq:Q-null-mu} verifies the boundary-null hypothesis of
\cref{lem:varying-partition-entropy}.  If $\beta_M=0$, the corresponding
term is bounded by
$\beta_n^M\log\#\cQ$ and tends to zero; the first term in the next
display is then defined to be zero.  Since
$\Delta_r(g_n)\to\Delta_r(g)$ and
$\int w_n\,\dd\nu_n\to\int w\,\dd\mu$,
\begin{align}
 \limsup_n h_{\nu_n}(g_n)
 \le{}&\frac{\beta_M}{m}
 \Ent_{\xi^M/\beta_M}(\cQ_g^m)
 +(1-\beta_M)\Theta+\Delta_r(g)
 \notag\\
 &+\int w\,\dd\mu+1+\log R
 +\omega_M+\kappa_M+\frac{12m}{M}\log\#\cQ.
 \label{eq:entropy-before-limits}
\end{align}

Fix $m$ and let $M\to\infty$.  If $\beta>0$, then
$\beta_M>0$ for all sufficiently large $M$.  The finite measure
$\widehat\xi^M/\beta$ is dominated by both
$\widehat\xi^M/\beta_M$ and $\widehat\mu_1$, and has mass
$\beta_M/\beta$.  Hence
\begin{equation}\label{eq:selected-normalization-TV}
 \left\|\frac{\widehat\xi^M}{\beta_M}
       -\widehat\mu_1\right\|_{\TV}
 \le2\left(1-\frac{\beta_M}{\beta}\right)
 \longrightarrow0,
\end{equation}
where the total variation norm
$\|\sigma\|_{\TV}=|\sigma|(X)$ for a measure $\sigma$ on a space $X$.  Projection does not increase total variation.  Since Shannon entropy is
continuous on the finite probability simplex by Lemma \ref{lem:varying-partition-entropy},
\[
 \Ent_{\xi^M/\beta_M}(\cQ_g^m)
 \longrightarrow\Ent_{\mu_1}(\cQ_g^m).
\]
Thus, letting $M\to\infty$, \begin{equation}\label{eq:entropy-before-m}
\begin{split}
 \limsup_n h_{\nu_n}(g_n)
 \le{}&\frac{\beta}{m}\Ent_{\mu_1}(\cQ_g^m)
       +(1-\beta)\Theta+\Delta_r(g)\\
      &+\int w\,\dd\mu+1+\log R.
\end{split}
\end{equation}
The  invariance of $\mu_1$ gives
\[
 \lim_{m\to\infty}
 \frac{\beta}{m}\Ent_{\mu_1}(\cQ_g^m)
 =\beta h_{\mu_1}(g,\cQ)\le\beta h_{\mu_1}(g).
\]
This proves \eqref{eq:fixed-thm-entropy}. Finally, $\widehat \mu_1=\widehat \mu_1^u$ by Proposition \ref{prop:poincare-selection}.
\end{proof}


\section{Entropic continuity}\label{sec:ergodic-continuity}

To return from $f_n^p$ to $f_n$, we select one cyclic component, apply
\cref{thm:fixed-iterate-decomposition}, and average the base measure and its
unstable lift over the cycle.  The component is chosen to minimize the
localization error.  This error vanishes in the large-iterate limit, while
the finite-smoothness loss is removed by choosing $r$ before $p$.

\subsection{Cyclic decompositions and averaged lifts}

\begin{lemma}[Cyclic decomposition for an iterate]
\label{lem:cyclic-components}
Let $T:X\to X$ be an invertible measurable transformation of a standard
Borel space, and let $\nu$ be an ergodic $T$-invariant probability.
For $p\ge1$, the $T^p$-ergodic decomposition of $\nu$ consists of a finite
cycle $\nu_0,\ldots,\nu_{q-1}$, where $q\mid p$ and $T_*\nu_i=\nu_{i+1}$
cyclically.  In particular
\begin{equation}\label{eq:cyclic-base}
 \nu=\frac1p\sum_{j=0}^{p-1}T_*^j\nu_0
\end{equation}
and
\begin{equation}\label{eq:cyclic-entropy}
 h_{\nu_0}(T^p)=p h_\nu(T).
\end{equation}
If $T$ is a diffeomorphism and $\nu$ is hyperbolic with exactly $k$
positive Lyapunov exponents, then
\begin{equation}\label{eq:cyclic-lift}
 \wh\nu^u=\frac1p\sum_{j=0}^{p-1}F_{T*}^j\wh\nu_0^u,
\end{equation}
where $\wh\nu_0^u$ denotes the unstable lift of the
$T^p$-invariant measure $\nu_0$.
\end{lemma}

\begin{proof}
The cyclic structure of ergodic components under iteration follows from the ergodic decomposition theorem. Moreover,  the entropy power rule gives
$h_{\nu_i}(T^p)=h_\nu(T^p)=p h_\nu(T)$.
The splitting of $T^p$ is that of $T$, with all exponents multiplied by $p$,
so applying the unstable graph map to \eqref{eq:cyclic-base} gives
\eqref{eq:cyclic-lift}.
\end{proof}

\begin{lemma}[Entropy of a cyclic average]\label{lem:cyclic-average-entropy}
Let $f$ be invertible, $g=f^p$, and let $\tau$ be a $g$-invariant
probability.  Put
\begin{equation}\label{eq:cyclic-average-base}
 \overline\tau:=\frac1p\sum_{j=0}^{p-1}f_*^j\tau.
\end{equation}
Then
\begin{equation}\label{eq:cyclic-average-entropy}
 h_\tau(g)=p h_{\overline\tau}(f).
\end{equation}
\end{lemma}

\begin{proof}
Since $f^j$ commutes with $g=f^p$,
\[
 h_{f_*^j\tau}(g)=h_\tau(g)\quad(0\le j<p),
 \qquad
 f_*\overline\tau
 =\frac1p\left(\sum_{j=1}^{p-1}f_*^j\tau+f_*^p\tau\right)
 =\overline\tau.
\]
Affinity gives $h_{\overline\tau}(g)=h_\tau(g)$, and the power rule yields
\[
 h_\tau(g)=h_{\overline\tau}(f^p)=p h_{\overline\tau}(f).
\]
\end{proof}

\begin{lemma}[Ergodic decompositions of cyclic averages]
\label{lem:cyclic-average-decomposition}
Let $f$ be invertible, put $g=f^p$, and let $\zeta\in\mathcal M(g)$ have
$g$-ergodic decomposition probability $\tau_\zeta^g$.  For
$\rho\in\mathcal M_e(g)$, define
\begin{equation}\label{eq:cyclic-averaging-map}
 \mathscr A_p(\rho):=\frac1p\sum_{j=0}^{p-1}f_*^j\rho.
\end{equation}
Then $\mathscr A_p(\rho)\in\mathcal M_e(f)$,
\begin{equation}\label{eq:cyclic-component-entropy-scaling}
 h_\rho(g)=p h_{\mathscr A_p(\rho)}(f),
\end{equation}
and, for the cyclic average
$\overline\zeta=p^{-1}\sum_{j=0}^{p-1}f_*^j\zeta$,
\begin{equation}\label{eq:cyclic-decomposition-pushforward}
 (\mathscr A_p)_*\tau_\zeta^g=\tau_{\overline\zeta}^f.
\end{equation}
\end{lemma}

\begin{proof}
The measure $\mathscr A_p(\rho)$ is $f$-invariant.  If $B$ is
$f$-invariant, then it is $g$-invariant, so $\rho(B)\in\{0,1\}$ and
$\mathscr A_p(\rho)(B)=\rho(B)$.  Thus $\mathscr A_p(\rho)$ is
$f$-ergodic.  Equation \eqref{eq:cyclic-component-entropy-scaling} is
\cref{lem:cyclic-average-entropy} applied to $\rho$.  Finally,
\[
 \overline\zeta
 =\int\mathscr A_p(\rho)\,\dd\tau_\zeta^g(\rho).
\]
The measures in this integral are $f$-ergodic, so uniqueness of the ergodic
decomposition of $\overline\zeta$ gives
\eqref{eq:cyclic-decomposition-pushforward}.
\end{proof}

\begin{lemma}[Cyclic averaging of unstable lifts]
\label{lem:cyclic-average-lift}
Let $f\in\Diff^1(M)$, put $g=f^p$, and let $\tau$ be a
$g$-invariant probability whose ergodic components are hyperbolic with
exactly $k$ positive Lyapunov exponents.  Denote by $\wh\tau_g^u$ the
$G=F_f^p$-invariant measure obtained
by integrating their $g$-unstable lifts.  Define
\begin{equation}\label{eq:cyclic-average-lift-measures}
 \overline\tau:=\frac1p\sum_{j=0}^{p-1}f_*^j\tau,
 \qquad
 \wh R:=\frac1p\sum_{j=0}^{p-1}F_{f*}^j\wh\tau_g^u.
\end{equation}
Assume that almost every $f$-ergodic component of $\overline\tau$ is hyperbolic
and has exactly $k$ positive Lyapunov exponents.  Then
\begin{equation}\label{eq:cyclic-average-lift-identity}
 \wh R=\wh{\overline\tau}^{\,u}.
\end{equation}
\end{lemma}

\begin{proof}
Oseledets spaces for $f^p$ and $f$ coincide wherever both are defined;
finite shifts transport them by $df^j$.  Since
$f_*^j\tau\le p\overline\tau$, a full $\overline\tau$-measure regular
set has full $f_*^j\tau$-measure for every $0\le j<p$.  Consequently
\[
 F_{f*}^j\widehat\tau_g^u
   =\bigl(x\mapsto(x,E^u_{\overline\tau}(x))\bigr)_*(f_*^j\tau).
\]
Averaging and using \eqref{eq:cyclic-average-lift-measures} proves the identity.
\end{proof}

\begin{lemma}[Cyclic limits under varying maps]\label{lem:cyclic-limits}
Let $f_n\to f$ in $C^1$, fix $p$, and choose one
$f_n^p$-ergodic component $\nu_n^{(p)}$ of an ergodic
$f_n$-invariant probability $\nu_n$ which is hyperbolic of unstable
index $k$.  Denote its $f_n^p$-unstable lift by
$\wh\nu_n^{u,(p)}$.  If
\[
 \nu_n^{(p)}\to\mu^{(p)},
 \qquad
 \wh\nu_n^{u,(p)}\to\wh\eta^{(p)},
 \qquad
 \nu_n\to\mu,
 \qquad
 \wh\nu_n^u\to\wh\eta,
\]
then
\begin{equation}\label{eq:cyclic-limit-base}
 \mu=\frac1p\sum_{j=0}^{p-1}f_*^j\mu^{(p)}
\end{equation}
and
\begin{equation}\label{eq:cyclic-limit-lift}
 \wh\eta=\frac1p\sum_{j=0}^{p-1}F_{f*}^j\wh\eta^{(p)}.
\end{equation}
\end{lemma}

\begin{proof}
For each fixed $j<p$, $f_n^j\to f^j$ and $F_{f_n}^j\to F_f^j$ uniformly.
Pass to the limit in \eqref{eq:cyclic-base} and
\eqref{eq:cyclic-lift}.
\end{proof}

\subsection{Common submeasures and the finite-smoothness loss}

The next estimate turns the mass of a common submeasure into
total-variation control, and we also estimate the large-iterate loss by applying Yomdin theory.

\begin{lemma}[Common submeasure and total variation]
\label{lem:common-submeasure-TV}
If $P,Q$ are probabilities and a finite measure $R$ of mass $\beta$ satisfies
$R\le P,Q$, then,
\begin{equation}\label{eq:common-submeasure-TV}
 \|P-Q\|_{\TV}\le2(1-\beta).
\end{equation}
\end{lemma}

\begin{proof}
The positive measures $P-R$ and $Q-R$ both have mass $1-\beta$, so the
triangle inequality gives $\|P-Q\|_{\TV}\le2(1-\beta)$.
\end{proof}

For later use note that
\begin{equation}\label{eq:M-submultiplicative}
 \mathfrak M(f^p)\le\mathfrak M(f)^p.
\end{equation}
Therefore the loss in \eqref{eq:fixed-Delta} satisfies
\begin{equation}\label{eq:Yomdin-loss-divided}
 \frac{\Delta_r(f^p)}p
 \le\frac{\log C_{\rm Y}}p
 +\frac{3k^2}{r-1}\log\mathfrak M(f),
\end{equation}
and consequently
\begin{equation}\label{eq:Yomdin-loss-vanish}
 \lim_{r\to\infty}\limsup_{p\to\infty}
 \frac{\Delta_r(f^p)}p=0.
\end{equation}

\subsection{Proof of ergodic entropic continuity}

\begin{proof}[Proof of \cref{thm:core}]
Put
\[
 \lambda:=\lambda_{k-1}(f),
 \qquad
 h:=h_\mu(f).
\]
Invariance of $\mu$ follows from uniform convergence: for every
$\varphi\in C(M)$,
\[
 \int\varphi\circ f\,\dd\mu
 =\lim_n\int\varphi\circ f_n\,\dd\nu_n
 =\lim_n\int\varphi\,\dd\nu_n
 =\int\varphi\,\dd\mu.
\]
By the componentwise gap,
\begin{equation}\label{eq:positive-average-excess}
 h-\lambda
 =\int_{\mathcal M_e(f)}
   \bigl(h_\rho(f)-\lambda\bigr)\,\dd\tau_\mu(\rho)>0.
\end{equation}
The index $k$ is fixed by the hypotheses, both for the approximating
measures and almost every ergodic component of $\mu$.

\medskip\noindent\emph{Step 1: large iterates and cyclic components.}
Let $\wh\eta$ be an arbitrary subsequential limit of $\wh\nu_n^u$.  We
prove $\wh\eta=\wh\mu^u$.  Work along a subsequence converging to this
already fixed $\wh\eta$.  Fix numbers
\begin{equation}\label{eq:choose-a-vartheta}
 \lambda<a<\vartheta<h.
\end{equation}
Fix an integer $r\ge2$.  For all sufficiently large $p$, put
$g=f^p$, $g_n=f_n^p$, and $\Theta=p\vartheta>pa+1$.  By
\cref{lem:finite-iterate}, for all sufficiently large $n$,
\begin{equation}\label{eq:iterate-margin-main}
 A_{k-1}(g_n)\le pa.
\end{equation}
Define
\begin{equation}\label{eq:large-iterate-weight}
 w_{p,n}(x):=(\ell_k(f_n^p,x)-pa)^+,
 \qquad w_p(x):=(\ell_k(f^p,x)-pa)^+.
\end{equation}
For fixed $p$, $w_{p,n}\to w_p$ uniformly.  At $\mu$-almost every point,
Oseledets' theorem and the positivity of the first $k$ exponents give
\begin{equation}\label{eq:pointwise-lower-exterior-limit}
 \lim_{p\to\infty}\frac1p\ell_k(f^p,x)
 =\max_{0\le q<k}\sum_{i=1}^q\chi_i(x,f)
 =\sum_{i=1}^{k-1}\chi_i(x,f)
 \le\lambda<a.
\end{equation}
The empty sum is zero when $k=1$.  Choose $K\ge2$ bounding
$\|df\|_\infty$ and every $\|df_n\|_\infty$.  Since $a>\lambda\ge0$,
\[
 0\le w_{p,n}(x)\le p(k-1)\log K,
 \qquad 0\le w_p(x)\le p(k-1)\log K.
\]
Dominated convergence therefore yields
\begin{equation}\label{eq:localization-average-vanishes}
 \frac1p\int w_p\,\dd\mu\longrightarrow0.
\end{equation}
With
\begin{equation}\label{eq:large-iterate-error}
 R_p:=2+\lceil p(k-1)\log K\rceil,
 \qquad
 \varepsilon_p:=\frac1p\int w_p\,\dd\mu+\frac{1+\log R_p}{p},
\end{equation}
we thus have $\varepsilon_p\to0$.  The alphabet size $R_p$ depends only
linearly on $p$ before taking its logarithm.

For each $n$, choose among the finitely many ergodic $g_n$-components of
$\nu_n$ one, denoted by $\nu_n^{(p)}$, which minimizes the integral of
$w_{p,n}$.  Their equal-weight cyclic decomposition gives
\begin{equation}\label{eq:minimal-cycle-weight}
 \int w_{p,n}\,\dd\nu_n^{(p)}
 \le\int w_{p,n}\,\dd\nu_n.
\end{equation}
The Oseledets splitting of $f_n^p$ is the Oseledets splitting of $f_n$, with
all exponents multiplied by $p$; hence $\nu_n^{(p)}$ is hyperbolic with
exactly the same $k$ positive Lyapunov exponents.  By
\cref{lem:cyclic-components},
\begin{equation}\label{eq:iterate-component-entropy}
 h_{\nu_n^{(p)}}(g_n)=p h_{\nu_n}(f_n)\longrightarrow ph.
\end{equation}
Pass to a further subsequence so that
\begin{equation}\label{eq:iterate-component-limits}
 \nu_n^{(p)}\to\mu^{(p)},
 \qquad
 \wh\nu_n^{u,(p)}\to\wh\eta^{(p)}.
\end{equation}
The weight selection and uniform convergence of $w_{p,n}$ give
\begin{equation}\label{eq:limit-cycle-weight}
 \int w_p\,\dd\mu^{(p)}\le\int w_p\,\dd\mu.
\end{equation}
Uniform convergence of $g_n$ and of the induced Grassmann maps implies that
$\mu^{(p)}$ and $\wh\eta^{(p)}$ are invariant under $g$ and $F_f^p$,
respectively.  By \cref{lem:cyclic-limits},
\begin{equation}\label{eq:average-limits-main}
 \mu=\frac1p\sum_{j=0}^{p-1}f_*^j\mu^{(p)},
 \qquad
 \wh\eta=\frac1p\sum_{j=0}^{p-1}F_{f*}^j\wh\eta^{(p)}.
\end{equation}
The first identity and \cref{lem:cyclic-average-entropy} give
\begin{equation}\label{eq:limit-component-entropy}
 h_{\mu^{(p)}}(g)=p h_\mu(f)=ph.
\end{equation}
Moreover,
\[
 \frac1p\mu^{(p)}\le\mu
\]
as $g$-invariant measures.  Powers multiply the Lyapunov exponents by $p$,
and invariant domination preserves full-measure properties.  Thus almost
every $g$-ergodic component of $\mu^{(p)}$ is hyperbolic of index $k$.
By \cref{lem:componentwise-gap-inheritance}, first for the power $g=f^p$
and then for this domination, these components also satisfy
\begin{equation}\label{eq:cyclic-limit-componentwise-gap}
 h_\rho(g)>p\lambda=p\lambda_{k-1}(f).
\end{equation}
All hypotheses of the localized fixed-iterate decomposition are now
verified.  The integral bound \eqref{eq:limit-cycle-weight}, rather than
the coarser domination $\mu^{(p)}\le p\mu$, is used for its error term.

\medskip\noindent\emph{Step 2: the fixed-iterate decomposition and the
uncaptured mass.}
All limits above are taken along a further subsequence of the original
subsequence defining the already fixed measure $\wh\eta$.  Apply
\cref{thm:fixed-iterate-decomposition} to $(g_n,\nu_n^{(p)})$ with
$a$ there replaced by $pa$, threshold $\Theta=p\vartheta$, and alphabet
bound $R_p$.  We obtain
$\beta=\beta(r,p)\in[0,1]$ and, when $\beta>0$, a $g$-invariant probability
$\mu_1^{(p)}$ such that
\begin{equation}\label{eq:iterate-domination-main}
 \beta\mu_1^{(p)}\le\mu^{(p)},
 \qquad
 \beta\wh\mu_1^{u,(p)}\le\wh\eta^{(p)},
\end{equation}
and
\begin{equation}\label{eq:iterate-entropy-main}
 ph
 \le\beta h_{\mu_1^{(p)}}(g)
 +(1-\beta)p\vartheta+\Delta_r(f^p)+p\varepsilon_p.
\end{equation}
If $\beta<1$, define
\[
 \mu_0^{(p)}:=
 \begin{cases}
  \mu^{(p)},&\beta=0,\\[1mm]
  \displaystyle
  \frac{\mu^{(p)}-\beta\mu_1^{(p)}}{1-\beta},&0<\beta<1.
 \end{cases}
\]
Let $\tau^{(p)}$ be the $g$-ergodic decomposition probability of
$\mu^{(p)}$.  For $\beta<1$, let $\tau_0^{(p)}$ be that of $\mu_0^{(p)}$
and put $\sigma_0^{(p)}:=(1-\beta)\tau_0^{(p)}$.  For $\beta=1$, set
$\sigma_0^{(p)}:=0$.  This is the finite ergodic-decomposition measure of
$\mu^{(p)}-\beta\mu_1^{(p)}$, with the zero-term convention when $\beta=0$.
Since the complementary invariant finite measure is dominated by
$\mu^{(p)}$,
\cref{lem:dominated-lifts} gives
\begin{equation}\label{eq:neutral-decomposition-domination}
 \sigma_0^{(p)}\le\tau^{(p)}.
\end{equation}
Entropy affinity and \eqref{eq:limit-component-entropy} give
\begin{equation}\label{eq:iterate-entropy-affinity}
 ph
 =\beta h_{\mu_1^{(p)}}(g)
 +\int h_\rho(g)\,\dd\sigma_0^{(p)}(\rho),
\end{equation}
where the first term is zero when $\beta=0$.  Comparing this identity with
\eqref{eq:iterate-entropy-main} yields
\begin{equation}\label{eq:neutral-entropy-excess-estimate}
 \frac1p\int h_\rho(g)\,\dd\sigma_0^{(p)}(\rho)
 \le(1-\beta)\vartheta+\frac{\Delta_r(f^p)}p+\varepsilon_p.
\end{equation}

Let $\mathscr A_p$ be the cyclic-averaging map from
\eqref{eq:cyclic-averaging-map}, and define the positive entropy excess
\[
 e(\eta):=h_\eta(f)-\lambda
 \qquad(\eta\in\mathcal M_e(f)).
\]
By \cref{lem:cyclic-average-decomposition},
\begin{equation}\label{eq:cyclic-excess-distribution}
 (\mathscr A_p)_*\tau^{(p)}=\tau_\mu,
 \qquad
 \frac1p h_\rho(g)=h_{\mathscr A_p(\rho)}(f).
\end{equation}
Consequently \eqref{eq:neutral-entropy-excess-estimate} becomes
\begin{equation}\label{eq:integrated-excess-estimate}
 \int e\bigl(\mathscr A_p(\rho)\bigr)
       \,\dd\sigma_0^{(p)}(\rho)
 \le(1-\beta)(\vartheta-\lambda)
    +\frac{\Delta_r(f^p)}p+\varepsilon_p.
\end{equation}
For every $t>0$, Markov's inequality,
\eqref{eq:neutral-decomposition-domination}, and
\eqref{eq:cyclic-excess-distribution} now give
\begin{align}
 1-\beta
 &\le
 \tau_\mu\{\eta:e(\eta)\le t\}
 +\frac{\vartheta-\lambda}{t}
 +\frac{\Delta_r(f^p)}{pt}+\frac{\varepsilon_p}{t}.
 \label{eq:componentwise-beta-estimate}
\end{align}
Indeed, split $\sigma_0^{(p)}$ over the sets where
$e\circ\mathscr A_p\le t$ and $e\circ\mathscr A_p>t$; the first part is
bounded by $\tau^{(p)}$ and the second by
$t^{-1}\int e\circ\mathscr A_p\,\dd\sigma_0^{(p)}$.

\medskip\noindent\emph{Step 3: averaging the hyperbolic component back to
$f$.}
For $\beta>0$, define
\begin{equation}\label{eq:bar-mu1}
 \overline\mu_1:=\frac1p\sum_{j=0}^{p-1}f_*^j\mu_1^{(p)},
 \qquad
 \wh R_1:=\frac1p\sum_{j=0}^{p-1}F_{f*}^j\wh\mu_1^{u,(p)}.
\end{equation}
Averaging the first inequality in
\eqref{eq:iterate-domination-main} and using
\eqref{eq:average-limits-main} gives
\begin{equation}\label{eq:averaged-base-domination}
 \beta\overline\mu_1\le\mu.
\end{equation}
Before averaging, apply \cref{lem:dominated-lifts} to the
$g$-invariant finite measures
$\beta\mu_1^{(p)}\le\mu^{(p)}$.  Almost every $g$-ergodic component of
$\mu_1^{(p)}$ therefore occurs among those of $\mu^{(p)}$, and is
hyperbolic with exactly $k$ positive Lyapunov exponents.  Apply the same
lemma to \eqref{eq:averaged-base-domination} for the map $f$.  Almost every
$f$-ergodic component of $\overline\mu_1$ occurs among those of $\mu$, and
is again hyperbolic with exactly $k$ positive exponents.  We may thus apply
\cref{lem:cyclic-average-lift} to $\tau=\mu_1^{(p)}$, obtaining
\begin{equation}\label{eq:averaged-lift-identification}
 \wh R_1=\wh{\overline\mu_1}^{\,u}.
\end{equation}
Averaging the second inequality in
\eqref{eq:iterate-domination-main} and using
\eqref{eq:average-limits-main} yields
\begin{equation}\label{eq:averaged-domination}
 \beta\wh R_1\le\wh\eta.
\end{equation}
Finally, \cref{lem:dominated-lifts}, applied to
\eqref{eq:averaged-base-domination}, gives
\begin{equation}\label{eq:averaged-domination-muu}
 \beta\wh R_1\le\wh\mu^u.
\end{equation}
When $\beta=0$, take the common submeasure to be zero.  Hence
\cref{lem:common-submeasure-TV} and
\eqref{eq:componentwise-beta-estimate} imply, for every $t>0$,
\begin{equation}\label{eq:TV-bound-main}
 \|\wh\eta-\wh\mu^u\|_{\TV}
 \le2\left(
 \tau_\mu\{\rho:0<h_\rho(f)-\lambda\le t\}
 +\frac{\vartheta-\lambda}{t}
 +\frac{\Delta_r(f^p)}{pt}+\frac{\varepsilon_p}{t}
 \right).
\end{equation}
The set in the first term may equivalently be written
$\{\rho:e(\rho)\le t\}$ modulo a $\tau_\mu$-null set.

\medskip\noindent\emph{Step 4: vanishing loss and convergence of lifts.}
The estimates above are valid for every choice
$\lambda<a<\vartheta<h$, every $r\ge2$, and all sufficiently large $p$.
By hypothesis, $e(\rho)>0$ for $\tau_\mu$-almost every $\rho$.  Therefore
\[
 \tau_\mu\{\rho:0<e(\rho)\le t\}\longrightarrow0
 \qquad(t\downarrow0).
\]
Given $\delta>0$, first choose $t>0$ so that the first term in
\eqref{eq:TV-bound-main} is sufficiently small.  Next choose
$\vartheta\in(\lambda,h)$ sufficiently close to $\lambda$, and then choose
$a\in(\lambda,\vartheta)$.  Next choose $r$ and then $p$ so that the
Yomdin term is sufficiently small, using \eqref{eq:Yomdin-loss-vanish}.
For these fixed choices of $t,a,\vartheta,r$, enlarge $p$ further so that
$\varepsilon_p/t$ is small, by \eqref{eq:large-iterate-error}, and
\eqref{eq:iterate-margin-main} holds.  The geometric scale may depend on
$p$, since the limits in \cref{thm:fixed-iterate-decomposition} are taken
with $p$ fixed.

Starting from the fixed subsequence $\wh\nu_n^u\to\wh\eta$, extract the
cyclic $f_n^p$-components and all subsequential limits used above.  This is
only a further subsequence, so its averaged lift still converges to the same
$\wh\eta$.  Thus \eqref{eq:TV-bound-main} holds for the originally fixed
limit $\wh\eta$.  Since $\delta$ is arbitrary,
$\wh\eta=\wh\mu^u$.  Every subsequential limit is the same, proving
\eqref{eq:core-lift-convergence}.

Finally, $f_n\to f$ in $C^1$ implies
\[
 \Phi_{f_n}\longrightarrow\Phi_f
 \quad\text{uniformly on }\wh M_k.
\]
Together with \eqref{eq:positive-sum} and
\eqref{eq:core-lift-convergence},
\[
 S^+(\nu_n,f_n)
 =\int\Phi_{f_n}\,\dd\wh\nu_n^u
 \longrightarrow
 \int\Phi_f\,\dd\wh\mu^u
 =S^+(\mu,f).
\]
This proves both assertions.
\end{proof}

\subsection{General invariant measures}\label{sec:general-measures}

We now pass from ergodic approximating measures to arbitrary invariant
probabilities.  We first record the compactness argument for their ergodic
decompositions; this argument uses entropy upper semicontinuity, but no
Lyapunov continuity.

For an invariant probability $\nu$, let $\tau_\nu$  be its ergodic decomposition probability  on $\mathcal M_e(f)$. The map $\rho\mapsto\tau_\rho$ is a Borel
probability kernel from $\mathcal M(f)$ to
$\mathcal P(\mathcal M_e(f))$.
 
 The bounded Borel function
\begin{equation}\label{eq:universal-positive-sum-function}
 \Sigma_g(x):=
 \max_{0\le q\le d}\limsup_{m\to\infty}\frac1m
       \log\|\wedge^q dg_x^m\|,
 \qquad \wedge^0dg_x^m=1,
\end{equation}
equals the sum of the positive Lyapunov exponents at every Oseledets-regular
point.  Hence
\begin{equation}\label{eq:Splus-universal-integral}
 S^+(\nu,g)=\int\Sigma_g\,\dd\nu.
\end{equation}
In particular, $S^+$ is Borel and harmonic on $\mathcal M(g)$: if a
probability $\tau$ on $\mathcal M(g)$ has barycenter $\nu$, then Fubini's
theorem gives
\begin{equation}\label{eq:Splus-ergodic-decomposition-general}
 S^+(\nu,g)=\int S^+(\rho,g)\,\dd\tau(\rho).
\end{equation}
We use the analogous harmonicity of metric entropy, which follows by
integrating the ergodic-decomposition formula for entropy.

\begin{lemma}[Entropy coordinates of ergodic decompositions]
\label{lem:entropy-graph-compactness}
Let $f_n\to f$ in $C^\infty$, let $\nu_n\in\mathcal M(f_n)$, and assume
\[
 \nu_n\to\mu,\qquad h_{\nu_n}(f_n)\to h_\mu(f).
\]
Write $\tau_n:=\tau_{\nu_n}$.  There are constants $C_h,C_S<\infty$ such
that the probabilities
\begin{equation}\label{eq:triple-graph-measure}
 \Gamma_n:=
 \bigl(\rho\mapsto(\rho,h_\rho(f_n),S^+(\rho,f_n))\bigr)_*\tau_n
\end{equation}
are supported on the compact space
$\mathcal P(M)\times[0,C_h]\times[0,C_S]$.  If $\Gamma_n\to\Gamma$ along
a subsequence and $\tau$ is the first marginal of $\Gamma$, then
\begin{equation}\label{eq:barycenter-general}
 \tau(\mathcal M(f))=1,\qquad
 \mu=\int\rho\,\dd\tau(\rho),
\end{equation}
\begin{equation}\label{eq:graph-entropy-equality}
 \mathfrak h=h_\rho(f)
 \quad\text{for }\Gamma\text{-almost every }(\rho,\mathfrak h,s),
\end{equation}
and 

 \begin{equation}\label{eq:integrated-ergodic-decomposition}
 \tau_\mu=\int\tau_\rho\,\dd\tau(\rho).
\end{equation}

More generally, suppose $E_n\subset\mathcal M_e(f_n)$ are Borel and the
restrictions of $\Gamma_n$ to graph points with base in $E_n$ converge to
a finite measure $\Gamma'$.  Then every point of
$\operatorname{supp}\Gamma'$ is a limit of actual graph points with bases
$\rho_j\in E_{n_j}$ and $n_j\to\infty$.
\end{lemma}

\begin{proof}
The entropy and Lyapunov-sum coordinates are Borel.  A common $C^1$ bound
and Ruelle's inequality give
\begin{equation}\label{eq:uniform-entropy-Splus-bounds}
 0\le h_\rho(f_n)\le C_h,
 \qquad 0\le S^+(\rho,f_n)\le C_S
 \quad(\rho\in\mathcal M(f_n)).
\end{equation}
The first marginals therefore satisfy $\tau_n\to\tau$.  Continuity of the
barycenter map gives the second identity in \eqref{eq:barycenter-general}.

Every open ball centered at a point of $\operatorname{supp}\Gamma$ has
positive $\Gamma$-mass.  Portmanteau's theorem and a decreasing sequence of
such balls give $n_j\to\infty$ and $\rho_j\in\mathcal M_e(f_{n_j})$ with
\begin{equation}\label{eq:triple-graph-approximating-sequence}
 (\rho_j,h_{\rho_j}(f_{n_j}),S^+(\rho_j,f_{n_j}))
 \longrightarrow(\rho,\mathfrak h,s).
\end{equation}
Uniform convergence implies $\rho\in\mathcal M(f)$, and joint upper
semicontinuity (Theorem~\ref{thm:joint-usc}) gives
$\mathfrak h\le h_\rho(f)$.  On the other hand,
weak convergence of $\Gamma_n$, entropy convergence, and harmonicity give
\[
 \int\mathfrak h\,\dd\Gamma
 =\lim_n h_{\nu_n}(f_n)
 =h_\mu(f)
 =\int h_\rho(f)\,\dd\tau(\rho).
\]
Thus the nonnegative difference $h_\rho(f)-\mathfrak h$ has zero integral,
proving \eqref{eq:graph-entropy-equality}.

Decompose each $\rho$ in \eqref{eq:barycenter-general} into its ergodic
components.  Uniqueness of the ergodic decomposition of $\mu$ gives
\eqref{eq:integrated-ergodic-decomposition}.  Finally, the same positive-ball
argument applied to the restricted graph measures proves the last assertion, and 
the selected points can be taken in $E_{n_j}$ because the restrictions are
carried by those graph points.
\end{proof}

\begin{theorem}[Fixed-index entropic continuity]
\label{thm:fixed-index-continuity}
Let $f_n,f\in\Diff^\infty(M)$ with $f_n\to f$ in $C^\infty$, and let
$\nu_n\in\mathcal M(f_n)$ and $\mu\in\mathcal M(f)$.  Suppose that
\begin{equation}\label{eq:fixed-index-continuity-assumptions}
 \nu_n\to\mu,\qquad h_{\nu_n}(f_n)\to h_\mu(f).
\end{equation}
Assume that, for some fixed $1\le k\le d-1$,
\begin{equation}\label{eq:asymptotic-index-concentration}
 \tau_{\nu_n}\bigl(\mathscr H_k(f_n)\bigr)\longrightarrow1,
\end{equation}
and that, for $\tau_\mu$-almost every $\rho$,
\begin{equation}\label{eq:fixed-index-componentwise-gap}
 \rho\in\mathscr H_k(f),\qquad h_\rho(f)>\lambda_{k-1}(f).
\end{equation}
Then
\begin{equation}\label{eq:fixed-index-Splus-continuity}
 S^+(\nu_n,f_n)\longrightarrow S^+(\mu,f).
\end{equation}
If the $\nu_n$ are ergodic, they eventually belong to $\mathscr H_k(f_n)$,
and their unstable lifts converge to $\widehat\mu^u$.
\end{theorem}

\begin{proof}
Fix an arbitrary subsequence and extract the graph limit $\Gamma$ from
\cref{lem:entropy-graph-compactness}.  Set $\lambda:=\lambda_{k-1}(f)$ and
\[
 \mathscr E_>:=
 \{\eta\in\mathscr H_k(f):h_\eta(f)>\lambda\}.
\]
It has full $\tau_\mu$-measure by
\eqref{eq:fixed-index-componentwise-gap}.  The identity
\eqref{eq:integrated-ergodic-decomposition} implies
$\tau_\rho(\mathscr E_>)=1$ for $\tau$-almost every $\rho$, where $\tau$ is the limit of $(\tau_{\nu_n})_n$.

Restrict $\Gamma_n$ to graph points with base in $\mathscr H_k(f_n)$.
The omitted mass tends to zero by
\eqref{eq:asymptotic-index-concentration}, so these restrictions have the
same limit $\Gamma$.  For $\Gamma$-almost every $(\rho,\mathfrak h,s)$,
\cref{lem:entropy-graph-compactness} therefore supplies a sequence
\[
 \rho_j\in\mathscr H_k(f_{n_j}),\qquad
 \rho_j\to\rho,\qquad h_{\rho_j}(f_{n_j})\to h_\rho(f),\qquad
 S^+(\rho_j,f_{n_j})\to s,
\]
and almost every ergodic component of $\rho$ lies in $\mathscr E_>$.
Apply \cref{thm:core} to obtain
\begin{equation}\label{eq:graph-support-identification}
 s=S^+(\rho,f)\quad\text{for }\Gamma\text{-almost every }
 (\rho,\mathfrak h,s).
\end{equation}
The third coordinate is bounded and continuous on the ambient compact
product.  Equations \eqref{eq:Splus-ergodic-decomposition-general} and
\eqref{eq:barycenter-general} now give
\[
 \lim_n S^+(\nu_n,f_n)
 =\int s\,\dd\Gamma
 =\int S^+(\rho,f)\,\dd\tau(\rho)
 =S^+(\mu,f).
\]
Since the initial subsequence was arbitrary, the full sequence converges.
If $\nu_n$ is ergodic, the left side of
\eqref{eq:asymptotic-index-concentration} is either zero or one, and hence
is eventually one.  The lift assertion then follows from \cref{thm:core}.
\end{proof}

\subsection{Proof of Theorem A}

\begin{proof}[Proof of \cref{thm:general-entropic-continuity}]
Since $f$ is entropy-dominant, \cref{lem:index}
selects a unique degree $k$ with $\lambda_k(f)>\lambda_{\mathrm{sec}}(f)$, while
\[
 \lambda_q(f)\le\lambda_{\mathrm{sec}}(f)\qquad(0\le q\le d-1,\ q\ne k).
\]
Moreover $\tau_\mu$-almost every $\rho$ has entropy larger than $\lambda_{\mathrm{sec}}(f)\ge \lambda_{k-1}(f) $, so that  such $\rho$ is hyperbolic of
index $k$.  
This verifies \eqref{eq:fixed-index-componentwise-gap}.

It remains to establish \eqref{eq:asymptotic-index-concentration}.  Suppose that it fails.  Pass to a
subsequence on which the mass of
$\mathcal M_e(f_n)\setminus\mathscr H_k(f_n)$ is at least $\delta>0$.
Apply \cref{lem:entropy-graph-compactness}, and pass further so that the
restrictions of $\Gamma_n$ to these complementary graph points converge
to a finite measure $\Gamma'$.  Then
\[
 0\le\Gamma'\le\Gamma,\qquad \Gamma'(\mathcal P(M)\times[0,C_h]\times[0,C_S])
 \ge\delta.
\]
For $\Gamma$-almost every $(\rho,\mathfrak h,s)$, the entropy coordinate
satisfies $\mathfrak h=h_\rho(f)$, and
\eqref{eq:integrated-ergodic-decomposition} transfers the original
componentwise gap to $\rho$.  Consequently $h_\rho(f)>\lambda_{\mathrm{sec}}(f)$ there.
Choose such a point in $\operatorname{supp}\Gamma'$.  The last assertion
of \cref{lem:entropy-graph-compactness} gives
\[
 \begin{gathered}
 \rho_j\in\mathcal M_e(f_{n_j})\setminus\mathscr H_k(f_{n_j}),
 \qquad \rho_j\to\rho,\\
 h_{\rho_j}(f_{n_j})\to h_\rho(f)>\lambda_{\mathrm{sec}}(f).
 \end{gathered}
\]

But by \cref{lem:index} such $\rho_j$ should belong to $\mathscr H_k(f_{n_j})$,
a contradiction.  Thus \eqref{eq:asymptotic-index-concentration} holds.
\Cref{thm:fixed-index-continuity} now proves
\eqref{eq:general-Splus}.
\end{proof}

\section{Stable SPR and the fixed-index criterion}
\label{sec:main-theorem-proof}

We prove a criterion in which the hyperbolicity of high-entropy
measures is assumed.  The exterior-growth gap in \cref{thm:main-intro}
naturally  implies this assumption.  

\begin{theorem}[A fixed-index criterion for stable SPR]
\label{thm:fixed-dimension-variant}
Let $U\subset\Diff^\infty(M)$ be nonempty and $C^\infty$-open.  Suppose
that there exist
\[
 k\in\{1,\ldots,d-1\},\qquad \chi>0,\qquad
 h_*<\inf_{g\in U}h_{\mathrm{top}}(g)
\]
such that, for every $g\in U$ and $\rho\in\mathcal M_e(g)$,
\begin{equation}\label{eq:criterion-high-entropy-gap}
 h_\rho(g)>h_*
 \quad\Longrightarrow\quad
 \chi_k(\rho,g)>\chi,\qquad \chi_{k+1}(\rho,g)<-\chi.
\end{equation}
Assume also that
\begin{equation}\label{eq:criterion-exterior-gap}
 h_{\mathrm{top}}(g)>\lambda_{k-1}(g)
 \qquad(g\in U).
\end{equation}
Then conclusions \textnormal{(1)}--\textnormal{(4)} of
\cref{thm:main-intro} hold with $U$ in place of $\mathcal S^\infty(M)$.
If $f\in U$ has a unique measure of maximal entropy, conclusion
\textnormal{(5)} also holds on a $C^\infty$ neighborhood of $f$.
\end{theorem}

Until the final subsection, fix $U,k,\chi,h_*$ as in this theorem.
In particular, every ergodic maximal measure belongs to $\mathscr H_k(g)$.

\subsection{Entropy continuity and strong positive recurrence}

\begin{lemma}[Entropy lower semicontinuity from a hyperbolic maximal measure]
\label{lem:entropy-lsc-gap}
If $g\in\Diff^\infty(M)$ has a hyperbolic ergodic measure of maximal
entropy and $h_{\mathrm{top}}(g)>0$, then $h_{\mathrm{top}}$ is lower
semicontinuous at $g$ in the $C^1$ topology relative to
$\Diff^\infty(M)$.
\end{lemma}

\begin{proof}
Let $\rho$ be such a measure.  Katok's horseshoe theorem \cite{Katok}
gives, for every $\delta>0$, a locally maximal hyperbolic set $K$ with
\[
 h_{\mathrm{top}}(g|_K)>h_\rho(g)-\delta
 =h_{\mathrm{top}}(g)-\delta.
\]
Its hyperbolic continuation $K_{g'}$ is conjugate to $K$ for every $g'$
sufficiently $C^1$-close to $g$.  Thus
\[
 h_{\mathrm{top}}(g')\ge h_{\mathrm{top}}(g'|_{K_{g'}})
 =h_{\mathrm{top}}(g|_K)>h_{\mathrm{top}}(g)-\delta.
\]
Let $\delta\downarrow0$.
\end{proof}

\begin{proposition}[Continuity of topological entropy]
\label{prop:entropy-continuity-neighborhood}
The map $g\mapsto h_{\mathrm{top}}(g)$ is continuous on $U$ in the
$C^\infty$ topology.
\end{proposition}

\begin{proof}
Newhouse upper semicontinuity implies existence of a maximal measure for
every $C^\infty$ diffeomorphism: take a maximizing sequence and a weak
limit.  Almost every ergodic component of a maximal measure is maximal.
For $g\in U$, these components are hyperbolic by
\eqref{eq:criterion-high-entropy-gap}, and
\eqref{eq:criterion-exterior-gap} gives $h_{\mathrm{top}}(g)>0$.
Hence \cref{lem:entropy-lsc-gap} gives lower semicontinuity.

For upper semicontinuity, let $g_n\to g$ in $C^\infty$ and pass to a
subsequence realizing the limsup of $h_{\mathrm{top}}(g_n)$.  Choose an
ergodic maximal measure $\rho_n$ of $g_n$ and extract $\rho_n\to\rho$.
By the joint upper semicontinuity theorem
(Theorem~\ref{thm:joint-usc}),
\[
 \limsup_n h_{\mathrm{top}}(g_n)
 =\limsup_n h_{\rho_n}(g_n)
 \le h_\rho(g)\le h_{\mathrm{top}}(g).
\]
The two semicontinuity statements prove continuity.
\end{proof}

\begin{proposition}[Near-maximal entropic continuity]
\label{prop:fixed-index-maximal-continuity}
Let $g_n,g\in U$ with $g_n\to g$ in $C^\infty$, and suppose that
$\nu_n\in\mathcal M(g_n)$ satisfy
\[
 \nu_n\to\mu,\qquad h_{\nu_n}(g_n)\to h_{\mathrm{top}}(g).
\]
Then $\mu$ is a  measure of maximal entropy for $g$ and
\[
 S^+(\nu_n,g_n)\to S^+(\mu,g).
\]
When the $\nu_n$ are ergodic, their unstable lifts also converge to
$\widehat\mu^u$.
\end{proposition}

\begin{proof}
Write $H_n:=h_{\mathrm{top}}(g_n)$ and $H:=h_{\mathrm{top}}(g)$.
The upper semicontinuity in Theorem \ref{thm:joint-usc} gives $h_\mu(g)=H$.
By \cref{prop:entropy-continuity-neighborhood}, $H_n\to H>h_*$.
If $\tau_n$ is the ergodic-decomposition probability of $\nu_n$, set
\[
 b_n:=\tau_n\{\rho:h_\rho(g_n)\le h_*\}.
\]
Entropy affinity yields
\[
 h_{\nu_n}(g_n)\le b_n h_*+(1-b_n)H_n,
 \qquad
 0\le b_n\le\frac{H_n-h_{\nu_n}(g_n)}{H_n-h_*}\longrightarrow0.
\]
All remaining components lie in $\mathscr H_k(g_n)$ by
\eqref{eq:criterion-high-entropy-gap}.  Almost every ergodic component of
$\mu$ is maximal, belongs to $\mathscr H_k(g)$, and has entropy
$H>\lambda_{k-1}(g)$.  Thus \cref{thm:fixed-index-continuity} applies.
Its last assertion gives the lift convergence in the ergodic case.
\end{proof}

\begin{proposition}[Stable strong positive recurrence]
\label{prop:stable-SPR}
Every $g\in U$ is strongly positively recurrent and has finitely many
ergodic measures of maximal entropy.
\end{proposition}

\begin{proof}
Fix $g\in U$.  We verify the hypotheses (EH) and (EC) of
\cite[Section~3.2, Theorem~3.1]{SPR2025} on $X=M$.
If $\rho_j$ are ergodic $g$-invariant probabilities with
$h_{\rho_j}(g)\to h_{\mathrm{top}}(g)$ and $\rho_j\to\rho$, then
$\rho$ is maximal.  By \eqref{eq:criterion-high-entropy-gap},
\[
 \chi_k(x,g)>\chi>-\chi>\chi_{k+1}(x,g)
 \quad\text{for }\rho\text{-almost every }x.
\]
This is (EH).  Applied to the constant sequence of maps,
\cref{prop:fixed-index-maximal-continuity} gives
$S^+(\rho_j,g)\to S^+(\rho,g)$, which is (EC).
\cite[Section~3.2, Theorem~3.1]{SPR2025}  gives SPR, and finiteness of the ergodic maximal
measures follows from \cite[Theorem~11.6(3)]{SPR2025}.
\end{proof}

\subsection{Uniform weighted Pesin blocks and homoclinic capture}

We recall the terminology needed in this subsection.  If an invariant
probability $\eta$ gives full mass to a Borel set $Y$, we say that $Y$
\emph{carries} $\eta$.  The Borel homoclinic relation of
\cite[Definitions~2.5 and~2.7]{SPR2025} is defined on the recurrent
nonuniformly hyperbolic set $\operatorname{NUH}'(g)$, where the canonical
global Pesin stable and unstable manifolds are available.  For
$x,y\in\operatorname{NUH}'(g)$, one writes $x\sim_g y$ when
\[
 W^s_g(x)\pitchfork g^mW^u_g(y)\ne\varnothing,
 \qquad
 W^u_g(x)\pitchfork g^nW^s_g(y)\ne\varnothing
\]
for some $m,n\in\mathbb Z$.  Its equivalence classes are the Borel
homoclinic classes.  If $Y$ is an invariant Borel set, we use the variational
notation
\[
 h_{\mathrm{top}}(g|Y)
 :=\sup\{h_\nu(g):\nu\in\mathcal M(g),\ \nu(Y)=1\}
\]
for its entropy.  Every hyperbolic invariant measure gives full mass to
$\operatorname{NUH}'(g)$, and every ergodic hyperbolic measure is carried by
a unique Borel homoclinic class; each class contains a dense family of
hyperbolic periodic orbits
\cite[Theorem~2.4 and Proposition~2.8]{SPR2025}.

We also briefly recall the weighted Pesin blocks used below.  For a map $T$,
a weighted $\ell$-Pesin block is a compact set $A$ whose full orbit
saturation $\bigcup_{j\in\mathbb Z}T^jA$ carries a $dT$-invariant splitting
\[
 TM=E^{cs}\oplus E^u,
 \qquad \dim E^u=\ell,
\]
with exponential estimates modified by the Birkhoff sums of a nonnegative
continuous weight $\varphi$.  More precisely, if
\[
 \Phi_m^\varphi(z):=\sum_{i=0}^{m-1}\varphi(T^iz),
 \qquad
 \Phi_{-m}^\varphi(z):=-\sum_{i=1}^{m}\varphi(T^{-i}z),
\]
then the defining bounds have the form
\begin{align*}
 \|dT^m|_{E^{cs}(T^jx)}\|
 &\le K e^{-m\chi+|j|\eps+\Phi_m^\varphi(T^jx)},\\
 \|dT^{-m}|_{E^u(T^jx)}\|
 &\le K e^{-m\chi+|j|\eps+\Phi_{-m}^\varphi(T^jx)}.
\end{align*}
The splitting is continuous on every finite saturation
$A_L=\bigcup_{|j|\le L}T^jA$: shifting by at most $L$ replaces $K$ by
$Ke^{L\eps}$.  See \cite[Definition~8 and the properties following it]{LTTTY}.
What matters here is that the local
weighted unstable manifolds are defined at every point of the orbit
saturation, have a uniform size on each finite orbit saturation, and vary in
the $C^1$ topology.  At regular points satisfying the relevant spectral-gap
and Birkhoff-regularity conditions, they coincide with the classical Pesin
unstable manifolds
\cite[Lemma~2.10 and Proposition~2.15]{LTTTY}.

For later reference, if the Lyapunov exponents of a map $T$ are listed in
nonincreasing order, put
\begin{equation}\label{eq:Gamma-ell-definition}
 \Gamma_\ell(T):=
 \left\{x:\begin{array}{l}
 x\text{ is Oseledets regular for }T,\\[-1mm]
 \chi_\ell(x,T)>\max\{\chi_{\ell+1}(x,T),0\}
 \end{array}\right\}.
\end{equation}
For an invariant probability $\eta$, the quantitative gap used in
\cite{LTTTY} is
\begin{equation}\label{eq:Hyp-ell-definition}
 \operatorname{Hyp}^{\ell}(\eta,T)
 :=\operatorname*{ess\,inf}_{\eta}
 \bigl(\chi_\ell(x,T)-\max\{\chi_{\ell+1}(x,T),0\}\bigr).
\end{equation}
Finally, for compact sets $K_n\subset M$, we use the sequential upper limit
\begin{equation}\label{eq:set-upper-limit-definition}
 \limsup_{n\to\infty}K_n
 :=\bigl\{x:\exists n_j\to\infty,\ x_{n_j}\in K_{n_j},\
                  x_{n_j}\to x\bigr\}.
\end{equation}
It is compact because $M$ is compact.

When the forward and inverse constructions below use
$T^u=g^{N_u}$ and $T^s=g^{-N_s}$, respectively, we write
$W^{u,\mathrm w}_{\mathrm{loc}}(x,g)$ for the local weighted unstable disk
constructed for $T^u$, and $W^{s,\mathrm w}_{\mathrm{loc}}(x,g)$ for the
local weighted unstable disk constructed for $T^s$.  The latter is a stable
disk for the original map $g$.  We use the analogous notation for the
approximating maps.

\begin{lemma}[Uniform two-sided weighted Pesin blocks]
\label{lem:uniform-two-sided-blocks}
Let $g_n,g\in U$, let $g_n\to g$ in $C^\infty$, and let $\rho_n$ be an
ergodic measure of maximal entropy for $g_n$.  Suppose that
$\rho_n\to\rho$.  For every $\delta>0$, after discarding finitely many
indices, there exist integers $N_u,N_s\ge1$, compact forward weighted blocks
$A_n^u,A^u$, compact inverse weighted blocks $A_n^s,A^s$, compact sets
$\Lambda_n,\Lambda\subset M$, and Borel sets
$\mathcal G_n,\mathcal G\subset M$ with the following properties.

The forward blocks are weighted $k$-Pesin blocks for
\[
 T_n^u:=g_n^{N_u},\qquad T^u:=g^{N_u},
\]
and the inverse blocks are weighted $(d-k)$-Pesin blocks for
\[
 T_n^s:=g_n^{-N_s},\qquad T^s:=g^{-N_s}.
\]
Moreover,
\begin{equation}\label{eq:uniform-block-masses}
 \rho_n(\Lambda_n)>1-\delta,
 \qquad
 \rho(\Lambda)>1-\delta,
 \qquad
 \limsup_{n\to\infty}\Lambda_n\subset\Lambda,
\end{equation}
and
\begin{equation}\label{eq:full-regular-homoclinic-sets}
 \rho_n(\mathcal G_n)=1,
 \qquad
 \rho(\mathcal G)=1,
 \qquad
 \mathcal G_n\subset\operatorname{NUH}'(g_n),
 \qquad
 \mathcal G\subset\operatorname{NUH}'(g).
\end{equation}
On $\Lambda_n$ and $\Lambda$ there are local weighted stable and unstable
manifolds of a common positive size.  At every point of
$\Lambda_n\cap\mathcal G_n$, respectively $\Lambda\cap\mathcal G$, these
weighted local disks lie in the classical Pesin stable and unstable
manifolds and give their usual local pieces after restriction.
If $x_n\in\Lambda_n$ and $x_n\to x\in\Lambda$, then, after restriction to
the common size,
\[
 W^{s,\mathrm w}_{\rm loc}(x_n,g_n)
 \longrightarrow W^{s,\mathrm w}_{\rm loc}(x,g),
 \qquad
 W^{u,\mathrm w}_{\rm loc}(x_n,g_n)
 \longrightarrow W^{u,\mathrm w}_{\rm loc}(x,g)
\]
in the $C^1$ topology.
\end{lemma}

\begin{proof}
By \cref{prop:entropy-continuity-neighborhood},
\[
 h_{\rho_n}(g_n)=h_{\rm top}(g_n)
 \longrightarrow h_{\rm top}(g).
\]
The upper semicontinuity in Theorem \ref{thm:joint-usc} implies that $\rho$ is a measure of maximal entropy for
$g$.  The hypotheses \eqref{eq:criterion-high-entropy-gap}--
\eqref{eq:criterion-exterior-gap} apply to all the $\rho_n$ and to almost
every ergodic component of $\rho$.  Thus \cref{thm:core} gives
\begin{equation}\label{eq:uniform-block-Splus}
 S^+(\rho_n,g_n)\longrightarrow S^+(\rho,g).
\end{equation}
Almost every ergodic component of $\rho$ satisfies
\eqref{eq:criterion-high-entropy-gap}, hence
$\operatorname{Hyp}^k(\rho,g)>0$ in the sense of
\eqref{eq:Hyp-ell-definition}.  Because the first $k$ exponents are exactly
the positive exponents, \eqref{eq:uniform-block-Splus} is the convergence of
the sums of the first $k$ Lyapunov exponents.

Fix $0<\alpha<1$.  The $C^\infty$ convergence gives a common
$C^{1+\alpha}$ bound.  In both applications below, choose the weighted-block
parameters with the smallness condition of \cite[Equation~(4)]{LTTTY},
using common norm bounds for the sequence.  This choice also gives the
controlled finite-saturation radii used below.  Applying
\cite[Propositions~2.11 and~2.15]{LTTTY} to the forward maps gives a common
iterate $N_u$ and compact weighted $k$-Pesin blocks $A_n^u,A^u$ of
arbitrarily large measure such that
\begin{equation}\label{eq:forward-weighted-upper-limit}
 \limsup_{n\to\infty}A_n^u\subset A^u.
\end{equation}
Their local weighted unstable manifolds have a common positive size and
vary in the $C^1$ topology.

For the stable aspect, put
\[
 S^-(\eta,h):=\int\sum_{\chi_i(x,h)<0}\chi_i(x,h)\,\dd\eta(x).
\]
Since all measures under consideration have index $k$,
\begin{equation}\label{eq:negative-sum-determinant}
 S^-(\rho_n,g_n)
 =\int\log|\det dg_n|\,\dd\rho_n-S^+(\rho_n,g_n).
\end{equation}
The determinant integrals converge by $C^1$ convergence and weak
convergence.  Thus \eqref{eq:uniform-block-Splus} implies
$S^-(\rho_n,g_n)\to S^-(\rho,g)$.  Equivalently, the sums of the first
$d-k$ exponents of $(g_n^{-1},\rho_n)$ converge to the corresponding sum
for $(g^{-1},\rho)$.  The uniform gap
\eqref{eq:criterion-high-entropy-gap} gives
$\operatorname{Hyp}^{d-k}(\rho,g^{-1})>0$.  Applying the same two
propositions to the inverse maps gives a common iterate $N_s$ and compact
weighted $(d-k)$-Pesin blocks $A_n^s,A^s$ such that
\begin{equation}\label{eq:inverse-weighted-upper-limit}
 \limsup_{n\to\infty}A_n^s\subset A^s.
\end{equation}
The weighted unstable manifolds for $T_n^s=g_n^{-N_s}$ are the weighted
stable manifolds used for $g_n$.

Let $\varphi^u$ and $\varphi^s$ be the  common continuous weights for $T^u$ and
$T^s$ supplied by \cite[Proposition~2.11]{LTTTY}.  Define
$\mathcal B_n^u,\mathcal B_n^s$ by the existence of the
corresponding backward Birkhoff averages, and define
$\mathcal B^u,\mathcal B^s$ analogously.   These sets have
full measure by Birkhoff's theorem. Set
\begin{align*}
 \mathcal G_n
 &:=\Gamma_k(g_n)\cap\Gamma_{d-k}(g_n^{-1})
       \cap\mathcal B_n^u\cap\mathcal B_n^s
       \cap\operatorname{NUH}'(g_n),\\
 \mathcal G
 &:=\Gamma_k(g)\cap\Gamma_{d-k}(g^{-1})
       \cap\mathcal B^u\cap\mathcal B^s
       \cap\operatorname{NUH}'(g).
\end{align*}
The uniform exponent gaps and
\cite[Theorem~2.4]{SPR2025} imply
\[
 \rho_n(\mathcal G_n)=1,\qquad \rho(\mathcal G)=1.
\]
Moreover, positive iterates preserve the relevant Lyapunov
regularity conditions. Thus the conditions in
\cite[Lemma~2.10]{LTTTY} are satisfied: on
$A_n^u\cap\mathcal G_n$ and $A^u\cap\mathcal G$, the global
forward weighted unstable manifolds coincide with the
classical unstable manifolds; on the corresponding inverse
blocks, the global weighted unstable manifolds coincide
with the classical stable manifolds. In particular, the
weighted local disks are classical Pesin manifold pieces
at all points of these good sets. At other points of the
blocks they remain well defined.

Choose the forward and inverse blocks with error $\delta/4$ and define
\begin{equation}\label{eq:Lambda-n-regular-closure}
 \Lambda_n:=\overline{A_n^u\cap A_n^s\cap\mathcal G_n}.
\end{equation}
Then $\rho_n(\Lambda_n)>1-\delta/2$.  Let
\[
 \Lambda:=\limsup_{n\to\infty}\Lambda_n.
\]
By \eqref{eq:forward-weighted-upper-limit}--
\eqref{eq:inverse-weighted-upper-limit}, this compact set is contained in
$A^u\cap A^s$.  It also satisfies $\rho(\Lambda)>1-\delta$.

Finally, $\Lambda_n$ and $\Lambda$ are contained in the
corresponding compact weighted blocks. The uniform graph
transform estimates of \cite[Proposition~2.15]{LTTTY}
therefore give weighted local stable and unstable manifolds
of a common positive size and their $C^1$ convergence
whenever $x_n\in\Lambda_n$ and $x_n\to x\in\Lambda$.
Together with the preceding identification on the full
measure good sets, this proves all the assertions.

\end{proof}

The next lemma keeps track of the domains used in a finite orbit study
of a homoclinic intersection.  A uniform lower bound on the size of local
disks at the original block alone would not suffice, since their sizes
can decrease along the orbit saturation.

\begin{lemma}[Finite orbit saturations and controlled exhaustion]
\label{lem:finite-weighted-saturations}
Use the blocks and iterates of \cref{lem:uniform-two-sided-blocks}.
For $\sigma\in\{u,s\}$ and $L\ge0$, set
\[
 A_{n,L}^\sigma:=\bigcup_{|j|\le L}(T_n^\sigma)^jA_n^\sigma,
 \qquad A_L^\sigma:=\bigcup_{|j|\le L}(T^\sigma)^jA^\sigma.
\]
There are numbers $r_L^\sigma>0$, independent of $n$, and restrictions
$D_{n,L}^\sigma(z)$, $D_L^\sigma(z)$ of the local weighted disks on these
sets, with domain radius $r_L^\sigma$ in their weighted graph coordinates,
such that the following hold.
For each fixed $L$, if $x_n\in A_n^\sigma$, $x_n\to x\in A^\sigma$,
and $|j|\le L$, then
\begin{equation}\label{eq:finite-saturation-Cone}
 D_{n,L}^\sigma((T_n^\sigma)^jx_n)
 \longrightarrow D_L^\sigma((T^\sigma)^jx)
 \quad\text{in }C^1.
\end{equation}
The disks are also $C^1$-continuous within each limit finite saturation.
Moreover, for $x\in A^\sigma$,
\begin{equation}\label{eq:controlled-weighted-exhaustion}
 W^{\sigma,\mathrm w}(x,g)
 =\bigcup_{m\ge0}(T^\sigma)^m
           \operatorname{int}D_m^\sigma((T^\sigma)^{-m}x).
\end{equation}
Here $W^{\sigma,\mathrm w}(x,g)$ denotes the global weighted unstable
manifold for $T^\sigma$, and the interior is relative to the disk.
For $\sigma=s$, this is a stable manifold for $g$.
\end{lemma}

\begin{proof}
Treat one direction, writing $T_n,T,A_n,A$ and omitting $\sigma$.
Let $K_b,\chi_b,\epsilon_b$ be the common weighted-block parameters, and
fix the H\"older exponent $0<\alpha<1$. 
Let $\Upsilon\ge2$ be a common $C^1$ bound for
$T_n^{\pm1}$ and $T^{\pm1}$.

As in \cite[Equation~(4)]{LTTTY}, choose
\[
 \tau_b=\frac{8(\log\Upsilon+1)}{\chi_b},
 \qquad
 \epsilon_b\le
 \frac{\alpha\chi_b^2}{1000(\log\Upsilon+1)}.
\]
Since $\chi_b\le\log\Upsilon$, these choices give
\begin{equation}\label{eq:graded-radius-gap}
 b_b:=\frac{4\tau_b\epsilon_b}{\alpha}
 \le \frac{32}{1000}\chi_b
 <\chi_b-3\epsilon_b.
\end{equation}

A shift of at most $L$ replaces the block constant by at most
$K_b e^{L\epsilon_b}$.  The weighted chart radius estimate
\cite[Equation~(10) and the estimate preceding
Proposition~2.15]{LTTTY} therefore gives a uniform domain
containing $B(0,2r_*e^{-b_bL})$, for some $r_*>0$
independent of $n$ and $L$. Set
\begin{equation}\label{eq:graded-radii}
 r_L:=r_*e^{-b_bL}
\end{equation}
and restrict the local weighted disks to this radius.
The uniform graph-transform estimates of
\cite[Proposition~2.15]{LTTTY}, applied to the finitely many
shifts, yield \eqref{eq:finite-saturation-Cone} and the
asserted $C^1$-continuity within each limit finite
saturation.

It remains to check that this particular restriction does not lose global
points.  Let $x\in A$ and let $y$ belong to its global weighted unstable
manifold for $T$.
By \cite[Proposition~2.9]{LTTTY}, there exists $m_0\ge0$ such that
$T^{-m_0}y$ lies in the original local weighted unstable disk through
$T^{-m_0}x$.  Local backward invariance and the contraction estimate in
\cite[Appendix~A.3, proof of Proposition~2.9]{LTTTY} imply, in the
weighted chart $\Psi_{T^{-m}x}$, that for $m\ge m_0$,
\begin{equation}\label{eq:graded-witness-contraction}
 \bigl|\Psi_{T^{-m}x}^{-1}(T^{-m}y)\bigr|_{T^{-m}x}
 \le C_{x,y} e^{-(\chi_b-3\epsilon_b)(m-m_0)}.
\end{equation}
Here $|\cdot|_z$ is the weighted coordinate norm.  
Since $b_b<\chi_b-3\epsilon_b$ and
$r_m=r_*e^{-b_bm}$, the ratio of the right hand side
to $r_m$ tends to zero. Hence
\[
 T^{-m}y\in
 \operatorname{int}D_m(T^{-m}x)
\]
for all sufficiently large $m$. This proves one inclusion
in \eqref{eq:controlled-weighted-exhaustion}. The reverse
inclusion follows from the invariance of the global
weighted unstable manifold.

Applying the same argument to $T^u$ and $T^s$ proves
the lemma.   
\end{proof}

\begin{lemma}[Robust  homoclinic relation]
\label{lem:robust-finite-data-homoclinic}
Assume the setting and notation of
\cref{lem:uniform-two-sided-blocks}.  Let $p$ be a hyperbolic periodic point
of $g$, and let $p_n$ be its continuation for $g_n$.  Let
$K\subset\Lambda\cap\mathcal G$ be compact, and suppose that every
$x\in K$ is homoclinically related to $p$.  Then there are an open
neighborhood $V$ of $K$ and $n_0\ge1$ such that, for every $n\ge n_0$ and
every
\[
 x_n\in V\cap\Lambda_n\cap\mathcal G_n,
\]
the point $x_n$ is homoclinically related to $p_n$.
\end{lemma}

\begin{proof}
Fix $x\in K$.  Since $x$ is homoclinically related to $p$,
there exist integers $a,b$ such that
\[
 W^u_g(x)\pitchfork W^s_g(g^ap)\ne\varnothing,
 \qquad
 W^s_g(x)\pitchfork W^u_g(g^bp)\ne\varnothing.
\]
Fix these two phases for the argument at $x$.
By \cite[Lemma~2.10]{LTTTY}, these are intersections of the
corresponding global weighted manifolds. By
\cref{lem:finite-weighted-saturations}, both intersections
are realized by finite iterates of local weighted disks.
Choose compact subdisks around the intersection points.

By the $C^1$-continuity of the finite weighted disks within
$\Lambda$, there is an open neighborhood
$O_x$ of $x$ in $M$ such that the two transverse
intersections persist for every $y\in O_x\cap \Lambda$.

Choose $V_x$ with $x\in V_x$ and $\overline{V_x}\subset O_x$.  We claim
that there is $n_x$ such that, for every $n\ge n_x$, every
$y_n\in V_x\cap\Lambda_n\cap\mathcal G_n$ is homoclinically related to
$p_n$.  Otherwise there are $n_j\to\infty$ and such points $y_{n_j}$ for
which one of the two intersections fails.  After extraction,
$y_{n_j}\to y\in\overline{V_x}\cap\Lambda$.  The finite-saturation
convergence from \cref{lem:finite-weighted-saturations}, together with the
$C^1$ convergence of the local stable manifolds of $g_{n_j}^ap_{n_j}$
and the local unstable manifolds of $g_{n_j}^bp_{n_j}$ to those of
$g^ap$ and $g^bp$, respectively, preserves the two compact transverse
intersections for all large $j$, a contradiction.

Since $y_n\in\mathcal G_n$, its global weighted
manifolds coincide with the classical Pesin manifolds by
\cite[Lemma~2.10]{LTTTY}. Thus the two intersections with the
corresponding periodic phases give the homoclinic relation in the sense of
\cite[Definition~2.5]{SPR2025}.

A finite family $V_{x_1},\ldots,V_{x_m}$ covers the compact set $K$.  Set
\[
 V:=\bigcup_{i=1}^mV_{x_i},
 \qquad
 n_0:=\max_{1\le i\le m}n_{x_i},
\]
where $n_{x_i}$ is the threshold obtained above.  The conclusion follows.
\end{proof}

\begin{lemma}[Homoclinic capture under perturbation]
\label{lem:homoclinic-capture}
In the setting of \cref{lem:uniform-two-sided-blocks}, suppose that
\begin{equation}\label{eq:component-positive-weight}
 \rho=\alpha\eta+(1-\alpha)\rho',
 \qquad 0<\alpha\le1,
\end{equation}
where $\eta$ is an ergodic measure of maximal entropy for $g$ and $\rho'$ is
$g$-invariant.  Let $p$ be a hyperbolic periodic point in the Borel
homoclinic class carrying $\eta$, and let $p_n$ be its hyperbolic
continuation for $g_n$.  Then, for all sufficiently large $n$, the Borel
homoclinic class of $p_n$ carries $\rho_n$.
\end{lemma}

\begin{proof}
Choose $0<\delta<\alpha/100$ and the blocks and full-measure sets in
\cref{lem:uniform-two-sided-blocks}.  From
\eqref{eq:component-positive-weight} and $\rho(\Lambda)>1-\delta$ we obtain
\[
 \eta(\Lambda)>1-\delta/\alpha>99/100.
\]
Since $\rho(\mathcal G)=1$ and $\rho\ge\alpha\eta$, one has
$\eta(\mathcal G)=1$.  The measure $\eta$ is carried by the Borel
homoclinic class of $p$ \cite[Proposition~2.8]{SPR2025}.  Hence inner
regularity gives a compact set
\begin{equation}\label{eq:compact-homoclinic-core}
 K\subset\Lambda\cap\mathcal G
       \cap\{x\in\operatorname{NUH}'(g):x\sim_g p\}
 \quad\text{with}\quad
 \eta(K)>2/3.
\end{equation}
Apply \cref{lem:robust-finite-data-homoclinic} to obtain an open
neighborhood $V$ of $K$.  Since $\rho\ge\alpha\eta$,
\[
 \rho(V)\ge\alpha\eta(V)>2\alpha/3.
\]
By Portmanteau and \eqref{eq:uniform-block-masses}, for all sufficiently
large $n$,
\[
 \rho_n(V\cap\Lambda_n)
 \ge \rho_n(V)-\rho_n(M\setminus\Lambda_n)>0.
\]
Because $\rho_n(\mathcal G_n)=1$, the set
$V\cap\Lambda_n\cap\mathcal G_n$ also has positive $\rho_n$-measure.  Every
point of this set is homoclinically related to $p_n$ by
\cref{lem:robust-finite-data-homoclinic}.  The Borel homoclinic class of
$p_n$ is invariant \cite[Proposition~2.8]{SPR2025}; thus ergodicity of $\rho_n$
 implies that this class carries $\rho_n$.
\end{proof}

\subsection{Stability of maximal measures}

\begin{proposition}[Upper semicontinuity of maximal multiplicity]
\label{prop:MME-number-usc}
The function $g\mapsto n_{\mathrm{MME}}(g)$ is upper semicontinuous on $U$ in
the $C^\infty$ topology.
\end{proposition}

\begin{proof}
Fix $g\in U$ and write its ergodic maximal measures as
$\eta_1,\ldots,\eta_q$; the number $q$ is finite by
\cref{prop:stable-SPR}.  Suppose, to the contrary, that there are
$g_n\to g$ for which $g_n$ has at least $q+1$ distinct ergodic maximal
measures
\[
 \rho_n^1,\ldots,\rho_n^{q+1}.
\]
After a diagonal extraction, assume $\rho_n^j\to\rho^j$ for every $j$.
By \cref{prop:entropy-continuity-neighborhood} and the joint upper semicontinuity in Theorem \ref{thm:joint-usc},
each $\rho^j$ is a maximal measure for $g$.  Hence
\begin{equation}\label{eq:limit-MME-convex-combination}
 \rho^j=\sum_{i=1}^q \alpha_{j,i}\eta_i,
 \qquad \alpha_{j,i}\ge0,
 \qquad \sum_i \alpha_{j,i}=1.
\end{equation}
For each $j$, choose $i(j)$ with $\alpha_{j,i(j)}>0$.  By the pigeonhole
principle, there are $j\ne\ell$ and an index $i$ such that
$\alpha_{j,i},\alpha_{\ell,i}>0$.

Choose a hyperbolic periodic point $p$ in the Borel homoclinic class of
$\eta_i$, and let $p_n$ be its continuation.  Applying
\cref{lem:homoclinic-capture} to the two sequences shows that, for all large
$n$, both $\rho_n^j$ and $\rho_n^\ell$ are carried by the Borel homoclinic
class of $p_n$.  Let $X_n$ denote this Borel homoclinic class.  Since
$\rho_n^j(X_n)=1$, the variational definition gives
\[
 h_{\rm top}(g_n|X_n)
 \ge h_{\rho_n^j}(g_n)=h_{\rm top}(g_n),
\]
whereas the reverse inequality is automatic.  Thus
$h_{\rm top}(g_n|X_n)=h_{\rm top}(g_n)$.  By \cite[Proposition 2.14]{SPR2025} (see also 
\cite{BCSMME}), the class carries at most one measure
attaining the topological  entropy.  Therefore
$\rho_n^j=\rho_n^\ell$, contradicting their choice.  Thus
$n_{\mathrm{MME}}(g_n)\le q$ for all sufficiently large $n$.
\end{proof}

\begin{proof}[Proof of \cref{thm:fixed-dimension-variant}]
By \cref{prop:stable-SPR,prop:entropy-continuity-neighborhood,prop:MME-number-usc},
conclusions \textnormal{(1)}--\textnormal{(3)} of \cref{thm:main-intro}
hold on $U$.  The sequential assertion of conclusion \textnormal{(4)}
is \cref{prop:fixed-index-maximal-continuity}.

Let
\[
 X_U:=\{(\mu,g):g\in U,\ \mu\in\mathcal M(g),\
                      h_\mu(g)=h_{\mathrm{top}}(g)\}.
\]
If $(\mu_n,g_n)\in X_U$ converges to $(\mu,g)\in\mathcal P(M)\times U$,
then entropy continuity gives
\[
 h_{\mu_n}(g_n)=h_{\mathrm{top}}(g_n)\to h_{\mathrm{top}}(g).
\]
By \cref{prop:fixed-index-maximal-continuity}, $(\mu,g)\in X_U$ and
$S^+(\mu_n,g_n)\to S^+(\mu,g)$.  Thus $X_U$ is relatively closed, and
$S^+$ is continuous on it.

Suppose now that $f\in U$ has a unique maximal measure.  Upper
semicontinuity of $n_{\mathrm{MME}}$ gives a $C^\infty$ neighborhood
$U_f\subset U$ on which $n_{\mathrm{MME}}(g)\le1$.  Existence of a maximal
measure gives equality; denote this measure by $\mu_g$.  If $g_n\to g$ in
$U_f$, every weak limit $\mu$ of $\mu_{g_n}$ satisfies $(\mu,g)\in X_U$
and hence equals $\mu_g$.  Compactness gives $\mu_{g_n}\to\mu_g$.  The identity
$h_{\mu_g}(g)=h_{\mathrm{top}}(g)$ and the continuity statements just
proved give the other two assertions in conclusion \textnormal{(5)}.
\end{proof}

\subsection{The proof of Theorem B}

\begin{proof}[Proof of \cref{thm:main-intro}] Let $f\in \mathcal S^\infty(M)$. It follows from the proof of Corollary  \ref{cor:openess} that there is $\chi>0$, an integer $k\in \{1,\cdots,d-1\}$ and a $C^1$ neighborhood $U$ of $f$ such that  any $g\in U\cap \mathrm{Diff}^\infty(M)$ lies in $\mathcal S^\infty(M)$ and satisfies
\begin{itemize}
\item $\lambda_{\mathrm{sec}}(g)<\chi + \lambda_{\mathrm{sec}}(f)$,
\item $ h_{\mathrm{top}}(g)> h_{\mathrm{top}}(f)-\chi$,
\item $\lambda_k(g)\geq h_{\mathrm{top}}(g)>4\chi+ \lambda_{\mathrm{sec}}(g)> \lambda_{k-1}(g)$.
\end{itemize}
Let $h_*=2\chi+\lambda_{\mathrm{sec}}(f) <\inf_{g\in U}h_{\mathrm{top}}(g)$. Then for any $g\in U$ and for any $\rho\in \mathcal M_e(g)$ with $h_\rho(g)>h_*$ we have $h_\rho(g)>\lambda_{\mathrm{sec}}(g)+\chi$ and therefore $\chi_k(\rho,g),-\chi_{k+1}(\rho,g)>\chi$ by Lemma \ref{lem:index}.  
We may therefore conclude the proof of \cref{thm:main-intro} by applying   \cref{thm:fixed-dimension-variant}. 
\end{proof}


\section{Examples and applications}
\label{sec:examples}

We first prove the two sharpness examples from
\cref{subsec:intro-sharpness}, keeping the diffeomorphism fixed while
varying the invariant measure. We then establish the applications to
Newhouse domains and products of standard maps. Throughout this section,
we extend \eqref{eq:def-lambda-fi} to the top exterior degree by the same
formula; in particular, $\lambda_d(f)=0$ for volume-preserving maps.

The following standard facts will be used in both examples and in the
standard-map application. Recall the cat map $A$ and its entropy $h_A$
from \cref{subsec:intro-sharpness}.

\begin{lemma}[Products and the cat map]
\label{lem:sharp-product}
For diffeomorphisms $f_i:M_i\to M_i$ of closed manifolds,
\begin{align*}
 \lambda_q(f_1\times f_2)
 &=\max_{i+j=q}\{\lambda_i(f_1)+\lambda_j(f_2)\},\\
 h_{\mathrm{top}}(f_1\times f_2)
 &=h_{\mathrm{top}}(f_1)+h_{\mathrm{top}}(f_2),\\
 h_{\nu_1\times\nu_2}(f_1\times f_2)
 &=h_{\nu_1}(f_1)+h_{\nu_2}(f_2).
\end{align*}
Here the maximum ranges over the admissible exterior degrees, and the
$\nu_i$ are invariant probabilities. If $\nu_1$ is weakly mixing and
$\nu_2$ is ergodic, their product is ergodic. Moreover, Haar probability
$m_A$ is mixing for every $A^m$, $m\ge1$, and
\begin{gather*}
 h_{\mathrm{top}}(A^m)=h_{m_A}(A^m)=mh_A,\\
 (\lambda_0(A^m),\lambda_1(A^m),\lambda_2(A^m))=(0,mh_A,0).
\end{gather*}
Its Lyapunov exponents for $A^m$ are $\pm mh_A$.
\end{lemma}

\begin{proof}
In the orthogonal product metric, the decomposition
\[
 \wedge^q(E_1\oplus E_2)
 =\bigoplus_{i+j=q}\wedge^iE_1\otimes\wedge^jE_2
\]
gives the exterior formula by taking norms, suprema, and exponential
growth rates. The entropy and mixing assertions are standard
\cite{Walters}; the cat-map identities follow from its eigenvalues
$e^{\pm h_A}$ and the entropy formula for toral automorphisms
\cite{BowenEntropy,KatokHasselblatt}.
\end{proof}

\subsection{A uniform-gap discontinuity at the metric threshold}
\label{subsec:sharp-measure-proof}

\begin{proof}[Proof of \cref{prop:sharp-measure}]
We first construct periodic measures of an area-preserving surface map
whose exponents have a uniform gap but are discontinuous at their weak
limit. A power of the cat map then places their entropy at the required
threshold.

\medskip
\noindent\emph{Conservative realization.}
Fix $\sigma>1$, $r>1$, $\alpha=r^{-1}$ and $d>0$, and consider
\begin{align}
 T_0(x,y)&=(\sigma^{-1}x,\sigma y),
 \label{eq:sharp-T0}\\
 T_1(x,1+v)&=(1-rv,\alpha x+dv^2).
 \label{eq:sharp-T1}
\end{align}
Both maps have determinant one. Choose $\sigma^{-1}<\beta<1$ and
$\eta>0$ sufficiently small, and let $D_0$ be a smoothly rounded disk with
\begin{equation}\label{eq:sharp-D0-sandwich}
 [-\eta,1+\eta]\times[-\eta,\beta-2\eta]
 \subset\operatorname{int}D_0
 \subset[-2\eta,1+2\eta]\times[-2\eta,\beta-\eta].
\end{equation}
Then $(0,1)\in\operatorname{int}T_0(D_0)$, whereas a neighborhood of
$(1,0)$ lies in $\operatorname{int}D_0\setminus T_0(D_0)$. A sufficiently
small disk $D_1$ about $(0,1)$ therefore satisfies
\begin{equation}\label{eq:sharp-disk-configuration}
 D_1\subset\operatorname{int}T_0(D_0)\cap\{y>\beta\},
 \qquad
 T_1(D_1)\subset\operatorname{int}D_0\setminus T_0(D_0).
\end{equation}
In particular, the source disks $D_0,D_1$ are disjoint, as are their
images $T_0(D_0),T_1(D_1)$. Place all four
disks in a common coordinate chart of $S^2$.
The prescribed local maps extend to an orientation-preserving
diffeomorphism $\widetilde T\in\Diff^\infty(S^2)$ by
extending their boundary collar maps across the complementary
annulus \cite[Chapter~8]{Hirsch}.

Choose an area form $\omega$ equal to $dx\wedge dy$ near
all four disks. Since the local maps preserve $dx\wedge dy$,
the form $\widetilde T^*\omega-\omega$ is supported in
the complementary annulus and has integral zero.
It therefore admits a compactly supported primitive there.
The relative Moser argument gives a diffeomorphism $h$,
equal to the identity near $D_0\cup D_1$, such that
\[
 h^*\widetilde T^*\omega=\omega
 \]
\cite{Moser,Cannas}. Thus $T:=\widetilde T\circ h$
preserves $\omega$ and agrees with $T_i$ near $D_i$.
In particular, the prescribed saddle at $q=(0,0)$ is
unchanged and has Lyapunov exponents $\pm\gamma$,
where $\gamma=\log\sigma$.

\medskip
\noindent\emph{Periodic measures and their exponents.}
Put $\varepsilon_n=\sigma^{-n}$.  Solving the fixed-point
equations for $T_0^n\circ T_1$ gives
\[
 x_n=\varepsilon_n(1-rv_n),\qquad
 dv_n^2-2\varepsilon_nv_n-(1-\alpha)\varepsilon_n=0.
\] 
Taking the positive root, we obtain
\begin{equation}\label{eq:sharp-v-asymptotic}
 v_n=\frac{\varepsilon_n+
 \sqrt{\varepsilon_n^2+d(1-\alpha)\varepsilon_n}}{d}
 =\sqrt{\frac{1-\alpha}{d}}\,\sigma^{-n/2}
 +O(\sigma^{-n}).
\end{equation} 
For all large $n$, $z_n=(x_n,1+v_n)\in\operatorname{int}D_1$,
and its successive iterates are
\begin{equation}\label{eq:sharp-local-orbit-points}
 T_0^jT_1(z_n)=
 \bigl(\sigma^{-j}(1-rv_n),\sigma^{j-n}(1+v_n)\bigr),
 \qquad 0\le j\le n.
\end{equation}
For $j<n$, these points belong to $\operatorname{int}D_0$,
since their first coordinates lie in $(0,1)$ and their second
coordinates are bounded by
$\sigma^{-1}(1+v_n)<\beta-2\eta$. At $j=n$ the orbit returns
to $z_n$. Thus $z_n$ has exact period $n+1$ for $T$.

Let $\nu_n$ be its periodic probability. Given a neighborhood $U$ of
$q$, choose $J$ large. Formula \eqref{eq:sharp-local-orbit-points} puts
all points with $J\le j\le n-J$ in $U$, for every sufficiently large
$n$. Hence $\nu_n(U^c)\le(2J+1)/(n+1)$ and
$\nu_n\to\delta_q$.

The derivative of the full return is
\begin{equation}\label{eq:sharp-monodromy}
 M_n=D(T_0^n\circ T_1)(z_n)
 =\begin{pmatrix}
 0&-r\sigma^{-n}\\
 \alpha\sigma^n&2dv_n\sigma^n
 \end{pmatrix}.
\end{equation}
Its determinant is one and, by \eqref{eq:sharp-v-asymptotic},
\[
 \operatorname{tr}M_n
 =2\sqrt{d(1-\alpha)}\,\sigma^{n/2}(1+o(1)).
\]
For large $n$, its eigenvalues are positive and reciprocal, with the
larger one asymptotic to its trace. Dividing its logarithm by the period
$n+1$ gives
\begin{equation}\label{eq:sharp-second-periodic-exponents}
 \operatorname{Lyap}(\nu_n,T)=\{\gamma_n,-\gamma_n\},
 \qquad \gamma_n\longrightarrow\frac\gamma2.
\end{equation}
Discarding finitely many terms ensures $\gamma_n\ge\gamma/4$.

\medskip
\noindent\emph{Product with a cat map.}
Choose a metric on $S^2$ with area form $\omega$, and set
$L=\lambda_1(T)$. Then $L\ge\gamma>0$ and $\lambda_2(T)=0$.
Choose $m\ge1$ with $H=mh_A>L$, put $B=A^m$, and let $m_B$ be Haar
probability. Define
\begin{equation}\label{eq:sharp-second-product}
 F=B\times T,
 \qquad \mu_n=m_B\times\nu_n,
 \qquad \mu=m_B\times\delta_q.
\end{equation}
The map $F$ preserves product volume and $\mu_n\to\mu$. Since $m_B$ is
mixing, these product measures are ergodic and have entropy $H$ by
\cref{lem:sharp-product}. Their spectra are
\[
 \operatorname{Lyap}(\mu_n,F)=\{H,\gamma_n,-\gamma_n,-H\},
 \qquad
 \operatorname{Lyap}(\mu,F)=\{H,\gamma,-\gamma,-H\}.
\]
Consequently, their unstable dimension is two, all exponents have
absolute value at least $\min\{H,\gamma/4\}$, and
\begin{equation}\label{eq:sharp-second-Splus-jump}
 S^+(\mu_n,F)\longrightarrow H+\frac\gamma2
 <H+\gamma=S^+(\mu,F).
\end{equation}
Finally, the exterior-growth formula and
\[
 (\lambda_0(B),\lambda_1(B),\lambda_2(B))=(0,H,0),
 \qquad
 (\lambda_0(T),\lambda_1(T),\lambda_2(T))=(0,L,0)
\]
give
\begin{equation}\label{eq:sharp-second-product-rates}
 \lambda_1(F)=H,\qquad
 \lambda_2(F)=H+L,\qquad
 \lambda_3(F)=H.
\end{equation}
Thus $\lambda_{\mathrm{sec}}(F)=H=h_{\mu_n}(F)=h_\mu(F)$, completing
the proof.
\end{proof}

\subsection{Critical equality for maximal measures}
\label{subsec:sharp-main-proof}

This example uses a zero-entropy surface map whose periodic saddle
measures converge to a hyperbolic fixed-point measure while their
exponents tend to zero.

\begin{proof}[Proof of \cref{prop:sharp-main}]
\emph{Resonant cycles near a separatrix.}
Fix $c\in(0,h_A)$ and choose
\[
 (c/h_A)^2<\alpha<1,
 \qquad
 \tau=\frac{c}{4\pi^2\sqrt\alpha}.
\]
Consider the smooth vector field on $(\mathbb T^2,dx\wedge dy)$
given by
\[
 X(x,y)=
 \bigl(2\pi\alpha\sin(2\pi y),
       -2\pi\sin(2\pi x)\bigr).
\]
It is area-preserving and admits the first integral
\[
 E(x,y)=\cos(2\pi x)+\alpha\cos(2\pi y).
\]
Let $\phi^t$ be its flow and set $g_0=\phi^\tau$.
The equilibrium $q=(0,\tfrac12)$ has energy $E_0=1-\alpha$,
and
\[
 DX(q)=
 \begin{pmatrix}
 0&-4\pi^2\alpha\\
 -4\pi^2&0
 \end{pmatrix}.
\]
Consequently,
\begin{equation}\label{eq:sharp-q-spectrum-g0}
 \operatorname{Lyap}(\delta_q,g_0)=\{c,-c\},
 \qquad
 \|Dg_0\|_\infty\le e^{4\pi^2\tau}
 =e^{c/\sqrt\alpha}<e^{h_A}.
\end{equation}
The derivative bound follows from $\|DX\|_\infty\le4\pi^2$
and the variational equation.

For $E\in(E_0,1+\alpha)$, let $\gamma_E$ be the level curve
surrounding $(0,0)$ and let $T(E)$ be its primitive flow
period. The period is continuous on this regular energy interval
and tends to infinity as $E\downarrow E_0$. Indeed, writing
$\varepsilon=E-E_0$, the energy equation gives
\[
 |\dot x|
 =2\pi\sqrt{\alpha^2-(E-\cos(2\pi x))^2}
 \le C\sqrt{\varepsilon+x^2}
\]
near $q$. Hence a local traversal time is bounded below by
\[
 C^{-1}\int_{-\delta}^{\delta}
       (\varepsilon+x^2)^{-1/2}\,dx,
\]
which diverges as $\varepsilon\downarrow0$. Away from a fixed
neighborhood of $q$, the nearby level arcs form a compact smooth
family and the vector field is bounded away from zero, so the
remaining excursion time is uniformly bounded.

By continuity of $T$, choose $E_n\downarrow E_0$ and integers
$N_n\to\infty$ such that $T(E_n)=N_n\tau$. For
$x_n\in\gamma_{E_n}$, the points
\begin{equation}\label{eq:sharp-sample-orbit}
 p_{n,j}=\phi^{j\tau}(x_n),\qquad 0\le j<N_n,
\end{equation}
form an exact $N_n$-cycle of $g_0$. Their uniform probabilities
converge to $\delta_q$. Indeed, given a neighborhood $U$ of $q$,
choose $V$ with $q\in V$ and $\overline V\subset U$.
For sufficiently small $\delta_0>0$, a sample point outside
$U$ remains outside $V$ for times $|t|\le\delta_0$.
The corresponding intervals on the time circle of length
$N_n\tau$ are disjoint if $2\delta_0<\tau$. Since the time
spent by $\gamma_{E_n}$ outside $V$ is uniformly bounded,
the number of sample points outside $U$ is uniformly bounded.
Thus the periodic probabilities converge to $\delta_q$.

Let $\Gamma$ be the separatrix component of $E^{-1}(E_0)$
containing $q$. The curves $\gamma_E$ converge to $\Gamma$ in Hausdorff
distance: regular points are approached by the implicit function theorem,
and the energy equation on $x=0$ gives points tending to $q$.
Choose pairwise disjoint closed energy annuli $\mathcal A_n$ around
$\gamma_{E_n}$, accumulating only on $\Gamma$, such that
\begin{equation}\label{eq:sharp-rho-annuli}
 \rho_n:=\dist(\mathcal A_n,\Gamma)>0,
 \qquad
 \sup_{x\in\mathcal A_n}\dist(x,\Gamma)\longrightarrow0.
\end{equation}

\smallskip
\noindent\emph{Hyperbolization and smooth gluing.}
Choose symplectic action--angle coordinates
\[
 \Psi_n:[-\ell_n,\ell_n]\times\mathbb T
       \longrightarrow \mathcal A_n
\]
on a slightly larger annulus, with
\[
 g_0(I,\theta)=(I,\theta+\omega_n(I)),
 \qquad
 \omega_n(0)=\frac1{N_n},
 \qquad
 \Psi_n(0,j/N_n)=p_{n,j}.
\]
Here $\omega_n(I)=\tau/T(E(I))$.
Set
\[
 R_n(I,\theta)=(I,\theta+1/N_n),
 \qquad
 K_{0,n}(I)=\int_0^I
       \bigl(\omega_n(s)-1/N_n\bigr)\,ds.
\]
In these coordinates, $g_0$ is the composition of $R_n$
with the time-one map of the vector field
\[
 Y_{0,n}(I,\theta)=(0,K_{0,n}'(I)).
\]

Choose a cutoff $\chi_n(I)$ equal to one near zero and zero
near the boundary, and put
\[
 V_n(I,\theta)
 =\epsilon_n\chi_n(I)I\sin(2\pi N_n\theta),
 \qquad \epsilon_n>0.
\]
Set $K_n=K_{0,n}+V_n$ and consider the area-preserving
vector field
\[
 Y_n(I,\theta)=
 \bigl(-\partial_\theta K_n,\partial_I K_n\bigr).
\]
It commutes with $R_n$ and has hyperbolic equilibria
$z_{n,j}=(0,j/N_n)$, $0\le j<N_n$.
Let $\Phi_n^t$ denote its flow and define
\[
 G_n=R_n\circ\Phi_n^1.
\]
Since $R_n^{N_n}=\id$ and $R_n$ commutes with the flow,
\begin{equation}\label{eq:sharp-Gn-power}
 G_n^{N_n}=\Phi_n^{N_n},
 \qquad
 G_n(z_{n,j})=z_{n,j+1}\quad(\bmod N_n).
\end{equation}
At each protected equilibrium,
\[
 DY_n(z_{n,j})=
 \begin{pmatrix}
 -\kappa_n&0\\
 K_{0,n}''(0)&\kappa_n
 \end{pmatrix},
 \qquad
 \kappa_n=2\pi N_n\epsilon_n>0.
\]
Thus the derivative of the $N_n$-return has eigenvalues
$e^{\pm N_n\kappa_n}$, and the corresponding periodic
measure has Lyapunov exponents $\pm\kappa_n$.

Write $g_n=\Psi_nG_n\Psi_n^{-1}$. For each fixed $n$,
finite-time flows and their inverses depend continuously in
$C^r$ on the vector field in $C^r$. Since $\epsilon_n$ can
be chosen arbitrarily small, we may impose
$0<\kappa_n<1/n$ and, in a fixed atlas on $\mathbb T^2$,
\begin{equation}\label{eq:sharp-flat-estimate}
 \|g_n-g_0\|_{C^r(\mathcal A_n)}
 +\|g_n^{-1}-g_0^{-1}\|_{C^r(\mathcal A_n)}
 \le\rho_n^n,
 \qquad 0\le r\le n.
\end{equation}
The smallness is imposed after conjugation by each fixed
$\Psi_n$. We may also require
\[
 \|Dg_n-Dg_0\|_\infty
 <\frac12(e^{h_A}-\|Dg_0\|_\infty)
 \quad\text{on every }\mathcal A_n.
\]

Define
\begin{equation}\label{eq:sharp-global-g}
 g=g_n\text{ on }\mathcal A_n,
 \qquad
 g=g_0\text{ elsewhere}.
\end{equation}
Since $\chi_n=0$ near the boundary, each $g_n$ agrees
with $g_0$ on a boundary collar and preserves $\mathcal A_n$.
Thus $g$ is bijective, with the corresponding piecewise
inverse, and is smooth off $\Gamma$.

For every multi-index $\beta$ and every $m\ge1$,
\eqref{eq:sharp-flat-estimate} gives, near $\Gamma$,
\[
 |D^\beta(g-g_0)(x)|
 \le \dist(x,\Gamma)^m.
\]
Indeed, on $\mathcal A_n$ with $n\ge\max\{|\beta|,m\}$
and $\dist(x,\Gamma)<1$, we have
\[
 \rho_n^n\le\rho_n^m
 \le\dist(x,\Gamma)^m.
\]
The finitely many remaining annuli stay away from $\Gamma$.
Extending these derivatives by zero therefore proves
smoothness and flatness along $\Gamma$. The same argument
applies to the inverse. Hence $g$ is a $C^\infty$
area-preserving diffeomorphism with the same infinite jet
as $g_0$ on $\Gamma$.

\smallskip
\noindent\emph{Zero entropy and the product.}
By Young's theorem \cite{Young1977}, every continuous flow
on a compact surface has zero topological entropy. Since $g_0$ and each $\Phi_n^{N_n}|_{\mathcal A_n}$
are such time maps (extending the latter to a closed surface
if necessary), \eqref{eq:sharp-Gn-power} and the entropy
power rule imply that every ergodic $g$-invariant measure
has zero entropy. Thus the variational principle gives
\begin{equation}\label{eq:sharp-g-zero-entropy}
 h_{\mathrm{top}}(g)=0.
\end{equation}

The derivative bound and the unchanged saddle $q$ give
\begin{equation}\label{eq:sharp-L-range}
 c\le L:=\lambda_1(g)<h_A,
 \qquad \lambda_2(g)=0.
\end{equation}
The protected periodic orbit of $g_n$ is exactly
\eqref{eq:sharp-sample-orbit}. Its probability $\nu_n$ therefore satisfies
\begin{equation}\label{eq:sharp-nu-limit-spectrum}
 \nu_n\longrightarrow\delta_q,
 \qquad
 \operatorname{Lyap}(\nu_n,g)=\{\kappa_n,-\kappa_n\},
 \qquad \kappa_n\longrightarrow0.
\end{equation}

Finally, set
\[
 F=A\times g,
 \qquad
 \mu_n=m_A\times\nu_n,
 \qquad
 \mu=m_A\times\delta_q.
\]
The measures are ergodic because $m_A$ is mixing, and $\mu_n\to\mu$.
By \cref{lem:sharp-product} and
\eqref{eq:sharp-g-zero-entropy},
\[
 h_{\mathrm{top}}(F)=h_{\mu_n}(F)=h_\mu(F)=h_A.
\]
The factor spectra give those asserted in \cref{prop:sharp-main}, and
\[
 S^+(\mu_n,F)=h_A+\kappa_n\longrightarrow h_A
 \ne h_A+c=S^+(\mu,F).
\]
Finally, \cref{lem:sharp-product} and \eqref{eq:sharp-L-range} yield
\[
 \lambda_1(F)=h_A,
 \qquad \lambda_2(F)=h_A+L,
 \qquad \lambda_3(F)=h_A.
\]
Thus $\lambda_{\mathrm{sec}}(F)=h_A=h_{\mathrm{top}}(F)$, completing
the proof.
\end{proof}

\subsection{Newhouse domains}
\label{subsec:positive-examples}

We prove \cref{prop:newhouse-domains} by combining a controlled tangency
construction with Manning's entropy bound and Leal's persistence theorem.
A periodic saddle is \emph{sectionally dissipative} if it has one unstable
multiplier and every product of two multipliers of its period map has
modulus less than one. It is \emph{sectionally expansive} for $f$ if it is
sectionally dissipative for $f^{-1}$. 

Let $L\in\mathrm{GL}(3,\mathbb Z)$ be hyperbolic. Put
$k=\dim E_L^u$, $q=3-k$, and
\[
 h_L=\sum_{|\lambda|>1}\log|\lambda|,
 \qquad \chi=\chi_L.
\]
All eigenvalues are counted with multiplicity.
By taking a flat metric on $\mathbb T^3$,  
\begin{equation}\label{eq:newhouse-linear-separation}
 \|\wedge^qL\|=e^{h_L-\chi},
 \qquad h_L\ge2\chi.
\end{equation}
If $A=L$ for $k=1$ and $A=L^{-1}$ for $k=2$, then
$\dim E_A^u=1$ and $\|A|_{E_A^s}\|=e^{-\chi}$.
When $k=2$, moreover, $\|L|_{E_L^u}\|=e^{h_L-\chi}$.

Manning's first-homology bound \cite{Manning}, applied to $f$ when $k=1$
and to $f^{-1}$ when $k=2$, gives
\begin{equation}\label{eq:newhouse-Manning-bound}
 h_{\mathrm{top}}(f)\ge h_L
 \qquad\text{for every }f\simeq L.
\end{equation}
The following construction bounds the cost of creating a tangency
independently of $L$.

\begin{proposition}[Controlled tangency insertion]
\label{prop:newhouse-controlled-insertion}
There are universal constants $0<\beta_*<1$ and $C_*\ge1$ such that,
if $A\in\mathrm{GL}(3,\mathbb Z)$ is hyperbolic with one unstable
direction and $\beta=\|A|_{E_A^s}\|<\beta_*$ in the adapted metric,
there is $\varphi\in\Diff^\infty(\mathbb T^3)$ homotopic to $A$ with
a quadratic homoclinic tangency associated to a sectionally dissipative
fixed point, and
\begin{align}
 \|\wedge^2d\varphi\|_\infty
 &\le C_*\|\wedge^2A\|,
 \label{eq:newhouse-forward-insertion}\\
 \|d\varphi^{-1}\|_\infty
 &\le C_*\max\{1,\|A^{-1}|_{E_A^s}\|\}.
 \label{eq:newhouse-inverse-insertion}
\end{align}
\end{proposition}

\begin{proof}
Write $A=\operatorname{diag}(B,\varepsilon\Lambda)$ in the orthogonal
splitting $E_A^s\oplus E_A^u$, where $\|B\|=\beta$,
$\varepsilon\in\{-1,1\}$, and $\Lambda>1$.  Since $\Lambda|\det B|=1$, we have
$\Lambda\ge\beta^{-2}$. Fix $\sigma>4$ and put
$\kappa=\sigma/\Lambda$.

\emph{Compression.}
Choose nonnegative even bump functions $\eta_0,\eta_1$ with disjoint
supports in $(-1,1)$, with $0\le\eta_0\le1$, $\eta_0=1$ near zero,
and $\int\eta_1>0$. Define
\[
 r_\kappa'(t)=1-(1-\kappa)\eta_0(t)
 +(1-\kappa)\frac{\int\eta_0}{\int\eta_1}\eta_1(t),
 \qquad r_\kappa(0)=0.
\]
Then $r_\kappa=\id$ outside $(-1,1)$,
$r_\kappa'=\kappa$ near zero, and
$\kappa\le r_\kappa'\le C_0$ for a universal $C_0$.

Let $\zeta$ be a smooth cutoff supported in the unit stable
disk and equal to one near zero. The map
\[
 R(x,t)=\bigl(x,t+\zeta(x)(r_\kappa(t)-t)\bigr)
\]
is a compactly supported diffeomorphism isotopic to the
identity. After an isotropic rescaling, we may assume its
support lies in an arbitrarily small box. Its derivative
satisfies
\begin{equation}\label{eq:newhouse-dR}
 dR=\begin{pmatrix}I_2&0\\w&b\end{pmatrix},
 \qquad |w|\le C_0,\qquad
 \kappa\le b\le C_0,\qquad
 \|\wedge^2dR\|\le C_0,
\end{equation}
with a uniform constant $C_0$. Put $g=R\circ A$. Then
\begin{equation}\label{eq:newhouse-dg}
 dg=\begin{pmatrix}B&0\\wB&d\end{pmatrix},
 \qquad d=\varepsilon b\Lambda,\qquad |d|\ge\sigma.
\end{equation}

Fix $0<\eta<1/4$. On the unstable and stable cones
\[
 \mathcal C^u_\eta=\{|v|\le\eta|s|\},
 \qquad
 \mathcal C^s_\eta=\{|s|\le\eta|v|\},
\]
respectively, \eqref{eq:newhouse-dg} gives
\[
 \frac{|Bv|}{|wBv+ds|}
 \le\frac{\eta\beta}{\sigma-C_0\eta\beta},
 \qquad
 \frac{|d^{-1}(s-wv)|}{|B^{-1}v|}
 \le\frac{\beta(\eta+C_0)}{\sigma}.
\]
Choose $\beta_*>0$ sufficiently small, independently of
$A$, so that both bounds are less than $\eta$ and
\[
 \sigma-C_0\eta\beta_*>\sqrt{1+\eta^2},
 \qquad
 \beta_*^{-1}>\sqrt{1+\eta^2},
 \qquad
 \beta_*<\sigma^{-1}.
\]
The cone criterion then shows that $g$ is Anosov, with a
uniform lower bound $\alpha_*>0$ on the angle between its
stable and unstable directions. At the fixed point $p=0$,
\[
 dg(p)=\operatorname{diag}(B,\varepsilon\sigma).
\]
Since $\sigma\beta<1$, this fixed point is sectionally
dissipative.

\emph{Bending and homoclinic tangency.}
By the codimension-one Anosov classification
\cite{FranksAnosovTori,NewhouseCodimOne}, $g$ is topologically
conjugate to $A$. In particular, it is mixing, and the
transverse homoclinic points of $p$ are dense. Choose one,
$z$, outside the compression box, and a small ball $U$
around $z$ disjoint from that box and from every other point
of its homoclinic orbit.

In local coordinates near $z$, straighten the stable plaque
to $\{t=0\}$ and the unstable plaque to the $t$-axis.
The uniform angle bound allows these coordinates and their
inverses to have uniformly bounded first derivatives.
A fixed compactly supported local diffeomorphism, obtained
by a quarter-turn in the $(x_1,t)$-plane followed by a
quadratic shear, sends the unstable axis to
\[
 (s,0,cs^2),\qquad c>0,
\]
near the origin, while fixing the origin. Conjugating this
map by the straightening coordinates and restricting its
support to $U$ gives a diffeomorphism $S$, isotopic to the
identity, such that
\[
 \|dS\|_\infty+\|dS^{-1}\|_\infty
 +\|\wedge^2dS\|_\infty\le C_{\mathrm b},
\]
where $C_{\mathrm b}$ is universal.

Set $\varphi=S\circ g$. Since $U$ meets the homoclinic orbit
only at $z$, the stable plaque through $z$ is unchanged,
whereas the image of the unstable plaque through $g^{-1}z$
is bent into quadratic contact with it. Thus $\varphi$ has
a quadratic homoclinic tangency associated to $p$. The
homoclinic orbit and the derivative at $p$ are unchanged.

\emph{Derivative bounds.}
Since $g=R\circ A$ and the two perturbations have uniformly
bounded derivatives,
\[
 \|\wedge^2d\varphi\|_\infty
 \le C_{\mathrm b}C_0\|\wedge^2A\|.
\]

For the inverse estimate, on the compression box we have
\[
 A^{-1}dR^{-1}
 =
 \begin{pmatrix}
 B^{-1}&0\\
 -\delta w&\delta
 \end{pmatrix},
 \qquad
 \delta=\varepsilon\Lambda^{-1}b^{-1},
 \qquad |\delta|\le\sigma^{-1}.
\]
Hence its norm is bounded by
$C\max\{1,\|B^{-1}\|\}$ with a universal constant.
On the bending box, $R$ is the identity, so the inverse
derivative is bounded by
$C_{\mathrm b}\|A^{-1}\|$; outside both supports it is
simply $A^{-1}$. Since
\[
 \|A^{-1}\|
 =\max\{\|B^{-1}\|,\Lambda^{-1}\}
 \le\max\{1,\|B^{-1}\|\},
\]
we obtain
\[
 \|d\varphi^{-1}\|_\infty
 \le C_*\max\{1,\|A^{-1}|_{E_A^s}\|\}.
\]
Taking $C_*$ sufficiently large gives both required
estimates. Finally, $R$ and $S$ are isotopic to the identity,
so $\varphi\simeq A$.

\end{proof}

\begin{proof}[Proof of \cref{prop:newhouse-domains}]
Choose
\begin{equation}\label{eq:newhouse-universal-Delta}
 \Delta_{\mathrm{NH}}
 =1+\max\{-\log\beta_*,2\log C_*\}.
\end{equation}
If $\chi>\Delta_{\mathrm{NH}}$, put $h_L^0=h_L-\chi/2>0$.
Then $e^{-\chi}<\beta_*$ and
$C_*e^{h_L-\chi}<e^{h_L^0}$.

For $k=1$, apply \cref{prop:newhouse-controlled-insertion} to $A=L$.
The resulting map $F=\varphi$ has a sectionally dissipative tangency and
$\|\wedge^2dF\|_\infty<e^{h_L^0}$.
For $k=2$, apply it to $A=L^{-1}$ and take $F=\varphi^{-1}$.
Then $F$ has a sectionally expansive tangency and
\[
 \|dF\|_\infty
 \le C_*\|L|_{E_L^u}\|
 =C_*e^{h_L-\chi}<e^{h_L^0}.
\]
In either case $F\simeq L$, and there is a relative $C^1$ neighborhood
$\mathcal V_L$ of $F$ in $\Diff^\infty(\mathbb T^3)$ such that every
$f\in\mathcal V_L$ is homotopic to $L$ and
\begin{equation}\label{eq:newhouse-complementary-bound}
 \lambda_q(f)\le\log\|\wedge^qdf\|_\infty<h_L^0,
 \qquad q=3-k.
\end{equation}

Leal's theorem gives, in every $C^\infty$ neighborhood of a map with
a sectionally dissipative homoclinic tangency, an open persistence domain
in which maps with infinitely many nonhyperbolic H\'enon-like strange
attractors are dense
\cite[Theorem~2.1, Section~4.3, and proof of Theorem~A]{Leal}.
Apply this inside $\mathcal V_L$ for $k=1$, or inside
$\mathcal V_L^{-1}$ and then invert for $k=2$.
This produces the required Newhouse domain $\mathcal N_L\subset\mathcal V_L$
with the stated density of strange attractors or repellers.

For every $f\in\mathcal N_L$, Manning's bound and
\eqref{eq:newhouse-complementary-bound} imply
\[
 h_{\mathrm{top}}(f)\ge h_L>h_L^0>\lambda_q(f)
 \ge\min\{\lambda_1(f),\lambda_2(f)\}
 =\lambda_{\mathrm{sec}}(f).
\]
In fact the gap exceeds $\chi_L/2$.
By \cref{lem:index}, the only possible unstable index of an ergodic MME
is $k$, since the complementary degree $q$ lies below $h_L^0$.
Ergodic decomposition gives the same conclusion for all MMEs, and
\cref{thm:main-intro} supplies the stability assertions.
Finally $\chi_{L^m}=m\chi_L$, so every sufficiently large power of a
hyperbolic $L$ satisfies the hypothesis.
\end{proof}

\begin{remark}[Higher-dimensional extension]
The same construction extends to higher-dimensional tori
in the codimension-one setting. More precisely, for every
$d\ge3$, there exists $\Delta_{\mathrm{NH}}(d)>0$ such that
the conclusion of \cref{prop:newhouse-domains} remains valid
for hyperbolic $L\in\mathrm{GL}(d,\mathbb Z)$ satisfying
\[
 \dim E_L^u\in\{1,d-1\},
 \qquad
 \chi_L>\Delta_{\mathrm{NH}}(d).
\]
The case $\dim E_L^u=1$ gives the attracting version,
whereas the case $\dim E_L^u=d-1$ follows by applying
the construction to the inverse map and gives the
repelling version.
\end{remark}

\subsection{Products of standard maps}
\label{subsec:standard-products-proof}

\begin{proof}[Proof of \cref{thm:standard-map-products}]
Duarte's basic sets \cite[Theorem~A]{Duarte1994} are conjugate to full
shifts on $2n_k$ symbols, with $2n_k/(4k)\to1$. Combining this with
Obata's derivative estimate \cite[Section~2.2]{obatastandardmap}, we
obtain a universal $K_*>1$ such that
\begin{equation}\label{eq:standard-factor-estimates}
 h_{\mathrm{top}}(f_k)\ge\log k,
 \qquad
 0\le\lambda_1(f_k)<\log(4\pi k)
 \quad(k\ge K_*).
\end{equation}

Put $a_i=\lambda_1(f_{k_i})$. Area preservation gives
$\lambda_0(f_{k_i})=\lambda_2(f_{k_i})=0$. By
\cref{lem:sharp-product}, the largest exterior rate of $F$ is
$\lambda_m(F)=\sum_i a_i$. Every other degree omits at least one
$a_i$, and degree $m-1$ realizes each sum omitting one term. Hence
\begin{equation}\label{eq:standard-product-exterior}
 \lambda_{\mathrm{sec}}(F)
 =\max_j\sum_{i\ne j}a_i
 =\sum_i a_i-\min_i a_i.
\end{equation}
For $m=1$, this reads $\lambda_{\mathrm{sec}}(f_{k_1})=0$.
Entropy additivity and \eqref{eq:standard-factor-estimates} now give
\begin{equation}\label{eq:standard-product-gap}
 h_{\mathrm{top}}(F)-\lambda_{\mathrm{sec}}(F)
 \ge\min_i\log k_i-(m-1)\log(4\pi).
\end{equation}
Thus $k_0=\max\{K_*,e^\alpha(4\pi)^{m-1}\}$ has the required property,
and \cref{thm:main-intro} applies.
\end{proof}

\begin{proof}[Proof of \cref{cor:standard-product-uniqueness}]
Take $K_m$ larger than $k_0(1,m)$ and the threshold in
\cite[Theorem~A]{obatastandardmap}. Each factor then has a unique
Bernoulli maximal measure $\mu_i$, of entropy
$h_i=h_{\mathrm{top}}(f_{k_i})$. Their product is Bernoulli and maximal.
If $\nu$ is any maximal measure for $F$, with marginals $\nu_i$, then
\[
 \sum_i h_i=h_\nu(F)
 \le\sum_i h_{\nu_i}(f_{k_i})\le\sum_i h_i,
\]
so every $\nu_i=\mu_i$.

To identify the joining, choose finite-entropy Bernoulli generating
partitions $\cP_i$ for the factors and let
$\cP=\bigvee_i\pi_i^{-1}\cP_i$. This partition generates for $\nu$, so
the entropy formula for generating partitions \cite{Walters} gives
\[
 \sum_i h_i=h_\nu(F)
 \le\Ent_\nu(\cP)
 \le\sum_i\Ent_{\mu_i}(\cP_i)=\sum_i h_i.
\]
Equality makes the factor symbols independent at each time and the
joint symbol independent of its entire past. The joint process is
therefore Bernoulli with product symbol distribution, which implies
$\nu=\bigotimes_i\mu_i$. Finally,
\cref{thm:standard-map-products} and
\cref{thm:main-intro}\textnormal{(5)} give stable uniqueness and all
three continuity assertions.
\end{proof}

\small
\raggedright
\bibliographystyle{abbrv}
{\footnotesize\bibliography{library}}

\end{document}